\documentclass[preprint,10pt,3p]{elsarticle}
\usepackage{amssymb, amsmath, amsthm, amsfonts}
\usepackage{graphicx, subcaption}
\usepackage{multicol, multirow}
\usepackage{algorithm, algpseudocode}
\usepackage{color, hyperref}

\newtheorem{theorem}{\bf Theorem}[section]
\newtheorem{definition}{\bf Definition}[section]
\newtheorem{remark}{\bf Remark}[section]
\newtheorem{lemma}{\bf Lemma}[section]
\newtheorem{example}{\bf Example}[section]

\usepackage{mathrsfs}
\usepackage{cleveref}
\DeclareMathAlphabet{\mathbint}{U}{bbold}{m}{n}

\begin{document}

\begin{frontmatter}

\title{On Tensor Invertibility and Spectral Localization under the $M$-Product}
\author[aff1]{Aneesh Panchal}
\ead{aneeshp@iisc.ac.in}
\author[aff1]{Ratikanta Behera\corref{cor}}\cortext[cor]{Corresponding Author}
\affiliation[aff1]{organization={Department of Computational and Data Sciences, Indian Institute of Science}, city={Bangalore}, postcode={560012}, state={Karnataka}, country={IN}}
\ead{ratikanta@iisc.ac.in}

\begin{abstract}
The growth of multidimensional data demands rigorous multilinear algebraic frameworks. While $M$-product tensor transformations mimic matrix operations over a commutative ring of tubal scalars, the presence of zero-divisors has historically limited theoretical analysis. To overcome these limitations, in this paper, we introduce a novel tubal determinant to rigorously define tensor invertibility and resolve the zero-divisor effect (McCoy defect). In addition to this, we prove a tubal rank-nullity theorem using projective submodules and formalize tubal spectral analysis, yielding a novel Gershgorin hypersphere theorem. In addition to this, we find the Frobenius perturbation bounds for the pseudoinverse. Supporting the theoretical algebraic results, some numerical examples and visual illustrations are also provided in this paper. Finally, the theoretical results are applied experimentally for numerical illustrations. These results are essential and form the backbone for the theoretical analysis of tensor-based algorithms.
\end{abstract}

\begin{keyword}
Tubal tensor algebra, $M$-product, Tensor invertibility, Tubal rank-nullity theorem, Gershgorin hypersphere theorem.\\
\vspace{.35cm}
\end{keyword}
\end{frontmatter}


\section{Introduction}\label{sec:introduction}
The excessive growth of multidimensional data in modern engineering, including machine learning~\cite{ji2019survey}, signal processing~\cite{qi2017tensor}, and computer vision~\cite{panagakis2021tensor}, has driven the need for rigorous multilinear algebraic frameworks. While traditional tensor decompositions, such as the canonical polyadic~\cite{kruskal1977three, hong2020generalized}, tensor-train~\cite{oseledets2011tensor, gorodetsky2019function}, and Tucker decompositions~\cite{de2000multilinear, grasedyck2010hierarchical}, have proven highly effective for data compression and representation, they often suffer from a loss of the algebraic structures native to matrix theory. To bridge this gap, the $t$-product was introduced by Kilmer et al.~\cite{kilmer2011factorization, kilmer2013third}, allowing third-order tensors to be analyzed following the matrix structures. This framework was later generalized by Kernfeld et al.~\cite{kernfeld2015tensor} to the $M$-product, and subsequently extended to arbitrary $d$th-order tensors~\cite{qin2022low}, giving rise to a family of matrix-mimicking tensor-tensor products based on invertible linear transformations.

The $M$-product maps tensors into a module over a commutative ring of tubal scalars, $\mathbb{K}_M$. Under the $M$-product, classical matrix operations have been successfully generalized to tensors, leading to state-of-the-art computational results in robust tensor PCA~\cite{qiu2026robust, lin2024tensor}, tensor completion~\cite{wang2025robust, hou2025robust}, and optimization~\cite{feng2025learnable, wu2026guaranteed}. However, while the computational aspects of the $M$-product have been extensively explored, a rigorous, purely algebraic analysis remains underdeveloped, with only a few studies present in the literature by Avron and Mor~\cite{matrix-algebra-1} and Dunbar and Newman~\cite{matrix-algebra-2}. The primary mathematical challenge in tubal tensor algebra is that $\mathbb{K}_M$ does not form a field and consequently consists of zero-divisors. Without formalizing concepts like ideals, idempotents, and annihilators within the tubal domain, the theoretical analysis of tensor algorithms remains incomplete.

The primary mathematical hurdle in tubal tensor algebra is that $\mathbb{K}_M$ does not form a field, but it is a commutative ring equipped with zero-divisors. In the absence of formalized concepts such as ideals, idempotents, and annihilators within the tubal domain, fundamental notions like linear independence, determinants, and invertibility break down. Without these, proving the exactness, stability, or convergence of tensor-based algorithms remains heuristically driven rather than mathematically guaranteed. To address these gaps, this paper presents a comprehensive algebraic formalization of matrix-mimicking tubal tensor algebra. By leveraging commutative algebra~\cite{McCoy39, Hoffman-Kunze, ArtinRing27}, we rigorously define the necessary elements to analyze tensors over $\mathbb{K}_M$. The main contributions of the paper are as follows,
\begin{enumerate}
    \item Developed some of the foundational results for tensor algebra.
    \item Introduced determinant for tensors with novel theorems including invertibility of tensors and resolution of McCoy defect.
    \item Established tubal rank-nullity theorem for tubal tensor algebra. 
    \item Proposed spectral analysis for tensors with novel tubal spectral-divisor theorem and Gershgorin hypersphere theorem.
    \item Proved computational applicability of the tubal algebra through gradient descent, pseudoinverse bounds, and the Kaczmarz algorithm.
\end{enumerate}

For consistent notation, we have represented tensors as calligraphic `$\mathcal{A}$', tubal scalars as bold small alphabets `$\mathbf{a}$', tubal vectors as vectorized calligraphic `$\vec{\mathcal{A}}$', field scalars as small alphabets `$a$', field matrices as capital alphabets `$A$', fields as blackboard bold `$\mathbb{A}$', and sets, spaces and submodules as script `$\mathscr{A}$'.

The rest of the paper is organized as follows. \Cref{sec:notations-definitions} reviews the necessary prerequisites and existing literature in tubal algebra. \Cref{sec:matrix-mimic-algebra} develops the proposed tubal tensor algebra, including the tubal determinant and rank-nullity theorem. \Cref{sec:spectral-analysis} formalizes tubal spectral analysis and establishes the Gershgorin hypersphere theorem. \Cref{sec:applications} demonstrates the utility of these foundations through practical applications in tensor optimization and pseudo-inversion. Finally, \Cref{sec:conclusions} concludes the paper with future directions.

\section{Matrix-Mimicking Tensor Products}\label{sec:notations-definitions} 
Before getting deep into the definitions and properties of tubal tensor algebra, some of the basic definitions required for understanding the subsequent paper are provided in this subsection. We assume the basic results, including definitions of group, ring, and field, as prerequisites. We provide the necessary definitions required for understanding the subsequent sections. 

A \textit{left $\mathscr{R}$-module} is an abelian group $(\mathscr{M}, +)$ where, $(\mathscr{R}, +, \times)$ is a ring with identity $\mathbf{1}_{\mathscr{R}}$, equipped with a scalar multiplication operation $\ast:\mathscr{R} \times \mathscr{M} \to \mathscr{M}$, satisfying the following properties for all $\alpha, \beta \in \mathscr{R}$ and $\mu, \nu \in \mathscr{M}$,
\begin{enumerate}
    \item Distributivity over module addition: $\alpha\ast(\mu + \nu) = \alpha\ast \mu + \alpha\ast \nu$,
    \item Distributivity over ring addition: $(\alpha + \beta)\ast \mu = \alpha\ast \mu + \beta\ast \mu$,
    \item Associativity of scalar multiplication: $(\alpha\times\beta)\ast \mu = \alpha\times(\beta\ast \mu)$,
    \item Unital identity: $\mathbf{1}_{\mathscr{R}}\ast \mu = \mu$.
\end{enumerate}
A module is called a \textit{free module} if there exists a finite set of basis elements that spans the whole module. Let $\mathscr{M}$ be a free module over a ring $(\mathscr{R}, +,\times)$. A submodule $\mathscr{P} \subseteq \mathscr{M}$ is called a \textit{projective submodule} if it is a direct summand of a free module.

Let $\mathbb{F}$ be a field. An associative $\mathbb{F}$-algebra is a set $\mathscr{A}$ equipped with addition, scalar multiplication, and vector multiplication operations such that
\begin{enumerate}
    \item $(\mathscr{A}, +, \ast)$ forms a vector space over the field $\mathbb{F}$,
    \item $(\mathscr{A}, +, \times)$ forms a ring,
    \item For any $a \in \mathbb{F}$ and $\alpha, \beta \in \mathscr{A}$, we have, $a \ast (\alpha \times \beta) = (a \ast \alpha) \times \beta = \alpha \times (a \ast \beta)$.
\end{enumerate}

Now, we are equipped with the required definitions and properties. The next sections will emphasize an explanation of existing discussion and research regarding tubal tensor algebra, followed by proposed definitions, properties, and theorems to understand the tubal algebra better.

\subsection{$d$th-order Matrix-Mimicking Tensor Products}\label{subsec:dth-order-definitions}
The definition and some of the properties for tensors of order three already exist in the literature~\cite{matrix-algebra-1, matrix-algebra-2}, but it still lacks for higher-dimensional tensors. Hence, in this subsection we extend the definitions required to high-dimensional tensors.

\begin{definition}[Tubal scalars]\label{def:tubal-scalars-d}
    For tubal scalars $\mathbf{a},\mathbf{b}\in\mathbb{F}^{1\times 1\times n_3\times n_4\times\dots\times n_d}$ there exists matrix $M\in \mathrm{GL}_N(\mathbb{F})$ where, $\mathrm{GL}_N(\mathbb{F}) = \{M\in\mathscr{M}_{N\times N}(\mathbb{F})~|~\mathrm{det}(M)\neq 0\}$ is general linear group over field $\mathbb{F}$, $N = n_3n_4\dots n_d$, $\mathbb{F} = \mathbb{R}$ or, $\mathbb{C}$ and $M = M_{n_d}\otimes\dots\otimes M_{n_4}\otimes M_{n_3}$ where, $M_{n_i}\in \mathrm{GL}_{n_i}(\mathbb{F})$ then, $\mathbb{K}_M \equiv (\mathbb{F}^{N}, +, \circ_M)$ forms a commutative ring with identity where, $\circ_M$ is defined as,
    \begin{displaymath}
    \mathbf{a}\circ_M \mathbf{b} = \texttt{devec}(M^{-1}(M\texttt{vec}(\mathbf{a})\odot M\texttt{vec}(\mathbf{b})))\in\mathbb{F}^{1\times 1\times n_3\times n_4\times\dots\times n_d},
    \end{displaymath}
    where, $\texttt{vec}$ and $\texttt{devec}$ are linear operators assuming column-major vectorization, and $\odot$ is element-wise multiplication. The equation is equivalent to $\mathbf{a}\circ_M \mathbf{b} = \mathfrak{M}^{-1}(\mathfrak{M}(\mathbf{a})\odot\mathfrak{M}(\mathbf{b}))$ where, $\mathfrak{M}(\mathbf{a}) = \mathbf{a}\times_3M_{n_3}\times_4M_{n_4}\dots\times_dM_{n_d}$ and $\mathfrak{M}^{-1}(\mathbf{a}) = \mathbf{a}\times_dM_{n_d}^{-1}\dots\times_4M^{-1}_{n_4}\times_3M^{-1}_{n_3}$.
\end{definition}

Here, we represent the additive identity of $\mathbb{K}_M$ as $\mathbint{0}$ and multiplicative identity as $\mathbint{1}$. It is clear from \Cref{def:tubal-scalars-d} that $\mathbb{K}_M$ is not a division ring. However, the pseudoinverse for every element exists. Hence, it is defined as follows.

\begin{definition}[{\cite[Lemma 6 and 7]{matrix-algebra-1}}]
    The pseudoinverse for a tubal scalar $\mathbf{a}\in\mathbb{K}_M$, is defined as $\mathbf{a}^\dagger = \mathfrak{M}^{-1}((\mathfrak{M}(\mathbf{a}))^\dagger)$ where, for each $x$ in $\mathfrak{M}(\mathbf{a})$, the corresponding $(\mathfrak{M}(\mathbf{a}))^\dagger$ is $x^{-1}$ if $x \neq 0$, and $0$ otherwise.
\end{definition}

\begin{remark}[{\cite[Definition 8]{matrix-algebra-1}}]\label{rem:Von-Neumann-Regular}
    It is important to note that $\mathbb{K}_M$ is Von-Neumann regular that is, for every $\mathbf{a}\in\mathbb{K}_M$, there exists an element $\mathbf{x}\in\mathbb{K}_M$ such that $\mathbf{a} = \mathbf{a}\circ_M\mathbf{x}\circ_M\mathbf{a}$.
\end{remark}

\begin{remark}\label{rem:semisimple-Artinian-ring}
    It is clear that $\mathbb{K}_M$ is isomorphic to $(\mathbb{F}^N,+,\odot)$, that is, the direct product of fields. Hence, using Wedderburn-Artin Theorem~\cite{Wedderburn08, ArtinRing27} and \Cref{rem:Von-Neumann-Regular}, $\mathbb{K}_M$ is semisimple Artinian ring.
\end{remark}

While being mathematically consistent, it is difficult to interpret the results with direct elements of $\mathbb{K}_M$. Hence, we will define a mapping which aids in the visualisation of results in the following sections.

\begin{definition}[{\cite[Proposition 2.1]{matrix-algebra-2}}]\label{def:transfomed-domain-mapping}
    According to \Cref{rem:semisimple-Artinian-ring}, there exists a linear mapping $\Phi:\mathbb{K}_M\to \mathbb{F}^N$ such that $\Phi(\mathbf{a}) = M\texttt{vec}(\mathbf{a})$ for $\mathbf{a}\in\mathbb{K}_M$. Also, define the inverse linear mapping $\Phi^{-1}:\mathbb{F}^N\to \mathbb{K}_M$ as $\Phi^{-1}(b) = \texttt{devec}(M^{-1}(b))$ for $b\in\mathbb{F}^N$.
\end{definition}

\begin{lemma}[{\cite[Theorem 4.1]{kernfeld2015tensor}}]\label{lem:free-module}
    $\mathbb{K}_M^k$ is a free module over $\mathbb{K}_M$ of rank $k$.
\end{lemma}

As a result of \Cref{lem:free-module}, multiplication of tubal scalars with tubal vectors and tensors is defined. For example, for tubal scalar $\mathbf{a}\in\mathbb{K}_M$ and tubal vector $\vec{\mathcal{Q}} = [\mathbf{q}_1, \mathbf{q}_2, \dots, \mathbf{q}_k]^\top\in\mathbb{K}_M^k$, the product $\mathbf{a}\circ_M \vec{\mathcal{Q}} = [\mathbf{a}\circ_M\mathbf{q}_1, \mathbf{a}\circ_M\mathbf{q}_2, \dots, \mathbf{a}\circ_M\mathbf{q}_k]^\top$ is well defined. The same will follow for tensors (also called tubal matrices) as well.

\begin{definition}[Matrix-mimicking tensor product of order $d$]\label{def:star-m-d}
    For tensors $\mathcal{A}\in\mathbb{F}^{n_1\times \ell\times n_3\times n_4\times\dots\times n_d} \equiv \mathbb{K}_M^{n_1\times \ell}$ and $\mathcal{B}\in\mathbb{F}^{\ell\times n_2\times n_3\times n_4\times\dots\times n_d}\equiv \mathbb{K}_M^{\ell\times n_2}$ there exists $M\in \mathrm{GL}_N(\mathbb{F})$ where, $N = n_3n_4\dots n_d$, $\mathbb{F} = \mathbb{R}$ or, $\mathbb{C}$ and $M = M_{n_d}\otimes\dots\otimes M_{n_4}\otimes M_{n_3}$ then, the tensor product $\mathcal{A}\star_M\mathcal{B}\in\mathbb{F}^{n_1\times n_2\times n_3\times n_4\times\dots\times n_d}$ is defined as,
    \begin{displaymath}
    (\mathcal{A}\star_M\mathcal{B})^{[i,j]} = \sum_{r = 1}^{\ell}\mathcal{A}^{[i,r]}\circ_M\mathcal{B}^{[r,j]}\in\mathbb{F}^{1\times 1\times n_3\times n_4\times\dots\times n_d},
    \end{displaymath}
    for $i\in[n_1]$ and $j\in[n_2]$ where, $\circ_M$ is defined in \Cref{def:tubal-scalars-d}.
\end{definition}

\begin{lemma}\label{lem:f-algebra}
    For the field $(\mathbb{F},+,\ast)$, the structure $(\mathbb{K}_M,+,\ast,\circ_M)$ forms an associative $\mathbb{F}$-algebra.
\end{lemma}
\begin{proof}
    Using \Cref{def:tubal-scalars-d} and the fact that $\mathfrak{M}$ and $\mathfrak{M}^{-1}$ are linear operators, we can easily prove that, $a \ast (\mathbf{a}\circ_M\mathbf{b}) = (a \ast \mathbf{a}) \circ_M \mathbf{b} = \mathbf{a}\circ_M (a \ast \mathbf{b})$ where, $a\in\mathbb{F}$ and $\mathbf{a}, \mathbf{b}\in\mathbb{K}_M$. Hence, the structure $(\mathbb{K}_M,+,\ast,\circ_M)$ forms an associative $\mathbb{F}$-algebra.
\end{proof}

The correspondence of \Cref{lem:f-algebra} with real-life applications can be seen in the literature, including gradient descent, principal component analysis, and many more, where scalars are multiplied with tensors elementwise directly.

From here onwards, we will consider only $\mathbb{F} = \mathbb{R}$ and leave the rest of the analysis, including complex and non-commutative quaternions, for future work.

\subsection{Analysis of Existing Tensor-Tensor Products}\label{sec:tensor-tensor-products}
In the current literature, there exist various types of tensor-tensor products, $t$-product~\cite{kilmer2011factorization, kilmer2013third}, $M$-product~\cite{kernfeld2015tensor}, $c$-product~\cite{kernfeld2015tensor}, $\phi$-product~\cite{song2020robust}, $f$-product~\cite{framelet-product}, $w$-product~\cite{panchal2026second}, and many more with their applications \cite{sahoo2024computations, KarmakaBehera2026, MiaoWei2020}. Despite having different analyses and usage, all of the tensor products follow the $M$-product with different transformation matrices $M$. For example, the $t$-product utilizes the DFT matrix, $c$-product utilizes the DCT matrix.

\begin{definition}[$M$-product~{\cite[Lemma 4.1]{kernfeld2015tensor}}]\label{def:mproduct}
    Let tensors $\mathcal{A}\in\mathbb{R}^{n_1\times \ell\times n_3\times\dots\times n_d}$ and $\mathcal{B}\in\mathbb{R}^{\ell\times n_2\times n_3\times\dots\times n_d}$, then \textbf{$\mathbf{{M}}$-product} is defined as,
    \begin{displaymath}
    \mathcal{A}\star_M\mathcal{B} = \mathfrak{M}^{-1}(\mathfrak{M}(\mathcal{A})\Delta\mathfrak{M}(\mathcal{B}))\in\mathbb{R}^{n_1\times n_2\times n_3\times\dots\times n_d},
    \end{displaymath}
    where $\Delta$ is the face-wise matrix product, $\mathfrak{M}(\mathcal{A}) = \mathcal{A}\times_{3}M_{n_3}\times_{4}M_{n_4}\dots\times_{d}M_{n_d}$, and $\mathfrak{M}^{-1}(\mathcal{A}) = \mathcal{A}\times_{d}M_{n_d}^{-1}\dots\times_{4}M_{n_4}^{-1}\times_{3}M_{n_3}^{-1}$, where, $M_{n_i}\in\mathbb{R}^{n_i\times n_i}$ is an invertible matrix.
\end{definition}

One of the reasons for the $M$-product to be so popular is due to the tubal algebra, which is matrix-mimicking. Now, some of the properties of $M$-products are well established and used in various applications. The crucial theorem in the literature for matrix-mimicking tensors is singular value decomposition (SVD), which is given by

\begin{theorem}[Singular value decomposition~{\cite[Theorem 5.1]{kernfeld2015tensor}}]\label{thrm:orig-svd}
    Let $\mathcal{A}\in\mathbb{R}^{n_1\times n_2\times n_3\times\dots\times n_d}$ be a $d$th-order tensor. Then, there exist tensors $\mathcal{U}\in\mathbb{R}^{n_1\times \ell\times n_3\times\dots\times n_d}$, $\Sigma\in\mathbb{R}^{\ell\times \ell\times n_3\times\dots\times n_d}$, and $\mathcal{V}\in\mathbb{R}^{n_2\times \ell\times n_3\times\dots\times n_d}$ such that $\mathcal{A}$ admits a low-rank approximation,
    \begin{displaymath}
        \mathcal{A}^{(r)} \approx \mathcal{U}\star_M \Sigma\star_M \mathcal{V}^\top = \sum_{i=1}^r \vec{\mathcal{U}_i} \star_M \vec{\sigma_i}\star_M \vec{\mathcal{V}_i}^\top ,
    \end{displaymath}
    where $\mathcal{U}$ and $\mathcal{V}$ are orthogonal tensors. The tensors $\vec{\mathcal{U}_i}\in\mathbb{R}^{n_1\times 1\times n_3\times\dots\times n_d}, \vec{\sigma_i}\in\mathbb{R}^{1\times 1\times n_3\times\dots\times n_d}, \vec{\mathcal{V}_i}\in\mathbb{R}^{n_2\times 1\times n_3\times\dots\times n_d}$ are the tubal vectors, for the decomposition to be represented as sum of single-rank tensors. Also, for any rank $r\leq\min\{n_1,n_2\}$, $\mathcal{A}^{(r)}$ satisfies Eckart-Young theorem as follows,
    \begin{displaymath}
    \|\mathcal{A} - \mathcal{A}^{(r)}\|_F\leq \frac{1}{\sigma_{min}(M)}\sqrt{\sum_{i = r+1}^{\min\{n_1,n_2\}} \|\vec{\sigma_i} \|_F^2 }.
    \end{displaymath}
\end{theorem}

\begin{definition}[Moore-Penrose Inverse]
    The pseudoinverse for a tensor $\mathcal{A}\in\mathbb{K}_M^{m\times n}$, is defined as $\mathcal{A}^\dagger = \sum_{i=1}^r \vec{\mathcal{V}_i} \star_M \vec{\sigma_i}^\dagger\star_M \vec{\mathcal{U}_i}^\top$, where $\mathcal{U}$ and $\mathcal{V}$ are orthogonal tensors. The pseudoinverse of a tensor $\mathcal{A}$ satisfies the Moore-Penrose conditions given by,
    \begin{displaymath}
        \mathcal{A}\star_M\mathcal{A}^\dagger\star_M\mathcal{A} = \mathcal{A},~ \mathcal{A}^\dagger\star_M\mathcal{A}\star_M\mathcal{A}^\dagger = \mathcal{A}^\dagger,~ (\mathcal{A}\star_M\mathcal{A}^\dagger)^\top = \mathcal{A}\star_M\mathcal{A}^\dagger,~ (\mathcal{A}^\dagger\star_M\mathcal{A})^\top = \mathcal{A}^\dagger\star_M\mathcal{A}
    \end{displaymath}
    Hence, it is also known as the Moore-Penrose inverse and is unique for a tensor $\mathcal{A}$.
\end{definition}

Kilmer et al.~\cite{kilmer2013third} gave some of the important definitions required for matrix-mimicking algebra in $t$-product settings, which are applicable directly to $M$-product settings as well. Some of the definitions, including rank, range space, and nullity, are provided here.

\begin{definition}[Multirank, multinullity and tubal rank~{\cite[Definition 4.5]{kilmer2013third}}]\label{def:tubal-rank}
    Let the tensor $\mathcal{A}\in\mathbb{R}^{n_1\times n_2\times p_3\times \dots\times p_d}$ decomposes into $\prod_{j=3}^d p_j$ matrices using $M$-transformation ($i$th matrix is denoted by $A_i\in\mathbb{R}^{n_1\times n_2}$). Then, multirank is defined as $\mathfrak{r} = \mathfrak{r}(\mathcal{A})\in\mathbb{N}^{\prod_j p_j}$, where, $\mathfrak{r}[i]\leq\min\{n_1,n_2\}$ is the rank of $A_i$. On similar grounds, multinullity is defined as $\mathfrak{n} = \mathfrak{n}(\mathcal{A})\in\mathbb{N}^{\prod_j p_j}$, where, $\mathfrak{n}[i] = n_2-\mathfrak{r}[i]$ is the nullity of $A_i$. Now, tubal rank can be defined using multirank as $\mathfrak{r}_{\mathfrak{t}} = \max({\mathfrak{r}(\mathcal{A})})$, which is equivalent to the number of non-zero singular tubes in the SVD decomposition. 
\end{definition}

Assuming for any vectorized tensor $\vec{\mathcal{X}}\in \mathbb{R}^{n_2\times 1\times n_3\times\dots\times n_d}$ we have singular value decomposition decomposition of the form,
\begin{equation}\label{eq:Ax-svd-expansion}
    \mathcal{A}\star_M \vec{\mathcal{X}} = \sum_{i=1}^{\min\{n_1, n_2\}} \vec{\mathcal{U}}_i \star_M \vec{\sigma}_i\star_M \underbrace{\vec{\mathcal{V}}_i^\top \star_M \vec{\mathcal{X}}}_{\vec{d}_i\in\mathbb{R}^{1\times 1\times n_3\times\dots\times n_d}}.
\end{equation}

The range and null space for the tensor space under tubal algebra using SVD defined in \Cref{thrm:orig-svd} are given as follows,

\begin{definition}[Range and null space~{\cite[Theorem 4.6]{kilmer2013third}}]\label{def:svd-range-null-space}
    Assuming the initial $k_1$ singular tubes are invertible, the next $k_2$ singular tubes are non-zero but not invertible, and the rest of the $\left(\min\{n_1,n_2\} - (k_1+k_2)\right)$ are zero tubes. Then, range space $\mathscr{R}(\mathcal{A})$ and null space $\mathscr{N}(\mathcal{A})$ are given by,
    \begin{displaymath}
        \mathscr{R}(\mathcal{A}) = \left\{\vec{\mathcal{U}}_1\star_M \vec{c}_1 + \dots + \vec{\mathcal{U}}_{k_1+k_2}\star_M \vec{c}_{k_1+k_2}~|~\vec{c}_i = \vec{\sigma}_i\star_M\vec{d}_i,~ k_1<i\leq k_1+k_2 \right\},
    \end{displaymath}
    \begin{displaymath}
        \mathscr{N}(\mathcal{A}) = \left\{\vec{\mathcal{V}}_{k_1+1}\star_M \vec{c}_{k_1+1} + \dots + \vec{\mathcal{V}}_{n_2}\star_M \vec{c}_{n_2}~|~ \vec{\sigma}_i\star_M\vec{c}_i = \vec{0},~ k_1<i\leq k_1+k_2 \right\}.
    \end{displaymath}
\end{definition}

\section{Matrix-Mimicking Algebra}\label{sec:matrix-mimic-algebra}
Following the definitions of the tensor-tensor product of order $d$ provided in \Cref{subsec:dth-order-definitions}, in this section, we will define and prove some of the important properties required for the analysis of the matrix-mimicking tensor systems, including determinant, rank, dependence, and invertibility of tensors.

In a standard totally ordered field, strict inequality is often defined simply as $\mathbf{a} \ge \mathbf{b}$ and $\mathbf{a} \neq \mathbf{b}$. But in the case of $\mathbb{K}_M$, this definition is insufficient due to the presence of zero-divisors. Hence, a definition for tubal inequality is required.

\begin{definition}[Strict tubal inequality]\label{def:tubal-inequality}
Let $\mathbf{a}, \mathbf{b} \in \mathbb{K}_M$. Then $\mathbf{a}$ is strictly greater than $\mathbf{b}$, denoted $\mathbf{a} \succ_{\mathbb{K}_M} \mathbf{b}$, if the following conditions hold,
\begin{enumerate}
    \item $\mathbf{a} \succeq_{\mathbb{K}_M} \mathbf{b}$, that is, there exists some $\mathbf{x} \in \mathbb{K}_M$ such that $\mathbf{a} - \mathbf{b} = \mathbf{x} \circ_M \mathbf{x}$.
    \item $(\mathbf{a} - \mathbf{b})$ is a unit in $\mathbb{K}_M$, that is, it is invertible, and its annihilator ideal is strictly trivial, $\mathrm{Ann}(\mathbf{a} - \mathbf{b}) = \{\mathbint{0}\}$, where, $\mathrm{Ann}(\mathbf{x}) = \{\kappa\in\mathbb{K}_M|\kappa\circ_M \mathbf{x} = \mathbint{0}\}$.
\end{enumerate}
\end{definition}

In the field $\mathbb{R}$, the fundamental elements are 1 and 0. Now, we define the same building blocks for a tensor, which we call fundamental idempotent tubal scalars.  

\begin{definition}[Fundamental idempotent tubal scalars]\label{def:funda-idem-tubal-scalar}
    $\mathbb{K}_M$ contains a set of $N$ fundamental idempotents $\{\mathbf{e}_1, \dots, \mathbf{e}_N\}$, satisfying,
    \begin{enumerate}
        \item $\mathbf{e}_i \circ_M \mathbf{e}_i = \mathbf{e}_i$,
        \item $\mathbf{e}_i \circ_M \mathbf{e}_j = \mathbint{0}$ for $i \neq j$,
        \item $\sum_{i=1}^N \mathbf{e}_i = \mathbint{1}$.
    \end{enumerate}
\end{definition}

Using \Cref{def:transfomed-domain-mapping}, the natural fundamental idempotent scalar set is given by, $\{\mathfrak{e}_1, \mathfrak{e}_2, \dots, \mathfrak{e}_N\} = \{\Phi^{-1}([1,0,\dots,0]), \Phi^{-1}([0,1,\dots,0]), \dots, \Phi^{-1}([0,0,\dots,1])\}$.

\begin{definition}[Natural basis]\label{def:tubal-natural-basis}
    Since $\mathbb{K}_M^{k}$ is a free module over $\mathbb{K}_M$, hence, tubal vectors $\{\vec{\mathcal{I}}_1, \vec{\mathcal{I}}_2, \dots, \vec{\mathcal{I}}_k\} = \{[\mathbint{1}, \mathbint{0}, \dots, \mathbint{0}]^\top, [\mathbint{0}, \mathbint{1}, \dots, \mathbint{0}]^\top, \dots, [\mathbint{0}, \mathbint{0}, \dots, \mathbint{1}]^\top\}$ forms the natural basis for $\mathbb{K}_M^{k}$ satisfying,
    \begin{enumerate}
        \item $\vec{\mathcal{I}}_i^\top\star_M\vec{\mathcal{I}}_i = \mathbint{1}$,
        \item $\vec{\mathcal{I}}_i^\top\star_M\vec{\mathcal{I}}_j = \mathbint{0}$ for $i\neq j$,
        \item $\mathfrak{D}(\vec{\mathcal{I}}_1, \vec{\mathcal{I}}_2, \dots, \vec{\mathcal{I}}_k) = \mathbint{1}$.
    \end{enumerate}
\end{definition}

\begin{figure}[htbp]
    \centering
    \includegraphics[width=0.8\linewidth]{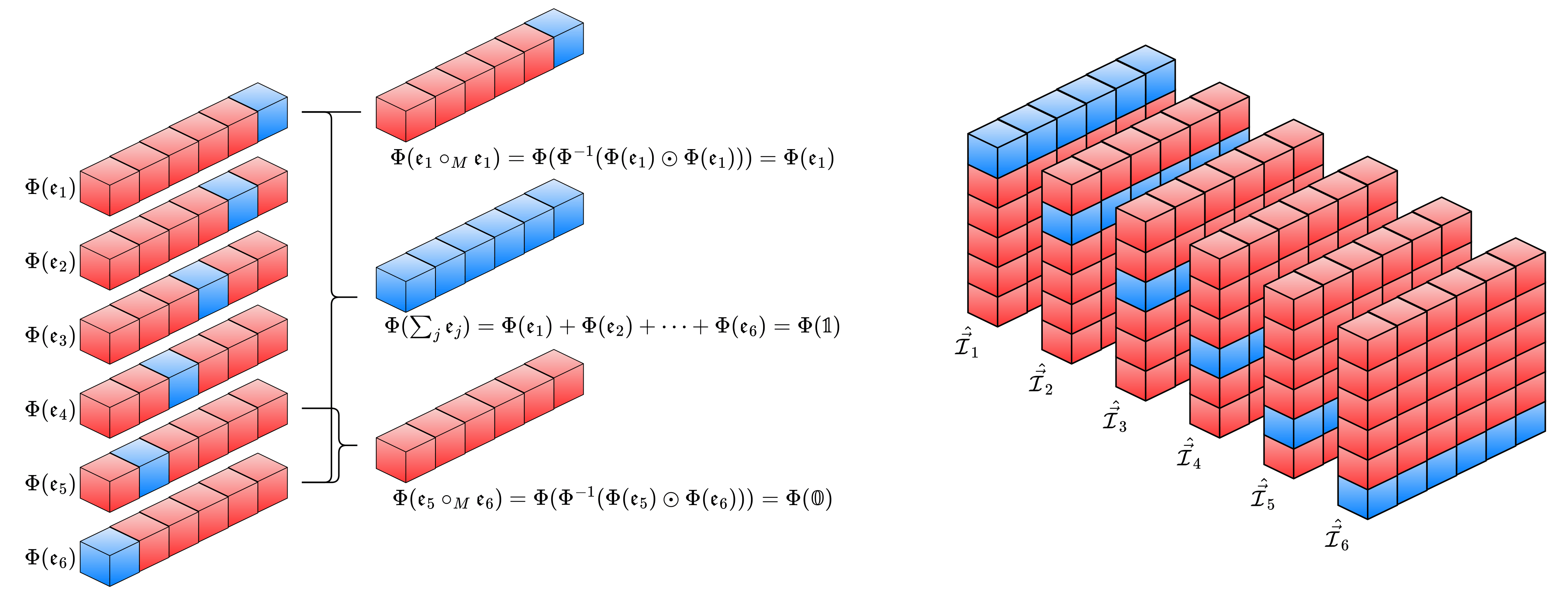}
    \caption{Visualization of fundamental idempotent scalars (left, \Cref{def:funda-idem-tubal-scalar}) and natural basis (right, \Cref{def:tubal-natural-basis}) for tubal tensor algebra.}
    \label{fig:Idempotents}
\end{figure}

The fundamental idempotent tubal scalars and natural basis for tubal tensor algebra can be visualized through \Cref{fig:Idempotents}.

\subsection{Determinant and Invertibility of Tensors}
The determinant is one of the basic metrics used to define the invertibility of a matrix of field elements. However, there is no such metric available for tensors. Hence, we start with basic definitions required for finding the determinant of a tensor.

\begin{definition}[$k$-tubal-linear function]
    Let $k\in\mathbb{N}$, and let function $\mathfrak{D}$ maps tensor $\mathcal{A}\in\mathbb{K}_M^{k\times k}$ to a tubal scalar, $\mathfrak{D}(\mathcal{A})$ in $\mathbb{K}_M$. Then we say, $\mathfrak{D}$ is $k$-tubal-linear if for every $i$, $1\leq i\leq k$, $\mathfrak{D}$ is linear function of $i$th tubal row when other $(k-1)$ tubal rows are kept fixed. 
\end{definition}

\begin{definition}[Alternating function]
    Let $\mathfrak{D}$ be a $k$-tubal-linear function. Then we can say, $\mathfrak{D}$ is an alternating function if,
    \begin{enumerate}
        \item $\mathfrak{D}(\mathcal{A}) = \mathbint{0}$, whenever two tubal-rows are equal.
        \item If $\mathcal{A}'$ is obtained by interchanging two tubal rows of $\mathcal{A}$, then $\mathfrak{D}(\mathcal{A}') = -\mathfrak{D}(\mathcal{A})$.
    \end{enumerate}
    where $\mathbint{0}$ is the additive identity in $\mathbb{K}_M$.
\end{definition}

Now, with the basic setup, we can define the determinant of a tensor using a matrix-mimicking framework. It is important to note that the determinant of a tensor is a tubal scalar which belongs to the ring $\mathbb{K}_M$ and not a field.

\begin{definition}[Determinant]\label{def:determinant}
    Let $k\in\mathbb{N}$, and let function $\mathfrak{D}$ maps tensor $\mathcal{A}\in\mathbb{K}_M^{k\times k}$ to a tubal scalar $\mathfrak{D}(\mathcal{A})$ in $\mathbb{K}_M$. Then $\mathfrak{D}$ is a determinant function if $\mathfrak{D}$ is $k$-tubal-linear, alternating, and $\mathfrak{D}(\mathcal{I}) = \mathbint{1}$, where $\mathbint{1}$ is the multiplicative identity in $\mathbb{K}_M$.
\end{definition}

According to \Cref{def:transfomed-domain-mapping}, computationally, determinant of an $d$th order tensor can be calculated as directly $\mathfrak{D}(\mathcal{A}) = \mathfrak{M}^{-1}(\mathrm{det}_\Delta(\mathfrak{M}(\mathcal{A}))$ which is similar to slice-wise determinant. The illustration showing the same is provided in \Cref{example:determinant}. It is important to note that the determinant can be a zero divisor of the ring $\mathbb{K}_M$ as well.

Now, we will define some properties and results for the determinant, including uniqueness, transpose, and others, showing its matrix-mimicking nature.

\begin{lemma}[Determinant uniqueness]\label{lem:det-pre}
    For tensor $\mathcal{A}\in\mathbb{K}_M^{k\times k}$, determinant function $\mathfrak{D}: \mathbb{K}_M^{k\times k}\to \mathbb{K}_M$, and arbitrary alternating $k$-tubal-linear function $\mathfrak{F}: \mathbb{K}_M^{k\times k}\to \mathbb{K}_M$, we have,
    \begin{displaymath}
    \mathfrak{F}(\mathcal{A}) = \mathfrak{D}(\mathcal{A})\circ_M\mathfrak{F}(\mathcal{I}).
    \end{displaymath}
\end{lemma}
\begin{proof}
Assuming a determinant function $\mathfrak{D}:\mathbb{K}_M^{k\times k}\to\mathbb{K}_M$ and an arbitrary alternating $k$-tubal-linear function $\mathfrak{F}: \mathbb{K}_M^{k\times k}\to \mathbb{K}_M$. Let tensor $\mathcal{A}\in\mathbb{K}_M^{k\times k}$ with tubal rows $\vec{\mathfrak{a}}_1, \vec{\mathfrak{a}}_2, \dots, \vec{\mathfrak{a}}_k$ and tubal rows of identity tensor $\mathcal{I}\in\mathbb{K}_M^{k\times k}$ be $\vec{\mathfrak{e}}_1, \vec{\mathfrak{e}}_2, \dots, \vec{\mathfrak{e}}_k$. Then it is clear that, $\vec{\mathfrak{a}}_i = \sum_{j=1}^k \mathcal{A}^{[i,j]}\star_M\vec{\mathfrak{e}}_j, 1\leq i\leq k$. Hence, the function $\mathfrak{F}$ can be written as $\mathfrak{F}(\mathcal{A}) = \sum_j\mathcal{A}^{[1,j]}\circ_M\mathfrak{F}(\vec{\mathfrak{e}}_j, \vec{\mathfrak{a}}_2, \dots, \vec{\mathfrak{a}}_k)$. Similarly, after generalization, we get,
\begin{align*}
    \mathfrak{F}(\mathcal{A}) &= \sum_{j_1, j_2,\dots, j_k}\mathcal{A}^{[1,j_1]}\circ_M\mathcal{A}^{[2,j_2]}\circ_M\dots\circ_M\mathcal{A}^{[k,j_k]}\circ_M\mathfrak{F}(\vec{\mathfrak{e}}_{j_1}, \vec{\mathfrak{e}}_{j_2}, \dots, \vec{\mathfrak{e}}_{j_k})\\
    &= \sum_\sigma \mathcal{A}^{[1,\sigma 1]}\circ_M\mathcal{A}^{[2,\sigma 2]}\circ_M\dots\circ_M\mathcal{A}^{[k,\sigma k]}\circ_M\mathfrak{F}(\vec{\mathfrak{e}}_{\sigma 1}, \vec{\mathfrak{e}}_{\sigma 2}, \dots, \vec{\mathfrak{e}}_{\sigma k})\\
    &= \sum_\sigma \mathrm{sgn}(\sigma) \circ_M \mathcal{A}^{[1,\sigma 1]}\circ_M\mathcal{A}^{[2,\sigma 2]}\circ_M\dots\circ_M\mathcal{A}^{[k,\sigma k]}\circ_M \mathfrak{F}(\mathcal{I}),
\end{align*}

where, $\mathrm{sgn}(\sigma) = +\mathbint{1}$ if $\sigma\in2\mathbb{N}$, else $-\mathbint{1}$. Also, we have assumed $\sigma$ to be a permutation of degree $k$, and the sequence $(\sigma 1,\sigma 2,\dots, \sigma k)$ can be obtained through a finite permutation of the sequence $(1,2,\dots,k)$. As, function $\mathfrak{F}$ is an alternating $k$-tubal-linear function, hence, $\mathfrak{F}(\vec{\mathfrak{e}}_{\sigma 1}, \vec{\mathfrak{e}}_{\sigma 2}, \dots, \vec{\mathfrak{e}}_{\sigma k}) = \mathrm{sgn}(\sigma)\circ_M \mathfrak{F}(\vec{\mathfrak{e}}_{1}, \vec{\mathfrak{e}}_{2}, \dots, \vec{\mathfrak{e}}_{k}) = \mathrm{sgn}(\sigma) \circ_M\mathfrak{F}(\mathcal{I})$. Hence, we can write,
\begin{displaymath}
    \mathfrak{F}(\mathcal{A}) = \mathfrak{D}(\mathcal{A}) \circ_M \mathfrak{F}(\mathcal{I}),
\end{displaymath}

where, $\mathfrak{D}(\mathcal{A}) = \sum_\sigma \mathrm{sgn}(\sigma) \circ_M \mathcal{A}^{[1,\sigma 1]}\circ_M\mathcal{A}^{[2,\sigma 2]}\circ_M\dots\circ_M\mathcal{A}^{[k,\sigma k]}$. Additionally, if we assume the property $\mathfrak{F}(\mathcal{I}) = \mathbint{1}$, then $\mathfrak{F}(\mathcal{A}) = \mathfrak{D}(\mathcal{A})$ for all tensors $\mathcal{A}$. Hence, a unique determinant function exists.
\end{proof}

\begin{lemma}\label{lem:det(ab)}
    For tensors $\mathcal{A},\mathcal{B}\in\mathbb{K}_M^{k\times k}$ and determinant function $\mathfrak{D}: \mathbb{K}_M^{k\times k}\to \mathbb{K}_M$, we have,
    \begin{displaymath}
    \mathfrak{D}(\mathcal{A}\star_M\mathcal{B}) = \mathfrak{D}(\mathcal{A})\circ_M\mathfrak{D}(\mathcal{B}).
    \end{displaymath}
\end{lemma}
\begin{proof}
    For tensors $\mathcal{A},\mathcal{B}\in\mathbb{K}_M^{k\times k}$ and determinant function $\mathfrak{D}: \mathbb{K}_M^{k\times k}\to \mathbb{K}_M$, assume an arbitrary function $\mathfrak{F}:\mathbb{K}_M^{k\times k}\to\mathbb{K}_M$ such that, $\mathfrak{F}(\mathcal{A}) = \mathfrak{D}(\mathcal{A}\star_M\mathcal{B})$. Let tensor $\mathcal{A}\in\mathbb{K}_M^{k\times k}$ have tubal rows $\vec{\mathfrak{a}}_1, \vec{\mathfrak{a}}_2, \dots, \vec{\mathfrak{a}}_k$, then $\mathfrak{F}(\mathcal{A}) = \mathfrak{F}(\vec{\mathfrak{a}}_1, \vec{\mathfrak{a}}_2, \dots, \vec{\mathfrak{a}}_k) = \mathfrak{D}(\vec{\mathfrak{a}}_1\star_M\mathcal{B}, \vec{\mathfrak{a}}_2\star_M\mathcal{B}, \dots, \vec{\mathfrak{a}}_k\star_M\mathcal{B})$. It is easy to show that $\mathfrak{F}$ is an alternating $k$-tubal-linear function because $\mathfrak{D}$ is a determinant function. Hence, applying \Cref{lem:det-pre}, we get, $\mathfrak{D}(\mathcal{A}\star_M\mathcal{B}) = \mathfrak{D}(\mathcal{A})\circ_M\mathfrak{D}(\mathcal{B})$.
\end{proof}

\begin{lemma}
    For tensor $\mathcal{A}\in\mathbb{K}_M^{k\times k}$ and determinant function $\mathfrak{D}: \mathbb{K}_M^{k\times k}\to \mathbb{K}_M$, we have, $\mathfrak{D}(\mathcal{A}^\top) = \mathfrak{D}(\mathcal{A})$.
\end{lemma}
\begin{proof}
    For tensor $\mathcal{A}\in\mathbb{K}_M^{k\times k}$ and determinant function $\mathfrak{D}: \mathbb{K}_M^{k\times k}\to \mathbb{K}_M$, we have,
    \begin{align*}
        &\mathfrak{D}(\mathcal{A}) = \sum_\sigma \mathrm{sgn}(\sigma) \circ_M \mathcal{A}^{[1,\sigma 1]}\circ_M\mathcal{A}^{[2,\sigma 2]}\circ_M\dots\circ_M\mathcal{A}^{[k,\sigma k]},\\
        &\mathfrak{D}(\mathcal{A}^\top) = \sum_\sigma \mathrm{sgn}(\sigma) \circ_M \mathcal{A}^{[\sigma 1,1]}\circ_M\mathcal{A}^{[\sigma 2,2]}\circ_M\dots\circ_M\mathcal{A}^{[\sigma k,k]}.
    \end{align*}
    As $\sigma$ is a permutation of degree $k$, there exists a permutation $\sigma^{-1}$ such that $\sigma\sigma^{-1}$ is identity permutation. Hence, $\mathrm{sgn}(\sigma)~\mathrm{sgn}(\sigma^{-1}) = 1$, and,
    \begin{displaymath}
        \mathcal{A}^{[\sigma 1,1]}\circ_M\mathcal{A}^{[\sigma 2,2]}\circ_M\dots\circ_M\mathcal{A}^{[\sigma k,k]} = \mathcal{A}^{[1,\sigma^{-1} 1]}\circ_M\mathcal{A}^{[2,\sigma^{-1} 2]}\circ_M\dots\circ_M\mathcal{A}^{[k,\sigma^{-1} k]}.
    \end{displaymath}
    As $\sigma^{-1}$ permutes over all combinations as well. Hence, we can write,
    \begin{displaymath}
        \mathfrak{D}(\mathcal{A}^\top) = \sum_\sigma \mathrm{sgn}(\sigma^{-1})\mathcal{A}^{[1,\sigma^{-1} 1]}\circ_M\mathcal{A}^{[2,\sigma^{-1} 2]}\circ_M\dots\circ_M\mathcal{A}^{[k,\sigma^{-1} k]} = \mathfrak{D}(\mathcal{A}).
    \end{displaymath}
\end{proof}

\begin{lemma}\label{lem:adjoint-determinant-rel}
    For tensor $\mathcal{A}\in\mathbb{K}_M^{k\times k}$ and determinant function $\mathfrak{D}: \mathbb{K}_M^{k\times k}\to \mathbb{K}_M$, we have $\texttt{adj}(\mathcal{A})\star_M\mathcal{A} = \mathfrak{D}(\mathcal{A})\star_M\mathcal{I},$ and, $\mathcal{A}\star_M\texttt{adj}(\mathcal{A}) = \mathfrak{D}(\mathcal{A})\star_M\mathcal{I}$.
\end{lemma}
\begin{proof}
    It is easy to show that for the tensor $\mathcal{A}\in\mathbb{K}_M^{k\times k}$ determinant function $\mathfrak{D}$, the following holds for any fixed column index $j$, 
    \begin{displaymath}
        \mathfrak{D}(\mathcal{A}) = \sum_{i=1}^k \mathrm{sgn}(i+j)\circ_M\mathcal{A}^{[i,j]}\circ_M \mathfrak{D}(\mathcal{A}^{[i|j]}) = \sum_{i=1}^k \mathcal{A}^{[i,j]}\circ_M \mathcal{C}^{[i,j]},
    \end{displaymath}
    
    where, $\mathcal{C}^{[i,j]} = \mathrm{sgn}(i+j)\circ_M \mathfrak{D}(\mathcal{A}^{[i|j]})$ can be called the $[i,j]$ tubal co-factor of $\mathcal{A}$. Also, it is clear that
    \begin{displaymath}
        \sum_{i=1}^k \mathcal{A}^{[i,\ell]}\circ_M\mathcal{C}^{[i,j]} = \delta_{j,\ell}(\mathbint{1},\mathbint{0}) \circ_M\mathfrak{D}(\mathcal{A}), \quad \text{where, } \delta_{j,\ell}(\mathbint{1},\mathbint{0}) = \begin{cases}
            \mathbint{1},\quad \text{if } j=\ell,\\
            \mathbint{0},\quad \text{else.}
        \end{cases}
    \end{displaymath}
    Now, we define the adjoint of $\mathcal{A}$ as, $\mathrm{adj}(\mathcal{A})^{[i,j]} = \mathcal{C}^{[i,j]}$. Hence, $\mathrm{adj}(\mathcal{A})\star_M\mathcal{A} = \mathfrak{D}(\mathcal{A})\star_M\mathcal{I}$.
    Now, applying transpose properties directly and using the above equation, we get $\mathcal{A}\star_M \mathrm{adj}(\mathcal{A}) = \mathfrak{D}(\mathcal{A})\star_M\mathcal{I}$.
\end{proof}

As we have defined the determinant, now we recall some of the basic definitions in commutative ring algebra, which we extend to tubal algebra for further analysis. 

\begin{definition}[Tubal McCoy rank]\label{def:mccoy-rank}
    Tensor $\mathcal{A}\in\mathbb{K}_M^{m\times n}$ is said to be of tubal McCoy rank, $\mathbint{r}$ is the largest $r\in\mathbb{N}$ such that, 
    \begin{displaymath}
    \mathrm{Ann}(I_r(\mathcal{A})) = \{\mathbint{0}\},
    \end{displaymath}
    where, $\mathrm{Ann}(I_\ell(\mathcal{A})) = \{\kappa\in\mathbb{K}_M|\kappa\circ_M x = \mathbint{0}, \forall x\in I_\ell(\mathcal{A})\}$ is the annihilator of the ideal, $I_\ell(\mathcal{A})$ of $\mathbb{K}_M$ generated using determinant of all $\ell\times\ell$ minors of tensor $\mathcal{A}$.
\end{definition}

\begin{definition}[$\ell$-tubal fitting ideal]\label{def:fitting-ideal}
    For a tensor $\mathcal{A}\in\mathbb{K}_M^{m\times n}$, the $\ell$-tubal fitting ideal is the ideal in $\mathbb{K}_M$ generated by the determinants of all $(m-\ell) \times (m-\ell)$ minors of $\mathcal{A}$. It is denoted by $\mathrm{Fitt}_\ell(\mathcal{A})$.
\end{definition}

Now, the tubal McCoy rank gives a basis for the analysis of rank and determinant in tubal algebra, as it gives the number of non-zero non-invertible tubes in the tensor.

\begin{theorem}[Invertibility of tensor]\label{thrm:invertibility_condition}
    If the tubal McCoy rank, $\mathbint{r}$, is equal to the tubal rank, $\mathfrak{r}_t$, of the tensor $\mathcal{A}\in\mathbb{K}_M^{k\times k}$, then the determinant of the tensor can't be a zero divisor. Moreover, if $\mathbint{r} = \mathfrak{r}_t = k$, then the tensor is invertible. Furthermore, the algebraic gap $\mathfrak{r}_t - \mathbint{r}$ represents the exact number of successive fitting ideals of $\mathcal{A}$ that are entirely composed of zero-divisors.
\end{theorem}
\begin{proof}
    Assume a tensor $\mathcal{A}\in\mathbb{K}_M^{k\times k}$. By \Cref{def:tubal-rank}, the tubal rank $\mathfrak{r}_t$ represents the total number of non-zero tubes which can be zero-divisors or invertible tubal scalars of $\mathbb{K}_M$, that is, the largest integer $\ell$ such that the ideal $I_\ell(\mathcal{A}) \neq \{\mathbint{0}\}$, where $I_\ell(\mathcal{A})$ is generated using determinant of all $\ell\times\ell$ minors of tensor $\mathcal{A}$. From \Cref{def:mccoy-rank}, the tubal McCoy rank $\mathbint{r}$ is the largest integer $\ell$ such that the annihilator $\mathrm{Ann}(I_\ell(\mathcal{A})) = \{\mathbint{0}\}$. 
    
    Assume $\mathbint{r} = \mathfrak{r}_t$. If $\mathfrak{r}_t < k$, then $I_k(\mathcal{A}) = \{\mathbint{0}\}$, that is, $\det(\mathcal{A}) = \mathbint{0}$, which is obvious. Alternatively, if $\mathbint{r} = \mathfrak{r}_t = k$, then $I_k(\mathcal{A})$, which is generated by $\det(\mathcal{A})$, satisfies $\mathrm{Ann}(I_k(\mathcal{A})) = \{\mathbint{0}\}$. This implies $\det(\mathcal{A})$ is not a zero-divisor and is invertible. Finally, from \Cref{def:fitting-ideal}, we know that $\mathrm{Fitt}_{k-\ell}(\mathcal{A}) = I_\ell(\mathcal{A})$. Hence, for any $\mathbint{r} < \ell\leq \mathfrak{r}_t$, the ideal satisfy, $I_\ell(\mathcal{A}) \neq \{\mathbint{0}\}$ (since $\ell \le \mathfrak{r}_t$) and $\mathrm{Ann}(I_\ell(\mathcal{A})) \neq \{\mathbint{0}\}$ (since $\ell > \mathbint{r}$). Hence, every ideal in this range is non-trivial but annihilated by some non-zero element, meaning it is composed entirely of zero-divisors. The number of indices satisfying this is exactly $\mathfrak{r}_t - \mathbint{r}$.
\end{proof}

\begin{figure}[htbp]
    \centering
    \includegraphics[width=0.75\linewidth]{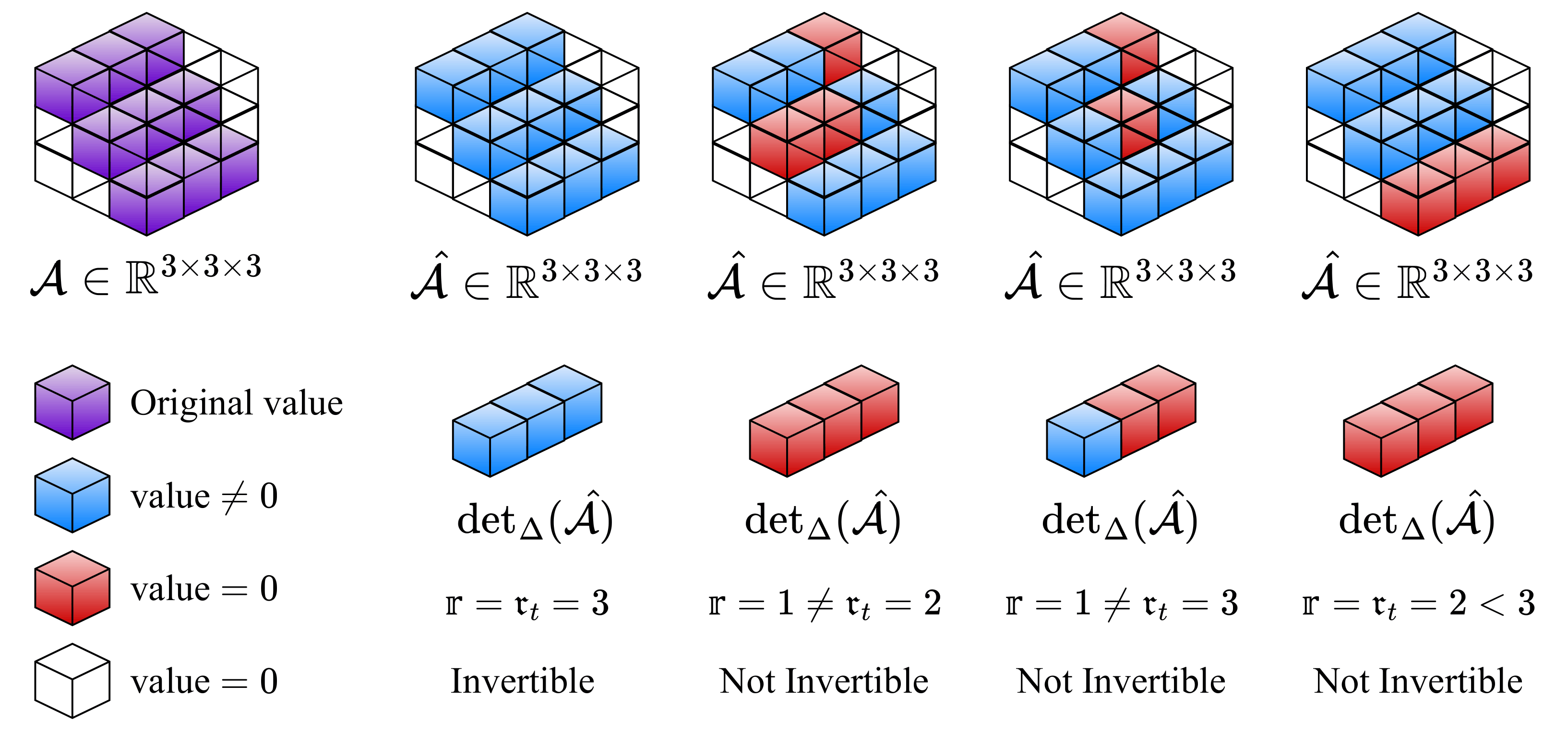}
    \caption{Different cases of determinants and their representations for visualization (according to \Cref{def:determinant} and \Cref{thrm:invertibility_condition}).}
    \label{fig:determinant-cases}
\end{figure}

To understand the working of the determinant and its properties, let's take an example for the same. For simplicity of calculations, we have assumed the $M$ matrix to be symmetric and orthogonal.
\begin{example}\label{example:determinant}
    Assuming the $M$ matrix is symmetric and orthogonal, such as
    \begin{displaymath}
    M^{-1} = M^\top = M = \begin{bmatrix}
        0.5&0.5&0.5&0.5\\
        0.5&-0.5&0.5&-0.5\\
        0.5&0.5&-0.5&-0.5\\
        0.5&-0.5&-0.5&0.5
    \end{bmatrix}.
    \end{displaymath}
    Now, assuming $5$ different cases,
    \begin{enumerate}
        \item Zero determinant case,
        \begin{displaymath}
        \mathcal{A}[:,:,1] = \mathcal{A}[:,:,2] = \mathcal{A}[:,:,3] = \mathcal{A}[:,:,4] = \begin{bmatrix}
            0.5&0.5\\1&1
        \end{bmatrix}.
        \end{displaymath}
        Here, after transformation $\hat{\mathcal{A}} = \mathfrak{M}(\mathcal{A})$ is given by, 
        \begin{displaymath}
        \hat{\mathcal{A}}[:,:,1] = \begin{bmatrix}
            1&1\\2&2
        \end{bmatrix},\quad\hat{\mathcal{A}}[:,:,2] = \hat{\mathcal{A}}[:,:,3] = \hat{\mathcal{A}}[:,:,4] = \begin{bmatrix}
            0&0\\0&0
        \end{bmatrix}.
        \end{displaymath}
        Here, we have $\mathrm{det}_\Delta(\hat{\mathcal{A}}) = [0,0,0,0]^\top$. 
        Hence, $\mathfrak{D}(\mathcal{A}) = \mathfrak{M}^{-1}(\mathrm{det}_\Delta(\hat{\mathcal{A}})) = [0,0,0,0]^\top$. In this case, $\mathbint{r} = 0$ and $\mathfrak{r}_t = 1$, and the tensor $\mathcal{A}$ is not invertible, which is also evident from the formulation of $\hat{\mathcal{A}}$.
        
        \item Non-zero determinant case,
        \begin{displaymath}
        \mathcal{A}[:,:,1] = \begin{bmatrix}
            2.5&0\\0&2.5
        \end{bmatrix},\quad \mathcal{A}[:,:,2] = \mathcal{A}[:,:,3] = \mathcal{A}[:,:,4] = \begin{bmatrix}
            0.5&0\\0&0.5
        \end{bmatrix}.
        \end{displaymath}
        Here, after transformation $\hat{\mathcal{A}} = \mathfrak{M}(\mathcal{A})$ is given by, 
        \begin{displaymath}
        \hat{\mathcal{A}}[:,:,1] = \begin{bmatrix}
            2&0\\0&2
        \end{bmatrix},\quad \hat{\mathcal{A}}[:,:,2] = \hat{\mathcal{A}}[:,:,3] = \hat{\mathcal{A}}[:,:,4] = \begin{bmatrix}
            1&0\\0&1
        \end{bmatrix}.
        \end{displaymath}
        Here, we have $\mathrm{det}_\Delta(\hat{\mathcal{A}}) = [4,1,1,1]^\top$.
        Hence, $\mathfrak{D}(\mathcal{A}) = \mathfrak{M}^{-1}(\mathrm{det}_\Delta(\hat{\mathcal{A}})) = [3.5,1.5,1.5,1.5]^\top$. In this case, $\mathbint{r} = \mathfrak{r}_t = k = 2$, hence, the tensor $\mathcal{A}$ is invertible, which is obvious as the slices of $\hat{\mathcal{A}}$ are diagonal.

        \item Zero divisor case,
        \begin{displaymath}
        \mathcal{A}[:,:,1] = \mathcal{A}[:,:,3] = \begin{bmatrix}
            2&0\\0&1
        \end{bmatrix},\quad \mathcal{A}[:,:,2] = \mathcal{A}[:,:,4] = \begin{bmatrix}
            0&0\\0&0
        \end{bmatrix}.
        \end{displaymath}
        Here, after transformation $\hat{\mathcal{A}} = \mathfrak{M}(\mathcal{A})$ is given by, 
        \begin{displaymath}
        \hat{\mathcal{A}}[:,:,1] = \hat{\mathcal{A}}[:,:,2] = \begin{bmatrix}
            2&0\\0&1
        \end{bmatrix},\quad \hat{\mathcal{A}}[:,:,3] = \hat{\mathcal{A}}[:,:,4] = \begin{bmatrix}
            0&0\\0&0
        \end{bmatrix}.
        \end{displaymath}
        Here, we have $\mathrm{det}_\Delta(\hat{\mathcal{A}}) = [2,2,0,0]^\top$.
        Hence, $\mathfrak{D}(\mathcal{A}) = \mathfrak{M}^{-1}(\mathrm{det}_\Delta(\hat{\mathcal{A}})) = [2,0,2,0]^\top$. In this case, $\mathbint{r} = 0$ and $\mathfrak{r}_t = 2$, and the tensor $\mathcal{A}$ is not invertible, which is also evident from the formulation of $\hat{\mathcal{A}}$.

        \item Non-invertible $\mathbint{r} = \mathfrak{r}_t$ case,
        \begin{displaymath}
        \mathcal{A}[:,:,1] = \begin{bmatrix}
            1&1\\1&1
        \end{bmatrix},\quad \mathcal{A}[:,:,2] = \mathcal{A}[:,:,3] = \mathcal{A}[:,:,4] = \begin{bmatrix}
            0&0\\0&0
        \end{bmatrix}.
        \end{displaymath}
        Here, after transformation $\hat{\mathcal{A}} = \mathfrak{M}(\mathcal{A})$ is given by, 
        \begin{displaymath}
        \hat{\mathcal{A}}[:,:,1] = \hat{\mathcal{A}}[:,:,2] = \hat{\mathcal{A}}[:,:,3] = \hat{\mathcal{A}}[:,:,4] = \begin{bmatrix}
            0.5&0.5\\0.5&0.5
        \end{bmatrix}.
        \end{displaymath}
        Here, we have $\mathrm{det}_\Delta(\hat{\mathcal{A}}) = [0,0,0,0]^\top$.
        Hence, $\mathfrak{D}(\mathcal{A}) = \mathfrak{M}^{-1}(\mathrm{det}_\Delta(\hat{\mathcal{A}})) = [0,0,0,0]^\top$. In this case, $\mathbint{r} = \mathfrak{r}_t = 1 < k$, hence, $\mathfrak{D}(\mathcal{A})$ is not a zero divisor. Also, the tensor $\mathcal{A}$ is not invertible.
    \end{enumerate}
\end{example}

As we have defined the determinant, tubal McCoy rank $\mathbint{r}$, and tubal rank $\mathfrak{r}_t$. Now we will map the McCoy defect ($\mathfrak{r}_t - \mathbint{r}$) formed due to the existence of zero divisors, which is required for a more precise understanding of the algebra.

\begin{theorem}[Resolution of McCoy defect]
For any tensor $\mathcal{A} \in \mathbb{K}_M^{k \times k}$ with a McCoy defect, that is, $\mathbint{r}(\mathcal{A}) < \mathfrak{r}_t(\mathcal{A})$, there exists a finite set of mutually orthogonal idempotent tubal scalars $\{\mathbf{e}_1, \dots, \mathbf{e}_p\} \subset \mathbb{K}_M$ satisfying $\sum_i \mathbf{e}_i = \mathbint{1}$, and an associated set of integers $\{r_1, \dots, r_p\}$ such that the projected tensors $\mathcal{A}_i = \mathbf{e}_i \star_M \mathcal{A}$ satisfy the following,
\begin{displaymath}
    \mathrm{Ann}(I_{r_i}(\mathcal{A}_i)) = \langle \mathbint{1} - \mathbf{e}_i \rangle \quad \text{and} \quad I_{r_i+1}(\mathcal{A}_i) = \{\mathbint{0}\},
\end{displaymath}
for all $1 \leq i \leq p$, where $\langle \mathbint{1} - \mathbf{e}_i \rangle$ is the principal ideal generated by $(\mathbint{1} - \mathbf{e}_i)$. Furthermore, the tubal rank satisfies $\mathfrak{r}_t(\mathcal{A}) = \max_{1 \leq i \leq p} \{r_i\}$.
\end{theorem}
\begin{proof}
Let $I_\ell(\mathcal{A})$ denote the ideal in $\mathbb{K}_M$ is generated using determinant of all $\ell\times\ell$ minors of tensor $\mathcal{A}$. Let $I_\ell(\mathcal{A}) = \langle \mathbf{u}_\ell \rangle$ where $\mathbf{u}_\ell$ is an idempotent. Now, as $I_{\ell+1}(\mathcal{A}) \subseteq I_\ell(\mathcal{A})$, the idempotents follows, $\mathbf{u}_\ell \circ_M \mathbf{u}_{\ell+1} = \mathbf{u}_{\ell+1}$. Let $\{r_1, r_2, \dots, r_p\}$ with $r_1 < r_2 < \dots < r_p$ be the complete set of integer indices where the sequence of idempotents strictly decreases, that is, $\mathbf{u}_{r_i} - \mathbf{u}_{r_i+1} \neq \mathbint{0}$, and where $\mathbf{u}_{r_p+1} = \mathbint{0}$. We define the elements, $\mathbf{e}_i = \mathbf{u}_{r_i} - \mathbf{u}_{r_i+1}$ for $1 \leq i \leq p$.

Since $\mathbf{u}_\ell$ is a descending chain of idempotents, it is trivial to verify that $\mathbf{e}_i \circ_M \mathbf{e}_j = \mathbint{0}$ for $i \neq j$. To ensure $\sum_i \mathbf{e}_i = \mathbint{1}$, if $\mathbf{u}_{r_1} \neq \mathbint{1}$, define an additional idempotent $\mathbf{e}_0 = \mathbint{1} - \mathbf{u}_{r_1}$. Let the projected tensor be $\mathcal{A}_i = \mathbf{e}_i \star_M \mathcal{A}$, then,
\begin{displaymath}
    I_\ell(\mathcal{A}_i) = \mathbf{e}_i \circ_M I_\ell(\mathcal{A}) = \mathbf{e}_i \circ_M \langle \mathbf{u}_\ell \rangle = \langle \mathbf{e}_i \circ_M \mathbf{u}_\ell \rangle.
\end{displaymath}
First, evaluating the ideal at index $r_i+1$,
\begin{displaymath}
    I_{r_i+1}(\mathcal{A}_i) = \langle\mathbf{e}_i \circ_M \mathbf{u}_{r_i+1}\rangle = \langle(\mathbf{u}_{r_i} - \mathbf{u}_{r_i+1}) \circ_M \mathbf{u}_{r_i+1}\rangle = \langle\mathbint{0}\rangle = \mathbint{0}.
\end{displaymath}
Now, evaluating the ideal at index $r_i$,
\begin{displaymath}
    I_{r_i}(\mathcal{A}_i) = \langle\mathbf{e}_i \circ_M \mathbf{u}_{r_i}\rangle = \langle(\mathbf{u}_{r_i} - \mathbf{u}_{r_i+1}) \circ_M \mathbf{u}_{r_i}\rangle = \langle\mathbf{u}_{r_i} - \mathbf{u}_{r_i+1}\rangle = \langle\mathbf{e}_i\rangle.
\end{displaymath}
Hence, $\mathrm{Ann}(I_{r_i}(\mathcal{A}_i)) = \langle \mathbint{1} - \mathbf{e}_i \rangle$. Now, by \Cref{def:tubal-rank}, the tubal rank $\mathfrak{r}_t(\mathcal{A})$ is the largest integer $q$ such that $I_q(\mathcal{A}) \neq \{\mathbint{0}\}$. Since $I_q(\mathcal{A}) = \langle \mathbf{u}_q \rangle$, hence it is maximum $q$ where $\mathbf{u}_q \neq \mathbint{0}$. Because $\mathbf{u}_{r_p+1} = \mathbint{0}$, the maximum non-zero index is exactly $r_p$. Thus, $\mathfrak{r}_t(\mathcal{A}) = r_p = \max_i \{r_i\}$.
\end{proof}

Assuming a signal perfectly exists in some channels, but is zero in others, giving rise to structural misalignment. Mathematically, this misalignment between the total active components (tubal rank) and the components that are active across all channels (tubal McCoy rank) is called the McCoy defect here. While for every $M\in\mathrm{GL}_M(\mathbb{F})$, the tensors form an algebra, for practical case analysis of the determinants, it plays a significant role. For e.g., if time-frequency analysis is required, then DFT or Lifting-matrices are better choices of $M$ to be taken into consideration. It is also important to note that the $M$-transformation works similarly to convolutions. Hence, direct analysis of determinants and other properties of tensors using arbitrary $M$ is not preferred.

\subsection{Linear Dependence and Tubal Rank-Nullity Theorem}
\begin{definition}[Linear independence]\label{def:linear-independent}
    Let $\mathscr{S}_{\mathrm{LI}} = \{\vec{\mathcal{X}}_1, \vec{\mathcal{X}}_2, \dots, \vec{\mathcal{X}}_m\}$ be set of columns where, $\vec{\mathcal{X}}_i\in\mathbb{K}_M^{n\times 1}$. Let $\mathbf{c}_1, \mathbf{c}_2,\dots, \mathbf{c}_m\in\mathbb{K}_M$ be the tubal scalars. Then the set $\mathscr{S}_{\mathrm{LI}}$ is linearly independent if the solution of $\sum_{i=1}^m \vec{\mathcal{X}}_i\star_M \mathbf{c}_i = \vec{\mathcal{O}}$ is $\mathbf{c}_1 = \mathbf{c}_2 = \dots = \mathbf{c}_n = \mathbint{0}$, where, $\vec{\mathcal{O}}$ is the tubal vector of scalars $\mathbint{0}$. 
\end{definition}

\begin{definition}[Linear dependence]\label{def:linear-dependent}
    Let $\mathscr{S}_{\mathrm{LD}} = \{\vec{\mathcal{X}}_1, \vec{\mathcal{X}}_2, \dots, \vec{\mathcal{X}}_m\}$ be set of columns where, $\vec{\mathcal{X}}_i\in\mathbb{K}_M^{n\times 1}$. Let $\mathbf{c}_1, \mathbf{c}_2,\dots, \mathbf{c}_m\in\mathbb{K}_M$ be the tubal scalars. Then the set $\mathscr{S}_{\mathrm{LD}}$ is linearly dependent if there exists a solution of $\sum_{i=1}^m \vec{\mathcal{X}}_i\star_M \mathbf{c}_i = \vec{\mathcal{O}}$ with atleast one coefficient $\mathbf{c}_i\neq \mathbint{0}$.
\end{definition}

As it is clear from \Cref{def:linear-dependent}, the coefficients $\{\mathbf{c}_1, \mathbf{c}_2,\dots, \mathbf{c}_m\}$ can't be all $\mathbint{0}$. So, the non-zero coefficient can be either invertible or a zero divisor. This gives rise to the sub-classification of linear dependence into two categories as follows, 
\begin{definition}[Linear absolute dependence]\label{def:linear-absolute-dependent}
    Let $\mathscr{S}_{\mathrm{LAD}} = \{\vec{\mathcal{X}}_1, \vec{\mathcal{X}}_2, \dots, \vec{\mathcal{X}}_m\}$ be set of columns where, $\vec{\mathcal{X}}_i\in\mathbb{K}_M^{n\times 1}$. Let $\mathbf{c}_1, \mathbf{c}_2,\dots, \mathbf{c}_m\in\mathbb{K}_M$ be the tubal scalars. Then the set $\mathscr{S}_{\mathrm{LAD}}$ is absolutely dependent if there exists a solution of $\sum_{i=1}^m \vec{\mathcal{X}}_i\star_M \mathbf{c}_i = \vec{\mathcal{O}}$ with atleast one coefficient $\mathbf{c}_i$ is invertible and no zero-divisors.
\end{definition}

\begin{definition}[Linear tubal dependence]\label{def:linear-tubal-dependent}
    Let $\mathscr{S}_{\mathrm{LTD}} = \{\vec{\mathcal{X}}_1, \vec{\mathcal{X}}_2, \dots, \vec{\mathcal{X}}_m\}$ be set of columns where, $\vec{\mathcal{X}}_i\in\mathbb{K}_M^{n\times 1}$. Let $\mathbf{c}_1, \mathbf{c}_2,\dots, \mathbf{c}_m\in\mathbb{K}_M$ be the tubal scalars. Then the set $\mathscr{S}_{\mathrm{LTD}}$ is tubal dependent if it is not linearly independent, and every non-trivial solution of $\sum_{i=1}^m \vec{\mathcal{X}}_i\star_M \mathbf{c}_i = \vec{\mathcal{O}}$ consists of at least one zero divisor. 
\end{definition}

\begin{figure}[htbp]
    \centering
    \includegraphics[width=0.8\linewidth]{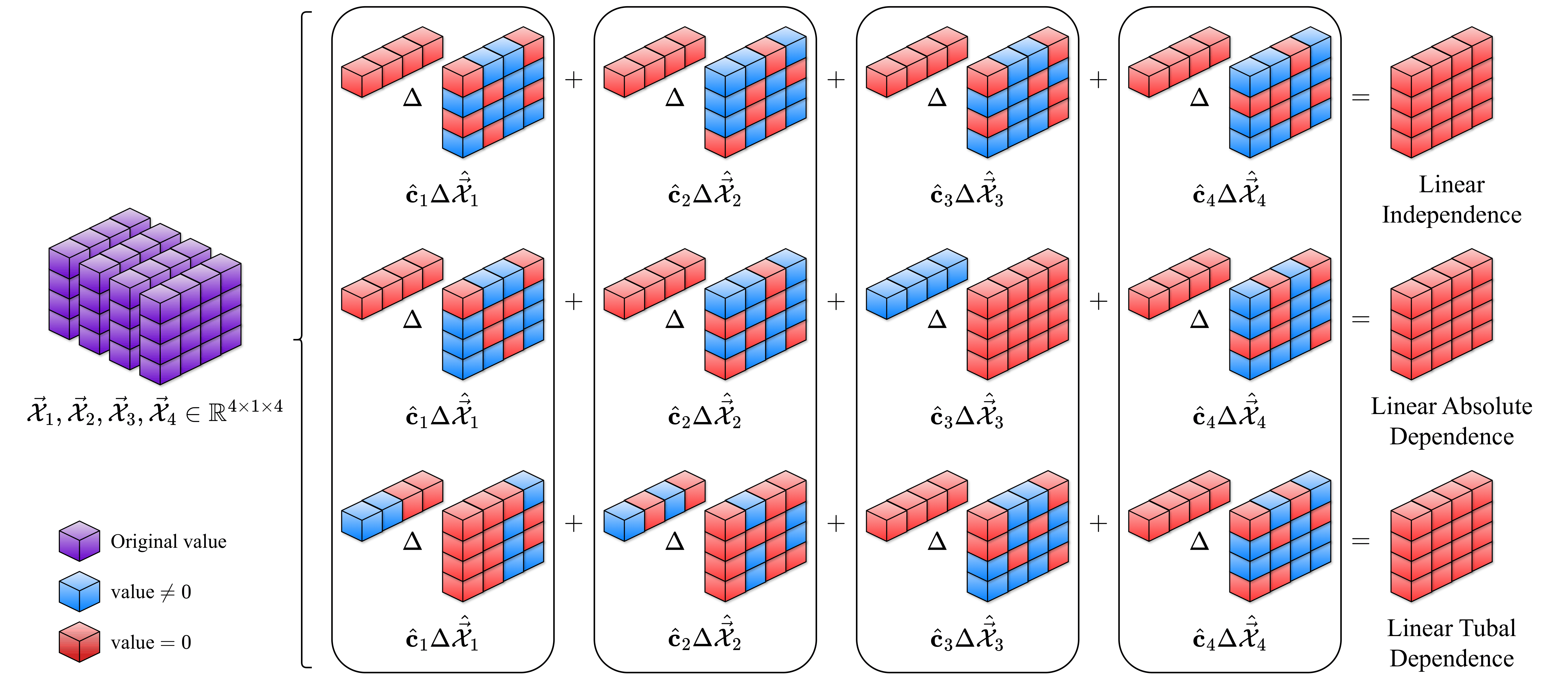}
    \caption{Visualization of linear independence (\Cref{def:linear-independent}), absolute dependence (\Cref{def:linear-absolute-dependent}), and tubal dependence (\Cref{def:linear-tubal-dependent}) under tubal tensor algebra.}
    \label{fig:independent-cases}
\end{figure}

It is important to note that due to the existence of zero-divisors, the classical Rank-Nullity Theorem will not hold for tubal tensor algebra. Hence, a more rigorous dimension analysis is required for finding the tubal rank-nullity theorem. 

\begin{definition}[Tubal dimension]\label{def:tubal-dimension}
    Let $\mathscr{W}$ be a finitely generated projective submodule of tubal vectors in $\mathbb{K}_M^{n}$. The tubal dimension of $\mathscr{W}$, denoted $\mathrm{Dim}_{\mathbb{K}_M}(\mathscr{W})$, is defined as the tubal trace of the unique orthogonal projector $\mathcal{P}_{\mathscr{W}} \in \mathbb{K}_M^{n \times n}$ onto $\mathscr{W}$, which is given by,
    \begin{displaymath}
    \mathrm{Dim}_{\mathbb{K}_M}(\mathscr{W}) = \texttt{t-Tr}(\mathcal{P}_{\mathscr{W}}) = \sum_{i=1}^n \mathcal{P}_{\mathscr{W}}^{[i,i]} \in \mathbb{K}_M.
    \end{displaymath}
    where tubal trace os tensor $\mathcal{A}\in\mathbb{K}_M^{n\times n}$ is defined as, $\texttt{t-Tr}(\mathcal{A}) = \sum_{i=1}^n \mathcal{A}^{[i,i]}$.
\end{definition}

Note that here we consider the projective submodule, not the free module, which is very restrictive, because zero vectors can also lie in this submodule. Due to the existence of zero-divisors, directly counting the number of tubal vectors as dependent and independent will not work. Hence, a different approach for finding the dimension is utilized for generalizing the rank-nullity theorem to tensors.

\begin{theorem}[Tubal rank-nullity theorem]\label{thrm:general-rank-nullity-thrm}
    Let a tensor $\mathcal{A} \in \mathbb{K}_M^{m \times n}$ with $\mathscr{R}(\mathcal{A})$ denote its range space and $\mathscr{N}(\mathcal{A})$ denote its null space. Then,
    \begin{displaymath}
    \mathrm{Dim}_{\mathbb{K}_M}(\mathscr{R}(\mathcal{A})) + \mathrm{Dim}_{\mathbb{K}_M}(\mathscr{N}(\mathcal{A})) = n \ast \mathbint{1}.
    \end{displaymath}
\end{theorem}
\begin{proof}
    Because $\mathbb{K}_M$ is Von Neumann regular, every tensor $\mathcal{A} \in \mathbb{K}_M^{m \times n}$ has a unique tubal Moore-Penrose pseudoinverse. The orthogonal projector onto the row space of $\mathcal{A}$ is given by $\mathcal{P}_{\mathrm{row}} = \mathcal{A}^\dagger \star_M \mathcal{A}$ and the orthogonal projector onto the null space $\mathscr{N}(\mathcal{A})$ is given by, $\mathcal{P}_{\mathrm{null}} = \mathcal{I}_n - \mathcal{A}^\dagger \star_M \mathcal{A}$, where $\mathcal{I}_n \in \mathbb{K}_M^{n \times n}$ is the identity tensor. Now applying the tubal trace, we get, $\texttt{t-Tr}(\mathcal{P}_{\mathrm{null}}) = \texttt{t-Tr}(\mathcal{I}_n - \mathcal{A}^\dagger \star_M \mathcal{A}) = \texttt{t-Tr}(\mathcal{I}_n) - \texttt{t-Tr}(\mathcal{A}^\dagger \star_M \mathcal{A})$. Now, according to the cyclic property of the tubal trace, $\texttt{t-Tr}(\mathcal{A}^\dagger \star_M \mathcal{A}) = \texttt{t-Tr}(\mathcal{A} \star_M \mathcal{A}^\dagger)$. As $\mathcal{A}\star_M\mathcal{A}^\dagger$ is the orthogonal projector onto range space $\mathscr{R}(\mathcal{A})$ and \Cref{def:tubal-dimension}, we have, $\mathrm{Dim}_{\mathbb{K}_M}(\mathscr{N}(\mathcal{A})) = n \ast \mathbint{1} - \mathrm{Dim}_{\mathbb{K}_M}(\mathscr{R}(\mathcal{A}))$.
\end{proof}

\begin{definition}[Tubal span]
    Let $\mathscr{S} = \{\vec{\mathcal{X}}_1, \vec{\mathcal{X}}_2, \dots, \vec{\mathcal{X}}_m\}$ be a finite set of tubal vectors in $\mathbb{K}_M^{n}$. The tubal span of $\mathscr{S}$, denoted as $\texttt{span}_{\mathbb{K}_M}(\mathscr{S})$, is defined as the set of all finite tubal linear combinations of the vectors in $\mathscr{S}$,
    \begin{displaymath}
    \texttt{span}_{\mathbb{K}_M}(\mathscr{S}) = \left\{ \sum_{i=1}^m \vec{\mathcal{X}}_i \star_M \mathbf{c}_i ~|~ \mathbf{c}_i \in \mathbb{K}_M \right\}.
    \end{displaymath}
    The set $\texttt{span}_{\mathbb{K}_M}(\mathscr{S})$ forms a finitely generated submodule of $\mathbb{K}_M^{n}$. If $\texttt{span}_{\mathbb{K}_M}(\mathscr{S}) = \mathscr{W}$, we say that $\mathscr{S}$ is a tubal generating set for the submodule $\mathscr{W}$.
\end{definition}

Unlike vector spaces over a field, the ring of tubal scalars $\mathbb{K}_M$ contains zero divisors. Consequently, while $\texttt{span}_{\mathbb{K}_M}(\mathscr{S})$ is well-defined as a finitely generated submodule, it is not guaranteed to be a free submodule. The presence of zero divisors implies that $\mathscr{W}$ may contain zero divisors, and a minimal generating set for $\mathscr{W}$ cannot strictly be assumed to form a basis.

\begin{definition}[Tubal range space]\label{def:tubal-range-space}
    If $\vec{\mathcal{A}}^{[:,1]}, \vec{\mathcal{A}}^{[:,2]}, \dots, \vec{\mathcal{A}}^{[:,n]}$ represent the tubal column vectors of tensor $\mathcal{A}\in\mathbb{K}_M^{m\times n}$, the range space is exactly their tubal span, which is given by 
    \begin{displaymath}
        \mathscr{R}(\mathcal{A}) = \texttt{span}_{\mathbb{K}_M}\{\{\vec{\mathcal{A}}^{[:,i]}\}_{i=1}^n\}.
    \end{displaymath}
\end{definition}

\begin{definition}[Tubal null space]\label{def:tubal-null-space}
    The null space of the tensor $\mathcal{A}\in\mathbb{K}_M^{m\times n}$ is set of all tubal vectors in $\mathbb{K}_M^{n\times 1}$ such that,
    \begin{displaymath}
        \mathscr{N}(\mathcal{A}) = \left\{ \vec{\mathcal{X}} \in \mathbb{K}_M^{n} ~|~ \mathcal{A} \star_M \vec{\mathcal{X}} = \vec{\mathcal{O}} \right\},
    \end{displaymath}
    where $\vec{\mathcal{O}} \in \mathbb{K}_M^{m \times 1}$ is the zero tubal vector with element being $\mathbint{0}$.
\end{definition}

The subspaces formed in \Cref{def:tubal-range-space} and \Cref{def:tubal-null-space} are equivalent to those formed in \Cref{def:svd-range-null-space}. By \Cref{def:tubal-range-space}, $\vec{\mathcal{Y}} \in \mathscr{R}(\mathcal{A})$ implies $\vec{\mathcal{Y}} = \mathcal{A} \star_M \vec{\mathcal{X}}$ for some $\vec{\mathcal{X}} \in \mathbb{K}_M^{n}$. Substituting the SVD given in \Cref{eq:Ax-svd-expansion} yields $\vec{\mathcal{Y}} = \sum_{i=1}^{\min\{m, n\}} \vec{\mathcal{U}}_i \star_M (\sigma_i \circ_M \mathbf{d}_i)$, where $\mathbf{d}_i = \vec{\mathcal{V}}_i^\top \star_M \vec{\mathcal{X}}$. Since $\sigma_i = \mathbint{0}$ for $i > k_1+k_2$, those terms vanishes to $\vec{\mathcal{O}}$, leaving the SVD range space. Similarly, $\vec{\mathcal{X}} \in \mathscr{N}(\mathcal{A})$ requires $\mathcal{A} \star_M \vec{\mathcal{X}} = \vec{\mathcal{O}}$, which simplifies to $\sum_{i} \vec{\mathcal{U}}_i \star_M (\sigma_i \circ_M \mathbf{d}_i) = \vec{\mathcal{O}}$. Due to the orthogonality of $\mathcal{U}$, this requires $\sigma_i \circ_M \mathbf{d}_i = \mathbint{0}$ for all $i$. For invertible tubes ($i \leq k_1$), $\mathbf{d}_i = \mathbint{0}$. For zero tubes ($i > k_1+k_2$), $\mathbf{d}_i$ can be invertible or non-invertible. For non-invertible non-zero tubes ($k_1 < i \leq k_1+k_2$), $\mathbf{d}_i$ must belong to $\mathrm{Ann}(\sigma_i)$. Multiplying by $\mathcal{V}$ reconstructs exactly the SVD null space.

\begin{example}\label{example:rank-nullity}
Evaluating three cases for a tensor $\mathcal{A} \in \mathbb{R}^{3 \times 3 \times 4}$, whose tubal columns are denoted as $\vec{\mathcal{X}}_1, \vec{\mathcal{X}}_2, \vec{\mathcal{X}}_3$.

\begin{enumerate}
    \item Linear independence case,
        \begin{displaymath}
            \mathcal{A}[:,:,1] = \begin{bmatrix}
                2&0&0\\0&2&0\\0&0&2
            \end{bmatrix},\quad
            \mathcal{A}[:,:,2] = \mathcal{A}[:,:,3] = \mathcal{A}[:,:,4] = \begin{bmatrix}
                0&0&0\\0&0&0\\0&0&0
            \end{bmatrix}.
            \end{displaymath}
            Here, after transformation $\hat{\mathcal{A}} = \mathfrak{M}(\mathcal{A})$ is given by, 
            \begin{displaymath}
            \hat{\mathcal{A}}[:,:,1] = \hat{\mathcal{A}}[:,:,2] = \hat{\mathcal{A}}[:,:,3] = \hat{\mathcal{A}}[:,:,4] = \begin{bmatrix}
                1&0&0\\0&1&0\\0&0&1
            \end{bmatrix}.
        \end{displaymath}
        Solving $\sum_{j=1}^3 \vec{\mathcal{X}}_j \star_M \mathbf{c}_j = \vec{\mathcal{O}}$ requires $\mathbf{c}_1 = \mathbf{c}_2 = \mathbf{c}_3 = \mathbint{0}$. 
        The orthogonal projector onto the range space is $\mathcal{P}_{\mathscr{R}} = \mathcal{I}_3$. Applying the \Cref{def:tubal-dimension}, we get, $\mathrm{Dim}_{\mathbb{K}_M}(\mathscr{R}) = \mathbint{1} + \mathbint{1} + \mathbint{1} = [6, 0, 0, 0]^\top$ and $\mathrm{Dim}_{\mathbb{K}_M}(\mathscr{N}) = \mathbint{0} = [0, 0, 0, 0]^\top$. Hence, $\mathrm{Dim}_{\mathbb{K}_M}(\mathscr{R}) + \mathrm{Dim}_{\mathbb{K}_M}(\mathscr{N}) = [6, 0, 0, 0]^\top = 3 \ast \mathbint{1}$, satisfying the \Cref{thrm:general-rank-nullity-thrm}.

    \item Linear absolute dependence case,
        \begin{displaymath}
            \mathcal{A}[:,:,1] = \begin{bmatrix}
                2&0&2\\0&2&2\\0&0&0
            \end{bmatrix},\quad
            \mathcal{A}[:,:,2] = \mathcal{A}[:,:,3] = \mathcal{A}[:,:,4] = \begin{bmatrix}
                0&0&0\\0&0&0\\0&0&0
            \end{bmatrix}.
            \end{displaymath}
            Here, after transformation $\hat{\mathcal{A}} = \mathfrak{M}(\mathcal{A})$ is given by, 
            \begin{displaymath}
            \hat{\mathcal{A}}[:,:,1] = \hat{\mathcal{A}}[:,:,2] = \hat{\mathcal{A}}[:,:,3] = \hat{\mathcal{A}}[:,:,4] = \begin{bmatrix}
                1&0&1\\0&1&1\\0&0&0
            \end{bmatrix}.
        \end{displaymath}
        It is clear that, $\vec{\mathcal{X}}_3$ is the sum of $\vec{\mathcal{X}}_1$ and $\vec{\mathcal{X}}_2$. Hence, $\mathbf{c}_1 = \mathbint{1}$, $\mathbf{c}_2 = \mathbint{1}$, $\mathbf{c}_3 = -\mathbint{1}$. As the coefficients are invertible, this is absolute dependence. Applying the \Cref{def:tubal-dimension}, we get,
        $\mathrm{Dim}_{\mathbb{K}_M}(\mathscr{R}) = 2 \ast \mathbint{1} = [4, 0, 0, 0]^\top$ and $\mathrm{Dim}_{\mathbb{K}_M}(\mathscr{N}) = 1 \ast \mathbint{1} = [2, 0, 0, 0]^\top$
        Hence, $\mathrm{Dim}_{\mathbb{K}_M}(\mathscr{R}) + \mathrm{Dim}_{\mathbb{K}_M}(\mathscr{N}) = [6, 0, 0, 0]^\top = 3 \ast \mathbint{1}$, satisfying the \Cref{thrm:general-rank-nullity-thrm}.

        \item Linear tubal dependence case,
        \begin{displaymath}
            \mathcal{A}[:,:,1] = \begin{bmatrix}
                2&0&1.5\\0&2&1.5\\0&0&0.5
            \end{bmatrix},\quad
            \mathcal{A}[:,:,2] = \mathcal{A}[:,:,3] = \mathcal{A}[:,:,4] = \begin{bmatrix}
                0&0&-0.5\\0&0&-0.5\\0&0&0.5
            \end{bmatrix}.
            \end{displaymath}
            Here, after transformation $\hat{\mathcal{A}} = \mathfrak{M}(\mathcal{A})$ is given by, 
            \begin{displaymath}
            \hat{\mathcal{A}}[:,:,1] = \begin{bmatrix}
                1&0&0\\0&1&0\\0&0&1
            \end{bmatrix},\quad \hat{\mathcal{A}}[:,:,2] = \hat{\mathcal{A}}[:,:,3] = \hat{\mathcal{A}}[:,:,4] = \begin{bmatrix}
                1&0&1\\0&1&1\\0&0&0
            \end{bmatrix}.
        \end{displaymath}
        It is clear that, $\vec{\mathcal{X}}_3$ depends on $\vec{\mathcal{X}}_1$ and $\vec{\mathcal{X}}_2$ \textit{only} on the sub-ideals $\mathbf{e}_i$ for $i \in \{2,3,4\}$. Hence, in the transform domain, the coefficients are $\hat{\mathbf{c}}_1 = \hat{\mathbf{c}}_2 = [0,1,1,1]^\top$ and $\hat{\mathbf{c}}_3 = [0,-1,-1,-1]^\top$. Hence $\mathbf{c}_1 = \mathbf{c}_2 = [1.5, -0.5, -0.5, -0.5]^\top$ and $\mathbf{c}_3 = [-1.5, 0.5, 0.5, 0.5]^\top$. As the coefficients are non-invertible and $\neq\mathbint{0}$, this is tubal dependence. Applying the \Cref{def:tubal-dimension}, we get,
        $\mathrm{Dim}_{\mathbb{K}_M}(\mathscr{R}) = \mathbf{e}_1 \ast 3 + (\mathbf{e}_2+\mathbf{e}_3+\mathbf{e}_4) \ast 2 = [4.5, 0.5, 0.5, 0.5]^\top$ and $\mathrm{Dim}_{\mathbb{K}_M}(\mathscr{N}) = \mathbf{e}_1 \ast 0 + (\mathbf{e}_2+\mathbf{e}_3+\mathbf{e}_4) \ast 1 = [1.5, -0.5, -0.5, -0.5]^\top$
        Hence, $\mathrm{Dim}_{\mathbb{K}_M}(\mathscr{R}) + \mathrm{Dim}_{\mathbb{K}_M}(\mathscr{N}) = [6, 0, 0, 0]^\top = 3 \ast \mathbint{1}$, satisfying the \Cref{thrm:general-rank-nullity-thrm}.
\end{enumerate}
\end{example}

\section{Spectral Analysis}\label{sec:spectral-analysis}
For spectral analysis, we are assuming only real-valued tubal scalars, leaving the rest of the analysis ($\mathbb{C}$ and $\mathbb{H}$) for future work.

\subsection{Tubal Eigenvalues and Eigenvectors}
\begin{definition}[Tubal eigenvalues and eigenvectors]\label{def:eigval-and-eigvect}
    Let tensor $\mathcal{A} \in \mathbb{K}_M^{k \times k}$. A tubal scalar $\lambda \in \mathbb{K}_M$ is called a tubal eigenvalue of $\mathcal{A}$ if there exists a tubal vector $\vec{\mathcal{X}} \in \mathbb{K}_M^{k \times 1}$ such that,
    \begin{displaymath}
        \mathcal{A} \star_M \vec{\mathcal{X}} = \lambda \circ_M \vec{\mathcal{X}},
    \end{displaymath}
    and the annihilator of the ideal generated by the components of $\vec{\mathcal{X}}$ is trivial, that is, $\mathrm{Ann}(\vec{\mathcal{X}}^{[1]}, \vec{\mathcal{X}}^{[2]}, \dots, \vec{\mathcal{X}}^{[k]}) = \{\mathbint{0}\}$. The vector $\vec{\mathcal{X}}$ is called a regular tubal eigenvector.
\end{definition}

\begin{definition}[Tubal eigensubmodule]
Let $\lambda \in \mathbb{K}_M$ be a valid tubal eigenvalue of $\mathcal{A} \in \mathbb{K}_M^{k \times k}$. The tubal eigensubmodule associated with $\lambda$, denoted as $E(\lambda)$, is the set of all tubal vectors $\vec{\mathcal{X}}$ satisfying $\mathcal{A} \star_M \vec{\mathcal{X}} = \lambda \circ_M \vec{\mathcal{X}}$. 
\end{definition}

In linear algebra over a field, an eigenspace excludes the zero vector. However, because $E(\lambda)$ is a module over the ring $\mathbb{K}_M$, it contains non-trivial vectors that possess zero-divisor annihilators. If $\vec{\mathcal{X}}$ is a regular tubal eigenvector in $E(\lambda)$, any scalar multiplication $\vec{\mathcal{Y}} = \mathbf{c} \circ_M \vec{\mathcal{X}}$, where $\mathbf{c}$ is a zero-divisor, will result in a degenerate tubal vector. Hence, while degenerate tubal vectors exist in $E(\lambda)$, they can't define the eigenvalue independently.

\begin{definition}[Tubal characteristic polynomial]\label{def:tubal-char-poly}
    For a tensor $\mathcal{A}\in\mathbb{K}_M^{k\times k}$, the tubal characteristic polynomial $p_{\mathcal{A}}(\lambda)$ is defined in polynomial ring $\mathbb{K}_M[\lambda]$ as,
    \begin{displaymath}
        p_{\mathcal{A}}(\lambda) = \mathfrak{D}(\lambda \star_M \mathcal{I} - \mathcal{A}) = \lambda^k + \mathbf{c}_{k-1}\circ_M\lambda^{k-1} + \dots + \mathbf{c}_1\circ_M\lambda + \mathbf{c}_0,
    \end{displaymath}
    where $\lambda, \mathbf{c}_i \in \mathbb{K}_M$, $\mathcal{I}$ is the identity tensor, and $\mathfrak{D}$ is the determinant function.
\end{definition}

It is important to note that the tubal characteristic polynomial $p_{\mathcal{A}}(\lambda)$ can exhibit ``partial'' roots due to the existence of zero-divisors. Specifically, if $p_{\mathcal{A}}(\lambda_0) = \mathbf{z}$ where $\mathbf{z} \neq \mathbint{0}$ and $\mathrm{Ann}(\mathbf{z}) \neq \{\mathbint{0}\}$ for a tubal eigenvalue $\lambda_0\in\mathbb{K}_M$. Then, the entire eigenspace associated with $\lambda_0$ strictly exhibits Linear Tubal Dependence. Also, in this case $\mathrm{Ann}(\vec{\mathcal{X}}^{[1]}, \vec{\mathcal{X}}^{[2]}, \dots, \vec{\mathcal{X}}^{[k]}) \neq \{\mathbint{0}\}$, that is, corresponding tubal eigenvector is not regular. 

\begin{remark}
    The polynomial equation $p_{\mathcal{A}}(\lambda) = \mathbint{0}$ may have up to $k^N$ distinct roots in $\mathbb{K}_M$ due to existing permutations under the decoupled nature of the transform domain rather than being bounded by $k$. This is due to the fact that an eigenvalue can be selected from $k$ possibilities across $N$ independent channels, yielding $k^N$ total combinations.
\end{remark}

While there are $k^N$ total eigenvalues, they do not form a totally ordered set under the partial ordering $\succeq$ due to the existence of zero-divisors. Hence, $k$ principal tubal eigenvalues are required for further analysis.
\begin{definition}[Principal tubal eigenvalues]
    We define the $k$ principal tubal eigenvalues to be systematically sorted, such that $\lambda_1^{(j)} \geq \lambda_2^{(j)} \geq \dots \geq \lambda_k^{(j)}$ for all $1\leq j\leq N$. Then, the $i$th principal tubal eigenvalue is constructed as,
    \begin{displaymath}
        \lambda_i = \sum_{j=1}^N \lambda_i^{(j)} \circ_M \mathbf{e}_j,
    \end{displaymath}
    where, $\mathbf{e}_j$ is the $j$th orthogonal idempotent tubal scalar. Now, these $k$ principal eigenvalues form a totally ordered chain in $\mathbb{K}_M$ as follows, $\lambda_1 \succeq \lambda_2 \succeq \dots \succeq \lambda_k$.
\end{definition}

\begin{theorem}[Tubal Cayley-Hamilton theorem]
    Every tensor $\mathcal{A}\in\mathbb{K}_M^{k\times k}$ satisfies its own tubal characteristic polynomial, that is,
    \begin{displaymath}
        p_{\mathcal{A}}(\mathcal{A}) = \mathcal{A}^k + \mathbf{c}_{k-1}\star_M\mathcal{A}^{k-1} + \dots + \mathbf{c}_1\star_M\mathcal{A} + \mathbf{c}_0\star_M\mathcal{I} = \mathcal{O},
    \end{displaymath}
    where $\mathcal{O}\in\mathbb{K}_M^{k\times k}$ is the zero tensor with each element being $\mathbint{0}$.
\end{theorem}
\begin{proof}
    As the tubal characteristic polynomial $p_{\mathcal{A}}(\lambda)$ is defined in polynomial ring $\mathbb{K}_M[\lambda]$. The adjoint tensor $\mathrm{adj}(\lambda \star_M \mathcal{I} - \mathcal{A})$ exists in $\mathbb{K}_M[\lambda]^{k\times k}$. Using \Cref{lem:adjoint-determinant-rel} and evaluating at $\lambda = \mathcal{A}$ yields the zero tensor, that is, $p_{\mathcal{A}}(\mathcal{A}) = \mathcal{O}$.
\end{proof}

Assuming a tubal eigenvalue $\lambda$ and corresponding regular tubal eigenvector $\vec{\mathcal{X}}$ of tensor $\mathcal{A}$. It is clear from \Cref{def:eigval-and-eigvect} that $\mathcal{A}\star_M\vec{\mathcal{X}} = \lambda\star_M\vec{\mathcal{X}}$ means $\lambda\star_M\vec{\mathcal{X}} - \mathcal{A}\star_M\vec{\mathcal{X}} = \vec{\mathcal{O}}$, which can be written as $(\lambda\star_M\mathcal{I} - \mathcal{A})\star_M\vec{\mathcal{X}} = \vec{\mathcal{O}}$. Now assume $\mathcal{B} = \lambda\star_M\mathcal{I} - \mathcal{A}$, multiply both sides by $\texttt{adj}(\mathcal{B})$, we get, $\texttt{adj}(\mathcal{B})\star_M\mathcal{B}\star_M\vec{\mathcal{X}} = \vec{\mathcal{O}}$. Now utilizing \Cref{lem:adjoint-determinant-rel}, we get, $\mathfrak{D}(\mathcal{B})\star_M\mathcal{I}\star_M\vec{\mathcal{X}} = \mathfrak{D}(\mathcal{B})\star_M\vec{\mathcal{X}} = \vec{\mathcal{O}}$. As $\vec{\mathcal{X}}$ is regular tubal eigenvector, hence, $\mathfrak{D}(\mathcal{B}) = \mathfrak{D}(\lambda\star_M\mathcal{I} - \mathcal{A}) = \mathbint{0}$. Hence from \Cref{def:tubal-char-poly} we can say that, $p_{\mathcal{A}}(\lambda) = \mathbint{0}$ for a tubal eigenvalue $\lambda$.

\begin{lemma}
    For any specific permutation of tubal eigenvalues of tensor $\mathcal{A}\in\mathbb{K}_M^{k\times k}$, the following holds,
    \begin{subequations}
        \begin{align}
        &\mathfrak{D}(\mathcal{A}) = {\prod}_{i=1}^k\lambda_i = \lambda_1\circ_M\lambda_2\circ_M\dots\circ_M\lambda_k,\label{eq:determinant-lambda-multiplication}\\
        &\texttt{t-Tr}(\mathcal{A}) = \sum_{i=1}^k\lambda_i = \lambda_1 + \lambda_2 + \dots + \lambda_k,\label{eq:trace-lambda-addition}
    \end{align}
    \end{subequations}
    where, $\lambda_i$ is the $i$th principal tubal eigenvalue.
\end{lemma}
\begin{proof}
    As the principal tubal eigenvalues $(\lambda_1,\lambda_2,\dots,\lambda_k)$ satisfies the tubal characteristic polynomial $p_{\mathcal{A}}(\lambda)$. Hence, we can write,
    \begin{displaymath}
        p_{\mathcal{A}}(\lambda) = \mathfrak{D}(\lambda\star_M\mathcal{I} - \mathcal{A}) = (\lambda - \lambda_1)\circ_M(\lambda - \lambda_2)\circ_M\dots\circ_M(\lambda - \lambda_k).
    \end{displaymath}
    Now here, substituting $\lambda = \mathbint{0}$, we get $\mathfrak{D}(-\mathcal{A}) = \mathrm{sgn}(k)\circ_M\lambda_1\circ_M\lambda_2\circ_M\dots\circ_M\lambda_k$. As we know, $\mathfrak{D}(-\mathcal{A}) = \mathrm{sgn}(k)\circ_M\mathfrak{D}(\mathcal{A})$, hence, we get, $\mathfrak{D}(\mathcal{A}) = \lambda_1\circ_M\lambda_2\circ_M\dots\circ_M\lambda_k$. Now expanding $\mathfrak{D}(\lambda\star_M\mathcal{I} - \mathcal{A})$ we get, 
    \begin{displaymath}
        \mathfrak{D}(\lambda\star_M\mathcal{I} - \mathcal{A}) = (\lambda - \mathcal{A}^{[1,1]})\circ_M(\lambda - \mathcal{A}^{[2,2]})\circ_M\dots\circ_M(\lambda - \mathcal{A}^{[k,k]}) + q_{\mathcal{A}}(\lambda),
    \end{displaymath}
    where, $q_{\mathcal{A}}(\lambda)$ is polynomial of degree atmost $k-2$. Now evaluating coefficients of $\lambda^{k-1}$ in both \Cref{eq:determinant-lambda-multiplication} and \Cref{eq:trace-lambda-addition}, we get, $\mathcal{A}^{[1,1]} + \mathcal{A}^{[2,2]} + \dots + \mathcal{A}^{[k,k]} = \texttt{t-Tr}(\mathcal{A}) = \lambda_1 + \lambda_2 + \dots + \lambda_k$.
\end{proof}

\begin{theorem}[Tubal spectral zero-divisor theorem]
Let $\mathcal{A} \in \mathbb{K}_M^{k \times k}$ have principal tubal eigenvalues, $\Lambda = \{\lambda_1, \lambda_2, \dots, \lambda_k\} \subset \mathbb{K}_M$. The determinant $\mathfrak{D}(\mathcal{A})$ is a zero-divisor in $\mathbb{K}_M$ if and only if at least one tubal eigenvalue $\lambda_i \in \Lambda$ is a zero-divisor. 

Moreover, if $\mathbint{r} < \mathfrak{r}_t = k$, the annihilator ideal of the determinant is equal to the sum of the annihilator ideals of the zero-divisor eigenvalues,
\begin{displaymath}
    \mathrm{Ann}(\mathfrak{D}(\mathcal{A})) = \sum_{\lambda_i \in \Lambda_{zd}} \mathrm{Ann}(\lambda_i),
\end{displaymath}
where $\Lambda_{zd}$ is the subset of eigenvalues that are zero-divisors.
\end{theorem}
\begin{proof}
For any element $\mathbf{x} \in \mathbb{K}_M$, let $\mathbf{u_x}$ be the unique idempotent such that the principal ideal $\langle \mathbf{x} \rangle = \langle \mathbf{u_x} \rangle$. It follows that the annihilator of $\mathbf{x}$ is, $\mathrm{Ann}(\mathbf{x}) = \langle \mathbint{1} - \mathbf{u_x} \rangle$. Using \Cref{eq:determinant-lambda-multiplication}, $\mathfrak{D}(\mathcal{A}) = \prod_{i=1}^k \lambda_i$. Suppose all eigenvalues $\lambda_i$ are invertible. Then $\mathfrak{D}(\mathcal{A})$ is also an invertible, meaning $\mathrm{Ann}(\mathfrak{D}(\mathcal{A})) = \{\mathbint{0}\}$, so it is not a zero-divisor. Conversely, suppose at least one eigenvalue $\lambda_j$ is a zero-divisor. Then there exists some $\mathbf{z} \neq \mathbint{0}$ such that $\mathbf{z} \circ_M \lambda_j = \mathbint{0}$, we get,
\begin{displaymath}
    \mathbf{z} \circ_M \mathfrak{D}(\mathcal{A}) = \mathbf{z} \circ_M \lambda_j \circ_M \left( \prod_{i \neq j} \lambda_i \right) = \mathbint{0} \circ_M \left( \prod_{i \neq j} \lambda_i \right) = \mathbint{0}.
\end{displaymath}
Since $\mathbf{z} \neq \mathbint{0}$, the determinant $\mathrm{Ann}(\mathfrak{D}(\mathcal{A})) = \{\mathbf{z}, \mathbint{0}\}\neq \{\mathbint{0}\}$ is annihilated and thus is a zero-divisor. Let $\mathbf{f}_i$ be the idempotent generating $\langle \lambda_i \rangle$. The ideal generated by the product of the eigenvalues is given by $\langle \mathfrak{D}(\mathcal{A}) \rangle = \langle \prod_{i=1}^k \mathbf{f}_i \rangle$. Hence, the annihilator of the determinant is
\begin{displaymath}
    \mathrm{Ann}(\mathfrak{D}(\mathcal{A})) = \left\langle \mathbint{1} - \prod_{i=1}^k \mathbf{f}_i \right\rangle.
\end{displaymath}

For any $j$, observe that $(\mathbint{1} - \mathbf{f}_j) \circ_M (\prod_{i=1}^k \mathbf{f}_i ) = \mathbint{0}$. This implies $(\mathbint{1} - \mathbf{f}_j) \in \mathrm{Ann}(\mathfrak{D}(\mathcal{A}))$, meaning $\sum_{i=1}^k \langle \mathbint{1} - \mathbf{f}_i \rangle \subseteq \mathrm{Ann}(\mathfrak{D}(\mathcal{A}))$. Conversely, we can express $(\mathbint{1} - \prod_{i=1}^k \mathbf{f}_i)$ algebraically as a linear combination of the individual $(\mathbint{1} - \mathbf{f}_i)$ terms. For $k=2$,
\begin{displaymath}
    \mathbint{1} - \mathbf{f}_1 \circ_M \mathbf{f}_2 = (\mathbint{1} - \mathbf{f}_1) + \mathbf{f}_1 \circ_M (\mathbint{1} - \mathbf{f}_2).
\end{displaymath}
By induction, $(\mathbint{1} - \prod_{i=1}^k \mathbf{f}_i) \in \sum_{i=1}^k \langle \mathbint{1} - \mathbf{f}_i \rangle$. 
Thus, $\mathrm{Ann}(\mathfrak{D}(\mathcal{A})) = \sum_{i=1}^k \mathrm{Ann}(\lambda_i)$. Since $\mathrm{Ann}(\lambda_i) = \{\mathbint{0}\}$ for any invertible $\lambda_i$, the sum reduces strictly to the subset of zero-divisor eigenvalues, $\Lambda_{zd}$, as $\mathrm{Ann}(\mathfrak{D}(\mathcal{A})) = \sum_{\lambda_i \in \Lambda_{zd}} \mathrm{Ann}(\lambda_i)$.
\end{proof}

\begin{example}
To demonstrate \Cref{def:eigval-and-eigvect}, we consider the same $4$ cases as \Cref{example:determinant}.

\begin{enumerate}
    \item After transformation, the tubal eigenvalues of the tensor and corresponding tubal eigenvectors are $\hat{\lambda}_1 = [3, 0, 0, 0]^\top$ and $\hat{\lambda}_2 = [0, 0, 0, 0]^\top$, which implies $\lambda_1 = [1.5, 1.5, 1.5, 1.5]^\top$ and $\lambda_2 = \mathbint{0}$.
    \begin{displaymath}
        \vec{\mathcal{Q}}_1 = \begin{bmatrix} [\frac{2}{\sqrt{5}}, 0, 0, 0]^\top \\ [\frac{2}{\sqrt{5}}, 0, 0, 0]^\top \end{bmatrix}, \quad
        \vec{\mathcal{Q}}_2 = \begin{bmatrix} [\frac{2}{\sqrt{2}}, 0, 0, 0]^\top \\ [-\frac{2}{\sqrt{2}}, 0, 0, 0]^\top \end{bmatrix}.
    \end{displaymath}
    Now here, it is clear that $\lambda_1\circ_M\lambda_2 = \mathbint{0} = \mathfrak{D}(\mathcal{A})$, that is, zero determinant case and non-invertible.

    \item After transformation, the tubal eigenvalues of the tensor and corresponding tubal eigenvectors are $\hat{\lambda}_1 = [2, 1, 1, 1]^\top$ and $\hat{\lambda}_2 = [2, 1, 1, 1]^\top$, which implies $\lambda_1 = [2.5, 0.5, 0.5, 0.5]^\top$ and $\lambda_2 = [2.5, 0.5, 0.5, 0.5]^\top$. Here we have, $\vec{\mathcal{Q}}_1 = \vec{\mathcal{I}}_1$ and $\vec{\mathcal{Q}}_2 = \vec{\mathcal{I}}_2$. Now here, it is clear that $\lambda_1\circ_M\lambda_2 = [3.5,1.5,1.5,1.5]^\top = \mathfrak{D}(\mathcal{A})$, that is, non-zero determinant case and invertible.

    \item After transformation, the tubal eigenvalues of the tensor and corresponding tubal eigenvectors are $\hat{\lambda}_1 = [2, 2, 0, 0]^\top$ and $\hat{\lambda}_2 = [1, 1, 0, 0]^\top$, which implies $\lambda_1 = [2, 0, 2, 0]^\top$ and $\lambda_2 = [1, 0, 1, 0]^\top$. Here we have, $\vec{\mathcal{Q}}_1 = \vec{\mathcal{I}}_1$ and $\vec{\mathcal{Q}}_2 = \vec{\mathcal{I}}_2$. Now here, it is clear that $\lambda_1\circ_M\lambda_2 = [2,0,2,0]^\top = \mathfrak{D}(\mathcal{A})$, that is, non-zero determinant case and non-invertible.

    \item After transformation, the tubal eigenvalues of the tensor and corresponding tubal eigenvectors are $\hat{\lambda}_1 = [1, 1, 1, 1]^\top$ and $\hat{\lambda}_2 = [0, 0, 0, 0]^\top$, which implies $\lambda_1 = [2, 0, 0, 0]^\top = \mathbint{1}$ and $\lambda_2 = \mathbint{0}$. 
    \begin{displaymath}
        \vec{\mathcal{Q}}_1 = \begin{bmatrix} [\sqrt{2}, 0, 0, 0]^\top \\ [\sqrt{2}, 0, 0, 0]^\top \end{bmatrix}, \quad
        \vec{\mathcal{Q}}_2 = \begin{bmatrix} [\sqrt{2}, 0, 0, 0]^\top \\ [-\sqrt{2}, 0, 0, 0]^\top \end{bmatrix}.
    \end{displaymath}
    Now here, it is clear that $\lambda_1\circ_M\lambda_2 = \mathbint{0} = \mathfrak{D}(\mathcal{A})$, that is, zero determinant case and non-invertible.
\end{enumerate}    
\end{example}

\begin{lemma}[Existence of regular orthonormal basis]\label{lem:ortho-basis-existence}
    Let $\mathcal{A} \in \mathbb{K}_M^{k \times k}$ be a symmetric tensor. There exists a set of $k$ regular, mutually orthogonal unit tubal vectors $\{\vec{\mathcal{Q}}_1, \dots, \vec{\mathcal{Q}}_k\}$ in $\mathbb{K}_M^{k \times 1}$ such that they form a free basis for the module $\mathbb{K}_M^k$, and each $\vec{\mathcal{Q}}_i$ is a regular tubal eigenvector of $\mathcal{A}$.
\end{lemma}
\begin{proof}
    From \Cref{def:funda-idem-tubal-scalar}, we know that $\mathbb{K}_M$ contains a complete set of $N$ fundamental orthogonal idempotents $\{\mathbf{e}_1, \dots, \mathbf{e}_N\}$ satisfying $\mathbf{e}_i \circ_M \mathbf{e}_j = \mathbint{0}$ for $i \neq j$ and $\sum_{i=1}^N \mathbf{e}_i = \mathbint{1}$. For any idempotent $\mathbf{e}_i$, the principal ideal $\langle \mathbf{e}_i \rangle$ behaves as a field $\mathbb{F}_m$ and hence, the projected tensor $\mathcal{A}_m = \mathbf{e}_m \star_M \mathcal{A}$ acts as a standard symmetric tensor over $\mathbb{F}_m^k$. Hence, by the classical Spectral Theorem over fields, $\mathcal{A}_m$ admits a complete orthonormal basis of eigenvectors $\{\vec{\mathcal{Q}}_1^{(m)}, \dots, \vec{\mathcal{Q}}_k^{(m)}\}$ in $\mathbb{F}_m^k$. The orthonormality over the local ideal is given by $(\vec{\mathcal{Q}}_i^{(m)})^\top \star_M \vec{\mathcal{Q}}_j^{(m)} = \delta_{i,j} \circ_M \mathbf{e}_m$. Constructing the global eigenvectors for the tensor $\mathcal{A}$, that is, $\vec{\mathcal{Q}}_i = \sum_{m=1}^N \vec{\mathcal{Q}}_i^{(m)}$, and verifying global orthonormality by cancelling cross-terms, we get, 
    \begin{displaymath}
        \vec{\mathcal{Q}}_i^\top \star_M \vec{\mathcal{Q}}_j = \sum_{m=1}^N (\vec{\mathcal{Q}}_i^{(m)})^\top \star_M \vec{\mathcal{Q}}_j^{(m)} = \sum_{m=1}^N \delta_{i,j} \circ_M \mathbf{e}_m = \delta_{i,j} \circ_M \mathbint{1}.
    \end{displaymath}
    Hence, $\vec{\mathcal{Q}}_i^\top \star_M \vec{\mathcal{Q}}_i = \mathbint{1}$ and $\vec{\mathcal{Q}}_i^\top \star_M \vec{\mathcal{Q}}_j = \mathbint{0}$ for $i \neq j$. Suppose there exists some tubal scalar $\mathbf{c} \in \mathbb{K}_M$ such that $\mathbf{c} \circ_M \vec{\mathcal{Q}}_i = \vec{\mathcal{O}}$. Multiplying both sides by the idempotent $\mathbf{e}_m$ gives, $\mathbf{c} \circ_M \vec{\mathcal{Q}}_i^{(m)} = \vec{\mathcal{O}}$. Since $\vec{\mathcal{Q}}_i^{(m)}$ is a unit vector in $\mathbb{F}_m^k$, its components cannot all be zero-divisors with respect to $\mathbf{e}_m$, which follows that $\mathbf{c} \circ_M \mathbf{e}_m = \mathbint{0}$ for every $m \in \{1, \dots, N\}$. Summing over all idempotents,
    \begin{displaymath}
        \sum_{m=1}^N (\mathbf{c} \circ_M \mathbf{e}_m) = \mathbf{c} \circ_M \left( \sum_{m=1}^N \mathbf{e}_m \right) = \mathbf{c} \circ_M \mathbint{1} = \mathbf{c} = \mathbint{0}.
    \end{displaymath}
    Since $\mathrm{Ann}(\vec{\mathcal{Q}}_i) = \{\mathbint{0}\}$, the tubal vector $\vec{\mathcal{Q}}_i$ is regular. Hence, $\{\vec{\mathcal{Q}}_1, \dots, \vec{\mathcal{Q}}_k\}$ constitutes a set of $k$ regular, mutually orthogonal tubal vectors in $\mathbb{K}_M^{k \times 1}$, spanning a free submodule of rank $k$. As this submodule contains $k$ linearly independent regular vectors, it strictly spans $\mathbb{K}_M^k$ entirely, forming a free orthonormal basis.
\end{proof}

\begin{theorem}[Tubal Gram-Schmidt orthogonalization]\label{thrm:tubal-gram-schmidt}
    Let $\mathscr{W}$ be a finitely generated projective submodule of $\mathbb{K}_M^{k \times 1}$ spanned by a finite set of tubal vectors $\mathscr{S} = \{\vec{\mathcal{X}}_1, \vec{\mathcal{X}}_2, \dots, \vec{\mathcal{X}}_m\}$. There exists a set of mutually orthogonal tubal vectors $\mathscr{Q} = \{\vec{\mathcal{Q}}_1, \vec{\mathcal{Q}}_2, \dots, \vec{\mathcal{Q}}_m\}$ that generates the same submodule $\mathscr{W}$, satisfying $\vec{\mathcal{Q}}_i^\top \star_M \vec{\mathcal{Q}}_j = \mathbint{0}$ for all $i \neq j$.
\end{theorem}
\begin{proof}
    For the base case, define $\vec{\mathcal{Q}}_1 = \vec{\mathcal{X}}_1$. For the inductive step, assume we have constructed a set of mutually orthogonal tubal vectors $\{\vec{\mathcal{Q}}_1, \dots, \vec{\mathcal{Q}}_{p}\}$ that spans the same submodule as $\{\vec{\mathcal{X}}_1, \dots, \vec{\mathcal{X}}_{p}\}$. Define the orthogonal projection of a tubal vector $\vec{\mathcal{X}}$ onto a tubal vector $\vec{\mathcal{Q}}$ using the tubal Moore-Penrose pseudoinverse as,
    \begin{displaymath}
        \texttt{proj}_{\vec{\mathcal{Q}}}(\vec{\mathcal{X}}) = \vec{\mathcal{Q}} \circ_M (\vec{\mathcal{Q}}^\top \star_M \vec{\mathcal{Q}})^\dagger \circ_M (\vec{\mathcal{Q}}^\top \star_M \vec{\mathcal{X}}).
    \end{displaymath}
    Now, constructing the $(p+1)$th tubal vector as,
    \begin{displaymath}
        \vec{\mathcal{Q}}_{p+1} = \vec{\mathcal{X}}_{p+1} - \sum_{i=1}^{p} \texttt{proj}_{\vec{\mathcal{Q}}_i}(\vec{\mathcal{X}}_{p+1}).
    \end{displaymath}
    Taking the tubal inner product of $\vec{\mathcal{Q}}_{p+1}$ with any previous tubal vector $\vec{\mathcal{Q}}_j$ where $1 \leq j \leq p$ and utlizing inductive hypothesis $\vec{\mathcal{Q}}_j^\top \star_M \vec{\mathcal{Q}}_i = \mathbint{0}$ for $i \neq j$, we get,
    \begin{displaymath}
        \vec{\mathcal{Q}}_j^\top \star_M \vec{\mathcal{Q}}_{p+1} =  \vec{\mathcal{Q}}_j^\top \star_M \vec{\mathcal{X}}_{p+1} - (\vec{\mathcal{Q}}_j^\top \star_M \vec{\mathcal{Q}}_j) \circ_M (\vec{\mathcal{Q}}_j^\top \star_M \vec{\mathcal{Q}}_j)^\dagger \circ_M (\vec{\mathcal{Q}}_j^\top \star_M \vec{\mathcal{X}}_{p+1}).
    \end{displaymath}
    Assume the idempotent tubal scalar $\mathbf{e}_j = (\vec{\mathcal{Q}}_j^\top \star_M \vec{\mathcal{Q}}_j) \circ_M (\vec{\mathcal{Q}}_j^\top \star_M \vec{\mathcal{Q}}_j)^\dagger$. Hence, we get, $\mathbf{e}_j \circ_M (\vec{\mathcal{Q}}_j^\top \star_M \vec{\mathcal{X}}_{p+1}) = \vec{\mathcal{Q}}_j^\top \star_M \vec{\mathcal{X}}_{p+1}$ as $\mathbf{e}_j$ supports $\vec{\mathcal{Q}}$. Now, after substitution, we get,
    \begin{displaymath}
        \vec{\mathcal{Q}}_j^\top \star_M \vec{\mathcal{Q}}_{p+1} = \vec{\mathcal{Q}}_j^\top \star_M \vec{\mathcal{X}}_{p+1} - \vec{\mathcal{Q}}_j^\top \star_M \vec{\mathcal{X}}_{p+1} = \mathbint{0}.
    \end{displaymath}
    By mathematical induction, the set $\mathscr{Q}$ is mutually orthogonal. Furthermore, because each $\vec{\mathcal{Q}}_i$ is formed entirely by finite linear combinations of $\mathscr{S}$, the tubal span is preserved, that is, $\texttt{span}_{\mathbb{K}_M}(\mathscr{Q}) = \texttt{span}_{\mathbb{K}_M}(\mathscr{S}) = \mathscr{W}$.
    Also, normalizing each resulting tubal vector in $\mathscr{Q}$ results in orthonormal basis generated by, 
    \begin{displaymath}
        \vec{\mathcal{U}}_i = 
        \vec{\mathcal{Q}}_i \circ_M \Phi^{-1}\left(\left(\sqrt{\Phi(\vec{\mathcal{Q}}_i^\top \star_M \vec{\mathcal{Q}}_i)}\right)^\dagger\right).
    \end{displaymath} 
\end{proof}

\begin{theorem}[Tensor-extended spectral theorem]\label{thrm:tensor-extended-spectral-thrm}
    Let $\mathcal{A} \in \mathbb{K}_M^{k \times k}$ be a symmetric tensor, such that $\mathcal{A}^\top = \mathcal{A}$. There exists an orthogonal tensor $\mathcal{Q} \in \mathbb{K}_M^{k \times k}$ (satisfying $\mathcal{Q}^\top \star_M \mathcal{Q} = \mathcal{Q} \star_M \mathcal{Q}^\top = \mathcal{I}$) and a $f$-diagonal tensor $\mathcal{D} \in \mathbb{K}_M^{k \times k}$ such that,
    \begin{displaymath}
    \mathcal{A} = \mathcal{Q} \star_M \mathcal{D} \star_M \mathcal{Q}^\top,
    \end{displaymath}
    where, $\mathcal{D}^{[i,i]}$ is the $i$th principal tubal eigenvalue of $\mathcal{A}$, and the columns of $\mathcal{Q}$ form a free orthonormal basis of regular tubal eigenvectors.
\end{theorem}
\begin{proof}
    Let $p_{\mathcal{A}}(\lambda) = \mathfrak{D}(\lambda \star_M \mathcal{I} - \mathcal{A})$ be the tubal characteristic polynomial. Because $\mathcal{A}$ is symmetric, eigenvectors corresponding to distinct principal tubal eigenvalues $\lambda_i \neq \lambda_j$ satisfy $(\lambda_i - \lambda_j) \circ_M (\vec{\mathcal{X}}_i^\top \star_M \vec{\mathcal{X}}_j) = \mathbint{0}$~{\cite[Proposition 4.3]{kernfeld2015tensor}}. Because $\mathbb{K}_M$ is a Von Neumann regular, assume the principal ideal generated by the difference of the eigenvalues can be represented by a unique idempotent $\mathbf{u} \in \mathbb{K}_M$, such that $\langle \lambda_i - \lambda_j \rangle = \langle \mathbf{u} \rangle$ due to which, $\mathrm{Ann}(\lambda_i - \lambda_j) = \langle\mathbint{1} - \mathbf{u}\rangle$, and hence, $\vec{\mathcal{X}}_i^\top \star_M \vec{\mathcal{X}}_j \in \langle \mathbint{1} - \mathbf{u} \rangle$. Now, on the ideal $\langle \mathbf{u} \rangle$, the eigenvalues are strictly distinct, and the inner product is perfectly $\mathbint{0}$. On the ideal $\langle \mathbint{1} - \mathbf{u} \rangle$, the inner product may be non-zero, but on this specific ideal, the difference $(\lambda_i - \lambda_j)$ is annihilated, meaning $\lambda_i = \lambda_j$. Consequently, if the inner product is non-zero, it only occurs in the sub-ideal where $\vec{\mathcal{X}}_i$ and $\vec{\mathcal{X}}_j$ share the exact same eigenvalue, meaning they belong to the same tubal eigensubmodule. Applying \Cref{thrm:tubal-gram-schmidt} to the basis tubal vectors spanning this shared eigensubmodule, which orthogonalizes the vectors on $\langle \mathbint{1} - \mathbf{u} \rangle$ without altering the orthogonality on $\langle \mathbf{u} \rangle$. Normalizing these orthogonalized eigenvectors such that $\vec{\mathcal{Q}}_i^\top \star_M \vec{\mathcal{Q}}_i = \mathbint{1}$ yields the columns of the orthogonal tensor $\mathcal{Q}$. The equivalence is given as, $\mathcal{A} \star_M \mathcal{Q} = \mathcal{Q} \star_M \mathcal{D}$. Now, because $\mathcal{Q} \star_M \mathcal{Q}^\top = \mathcal{I}$, right-multiplying by $\mathcal{Q}^\top$ yields, $\mathcal{A} = \mathcal{Q} \star_M \mathcal{D} \star_M \mathcal{Q}^\top$.
    
    This means regular tubal eigenvectors corresponding to distinct principal tubal eigenvalues are orthogonal. Normalizing these orthogonal eigenvectors, $\vec{\mathcal{X}}^\top \star_M \vec{\mathcal{X}} = \mathbint{1}$, yields the columns of the orthogonal tensor $\mathcal{Q}$. The equivalence $\mathcal{A} \star_M \mathcal{Q} = \mathcal{Q} \star_M \Lambda$. Because $\mathcal{Q}\star_M\mathcal{Q}^\top = \mathcal{I}$, right-multiplying by $\mathcal{Q}^\top$ yields $\mathcal{A} = \mathcal{Q} \star_M \mathcal{D} \star_M \mathcal{Q}^\top$.
\end{proof}

\begin{theorem}[Tensor-extended singular value decomposition]\label{thrm:tensor-extended-svd}
    For any tensor $\mathcal{A} \in \mathbb{K}_M^{m \times n}$, there exist orthogonal tensors $\mathcal{U} \in \mathbb{K}_M^{m \times m}$ and $\mathcal{V} \in \mathbb{K}_M^{n \times n}$, and an $f$-diagonal tensor $\mathcal{S} \in \mathbb{K}_M^{m \times n}$ such that,
    \begin{displaymath}
    \mathcal{A} = \mathcal{U} \star_M \mathcal{S} \star_M \mathcal{V}^\top = \sum_{i=1}^{\min\{m, n\}} \mathcal{U}^{[:,i]} \star_M \mathcal{S}^{[i,i]}\star_M \left(\mathcal{V}_i^\top\right)^{[i,:]},
    \end{displaymath}
    where, $\mathcal{S}^{[i,i]}$ is the $i$th tubal singular values of $\mathcal{A}$, satisfying $\mathcal{S}^{[i,i]} \succeq \mathbint{0}$.
\end{theorem}
\begin{proof}
    Construct the symmetric tensor $\mathcal{A}^\top \star_M \mathcal{A} \in \mathbb{K}_M^{n \times n}$. By \Cref{thrm:tensor-extended-spectral-thrm}, there exists an orthogonal tensor $\mathcal{V} = [\vec{\mathcal{V}}_1, \dots, \vec{\mathcal{V}}_n] \in \mathbb{K}_M^{n \times n}$ such that $\mathcal{V}^\top \star_M (\mathcal{A}^\top \star_M \mathcal{A}) \star_M \mathcal{V} = \Lambda$, where $\Lambda$ is an $f$-diagonal tensor of principal tubal eigenvalues $\lambda_i \succeq \mathbint{0}$ such that, $\lambda_1 \succeq \lambda_2 \succeq \dots \succeq \lambda_n$. Now, define the $f$-diagonal singular value tensor $\mathcal{S} \in \mathbb{K}_M^{m \times n}$ where the diagonal elements are the tubal square roots, that is, $\sigma_i = \mathcal{S}^{[i,i]} = \Phi^{-1}(\Phi(\sqrt{\lambda_i}))$. Assume a tubal vector $\vec{\mathcal{X}}_i = \mathcal{A} \star_M \vec{\mathcal{V}}_i \in \mathbb{K}_M^{m \times 1}$, then we have,
    \begin{displaymath}
        \vec{\mathcal{X}}_i^\top \star_M \vec{\mathcal{X}}_j = \vec{\mathcal{V}}_i^\top \star_M (\mathcal{A}^\top \star_M \mathcal{A}) \star_M \vec{\mathcal{V}}_j = \delta_{i,j} \circ_M \sigma_i^2.
    \end{displaymath}
    Assume a unique idempotent for each singular value as $\mathbf{e}_i = \sigma_i \circ_M \sigma_i^\dagger$ and tubal vector $\vec{\mathcal{W}}_i = \sigma_i^\dagger \circ_M \vec{\mathcal{X}}_i$. As, $\vec{\mathcal{W}}_i^\top \star_M \vec{\mathcal{W}}_j = \delta_{i,j} \circ_M \mathbf{e}_i$ implies the tubal vectors are orthonormal on $\langle \mathbf{e}_i \rangle$. Because these vectors may contain zero-divisors, they do not strictly span the free module $\mathbb{K}_M^m$. To construct a complete globally regular orthonormal basis, form an augmented set $\mathscr{S}_{\mathrm{full}}$ using \Cref{def:tubal-natural-basis}, such that $\mathscr{S}_{\mathrm{full}} = \{\vec{\mathcal{W}}_1, \dots, \vec{\mathcal{W}}_n, \vec{\mathcal{I}}_1, \dots, \vec{\mathcal{I}}_m\}$. Applying \Cref{thrm:tubal-gram-schmidt} to $\mathscr{S}_{\mathrm{full}}$ will sequentially accept the $\vec{\mathcal{W}}_i$ vectors on the local ideals, which guarantees the generation of a complete set of $m$ globally regular, mutually orthogonal unit tubal vectors $\{\vec{\mathcal{U}}_1, \dots, \vec{\mathcal{U}}_m\}$. Using properties of \Cref{thrm:tubal-gram-schmidt}, $\mathbf{e}_i \circ_M \vec{\mathcal{U}}_i = \vec{\mathcal{W}}_i$ which implies $\mathbf{e}_i \circ_M \vec{\mathcal{U}}_i = \sigma_i^\dagger \circ_M (\mathcal{A} \star_M \vec{\mathcal{V}}_i)$. Hence, we get, $\sigma_i \circ_M \vec{\mathcal{U}}_i = \mathbf{e}_i \circ_M (\mathcal{A} \star_M \vec{\mathcal{V}}_i)$. Crucially, on the complementary ideal $\langle \mathbint{1} - \mathbf{e}_i \rangle$, the singular value strictly evaluates to $\mathbint{0}$. Therefore, the relationship $\sigma_i \circ_M \vec{\mathcal{U}}_i = \mathcal{A} \star_M \vec{\mathcal{V}}_i$ holds as on that ideal. Thus, the relation holds unconditionally across all ideals of $\mathbb{K}_M$, that yeilds, $\mathcal{A} \star_M \vec{\mathcal{V}}_i = \sigma_i \circ_M \vec{\mathcal{U}}_i$. As $\mathcal{V}$ is a globally orthogonal tensor ($\mathcal{V} \star_M \mathcal{V}^\top = \mathcal{I}$), we get $\mathcal{A} = \mathcal{U} \star_M \mathcal{S} \star_M\mathcal{V}^\top$.
\end{proof}

\begin{figure}[htbp]
    \centering
    \includegraphics[width=0.8\linewidth]{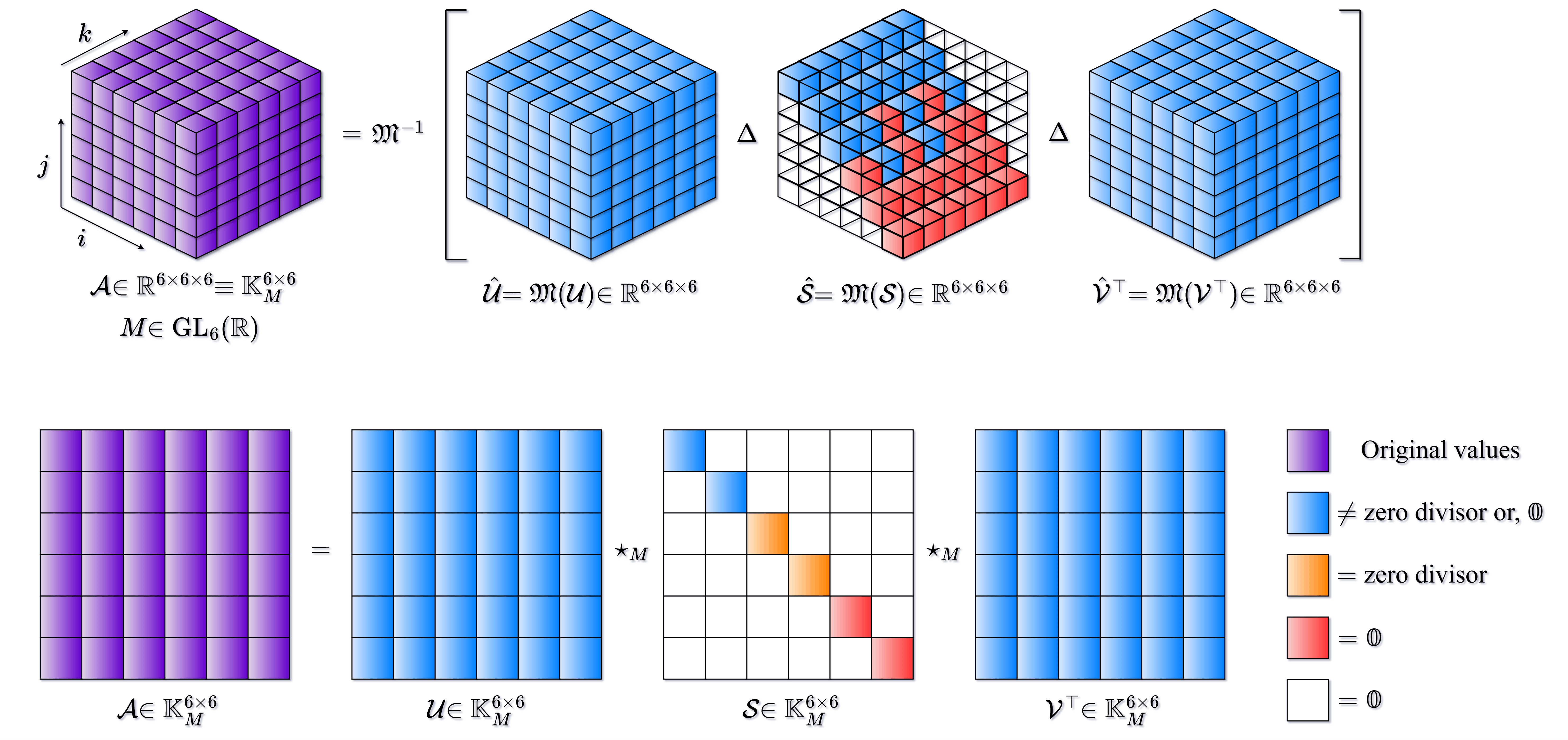}
    \caption{Visualization of Tensor-extended SVD of tensor $\mathcal{A}$ according to \Cref{thrm:tensor-extended-svd}.}
    \label{fig:svd-cases}
\end{figure}

\begin{example}\label{example:svd}
    To explicitly demonstrate \Cref{thrm:tensor-extended-svd} and \Cref{thrm:tubal-gram-schmidt}, we compute the tubal SVD for the three tensor cases defined in \Cref{example:rank-nullity}.

\begin{enumerate}
    \item For the invertible singular values case, the singular values are $\hat{\sigma}_1 = \hat{\sigma}_2 = \hat{\sigma}_3 = [1, 1, 1, 1]^\top$, which implies $\sigma_j = \mathbint{1}$ for $j = 1,2,3$. It is clear that here, $\mathcal{U} = \mathcal{S} = \mathcal{V} = \mathcal{I}$.

    \item For the zero singular value case, the singular values are $\hat{\sigma}_1 = [\sqrt{3}, \sqrt{3}, \sqrt{3}, \sqrt{3}]^\top$, $\hat{\sigma}_2 = [1, 1, 1, 1]^\top$ and $\hat{\sigma}_3 = [0, 0, 0, 0]^\top$, which implies $\sigma_1 = [2\sqrt{3}, 0, 0, 0]^\top$, $\sigma_2 = [2, 0, 0, 0]^\top = \mathbint{1}$ and $\sigma_3 = \mathbint{0}$. The right singular vectors map to $\hat{\mathcal{V}}[:,1,i] = [\frac{1}{\sqrt{6}}, \frac{1}{\sqrt{6}}, \frac{2}{\sqrt{6}}]^\top$, $\hat{\mathcal{V}}[:,2,i] = [\frac{-1}{\sqrt{2}}, \frac{1}{\sqrt{2}}, 0]^\top$, $\hat{\mathcal{V}}[:,3,i] = [\frac{1}{\sqrt{3}}, \frac{1}{\sqrt{3}}, \frac{-1}{\sqrt{3}}]^\top$ for $i = 1,2,3,4$. Constructing $\vec{\mathcal{W}}_j = \sigma_j^\dagger \circ_M (\mathcal{A} \star_M \vec{\mathcal{V}}_j)$, we get, $\hat{\mathcal{W}}[:,1,i] = [\frac{1}{\sqrt{2}}, \frac{1}{\sqrt{2}}, 0]^\top,$ $\hat{\mathcal{W}}[:,2,i] = [\frac{-1}{\sqrt{2}}, \frac{1}{\sqrt{2}}, 0]^\top$ for $i = 1,2,3,4$. As $\sigma_3 = \mathbint{0}$, $\vec{\mathcal{W}}_3$ collapses to $\vec{\mathcal{O}}$. Applying the \Cref{thrm:tubal-gram-schmidt} to project $\vec{\mathcal{I}}_3$ orthogonal to $\vec{\mathcal{W}}_1$ and $\vec{\mathcal{W}}_2$, we get $\hat{\mathcal{U}}[:,3,i] = [0, 0, 1]^\top$. Finally, 
    \begin{displaymath}
     \mathcal{U}[:,:,1] = \begin{bmatrix} \sqrt{2} & -\sqrt{2} & 0 \\ \sqrt{2} & \sqrt{2} & 0 \\ 0 & 0 & 2 \end{bmatrix},
        ~~\mathcal{S}[:,:,1] = \begin{bmatrix} 2\sqrt{3} & 0 & 0 \\ 0 & 2 & 0 \\ 0 & 0 & 0 \end{bmatrix}, 
        ~~\mathcal{V}[:,:,1] = \begin{bmatrix} \frac{2}{\sqrt{6}} & -\sqrt{2} & \frac{2}{\sqrt{3}} \\ \frac{2}{\sqrt{6}} & \sqrt{2} & \frac{2}{\sqrt{3}} \\ \frac{4}{\sqrt{6}} & 0 & \frac{-2}{\sqrt{3}} \end{bmatrix},
    \end{displaymath}
    \begin{displaymath}
        \mathcal{U}^{[:,:,i]} = \mathcal{S}^{[:,:,i]} = \mathcal{V}^{[:,:,i]} = \begin{bmatrix}
        0&0&0\\0&0&0\\0&0&0
    \end{bmatrix}, \quad \text{where, } i=2,3,4.
    \end{displaymath}

    \item For the zero-divisor singular value case, the singular values are $\hat{\sigma}_1 = [1, \sqrt{3}, \sqrt{3}, \sqrt{3}]^\top$, $\hat{\sigma}_2 = [1, 1, 1, 1]^\top$ and $\hat{\sigma}_3 = [1, 0, 0, 0]^\top$, which implies $\sigma_1 = [0.5(1+3\sqrt{3}), 0.5(1-\sqrt{3}), 0.5(1-\sqrt{3}), 0.5(1-\sqrt{3})]^\top$, $\sigma_2 = \mathbint{1}$ and $\sigma_3 = [0.5, 0.5, 0.5, 0.5]^\top$. The right singular vectors map to
    \begin{displaymath}
        \hat{\mathcal{V}}[:,:,1] = \begin{bmatrix} 1 & 0 & 0 \\ 0 & 1 & 0 \\ 0 & 0 & 1 \end{bmatrix},~~\hat{\mathcal{V}}[:,:,2] = \hat{\mathcal{V}}[:,:,3] = \hat{\mathcal{V}}[:,:,4] = \begin{bmatrix} \frac{1}{\sqrt{6}} & \frac{-1}{\sqrt{2}} & \frac{1}{\sqrt{3}} \\ \frac{1}{\sqrt{6}} & \frac{1}{\sqrt{2}} & \frac{1}{\sqrt{3}} \\ \frac{2}{\sqrt{6}} & 0 & \frac{-1}{\sqrt{3}} \end{bmatrix}.
    \end{displaymath}
    Constructing $\vec{\mathcal{W}}_j = \sigma_j^\dagger \circ_M (\mathcal{A} \star_M \vec{\mathcal{V}}_j)$, we get,
    \begin{align*}
        \hat{\mathcal{W}}[:,1,i] = \begin{cases} [1, 0, 0]^\top & i=1, \\ [\frac{1}{\sqrt{2}}, \frac{1}{\sqrt{2}}, 0]^\top & i = 2,3,4. \end{cases} \quad  
        \hat{\mathcal{W}}[:,2,i] = \begin{cases} [0, 1, 0]^\top & i=1, \\ [\frac{-1}{\sqrt{2}}, \frac{1}{\sqrt{2}}, 0]^\top & i=2,3,4.\end{cases}
    \end{align*}
    As $\sigma_3$ is zero-divisor, $\vec{\mathcal{W}}_3$ collapses to $\vec{\mathcal{O}}$ on sub-ideals. Applying the \Cref{thrm:tubal-gram-schmidt} to project $\vec{\mathcal{I}}_3$ orthogonal to $\vec{\mathcal{W}}_1$ and $\vec{\mathcal{W}}_2$, we get $\hat{\mathcal{U}}[:,3,i] = [0, 0, 1]^\top$. Finally,
    \begin{displaymath}
        \mathcal{U}[:,:,1] =
        \begin{bmatrix}
        \frac{1}{2} + \frac{3}{2\sqrt{2}} & -\frac{3}{2\sqrt{2}} & 0 \\
        \frac{3}{2\sqrt{2}} & \frac{1}{2} + \frac{3}{2\sqrt{2}} & 0 \\
        0 & 0 & 2
        \end{bmatrix},~~\mathcal{U}^{[:,:,i]} =
        \begin{bmatrix}
        \frac{1}{2} - \frac{1}{2\sqrt{2}} & \frac{1}{2\sqrt{2}} & 0 \\
        -\frac{1}{2\sqrt{2}} & \frac{1}{2} - \frac{1}{2\sqrt{2}} & 0 \\
        0 & 0 & 0
        \end{bmatrix},
    \end{displaymath}
    \begin{displaymath}
        \mathcal{S}[:,:,1] =
        \begin{bmatrix}
        \frac{1}{2} + \frac{3\sqrt{3}}{2} & 0 & 0 \\
        0 & 2 & 0 \\
        0 & 0 & \frac{1}{2}
        \end{bmatrix},~~\mathcal{S}^{[:,:,i]} =
        \begin{bmatrix}
        \frac{1}{2} - \frac{\sqrt{3}}{2} & 0 & 0 \\
        0 & 0 & 0 \\
        0 & 0 & \frac{1}{2}
        \end{bmatrix},
    \end{displaymath}
    
    \begin{displaymath}
        \mathcal{V}[:,:,1] =
        \begin{bmatrix}
        \frac{1}{2} + \frac{3}{2\sqrt{6}} & -\frac{3}{2\sqrt{2}} & \frac{3}{2\sqrt{3}} \\
        \frac{3}{2\sqrt{6}} & \frac{1}{2} + \frac{3}{2\sqrt{2}} & \frac{3}{2\sqrt{3}} \\
        \frac{3}{\sqrt{6}} & 0 & \frac{1}{2} - \frac{3}{2\sqrt{3}}
        \end{bmatrix},~~\mathcal{V}^{[:,:,i]} =
        \begin{bmatrix}
        \frac{1}{2} - \frac{1}{2\sqrt{6}} & \frac{1}{2\sqrt{2}} & -\frac{1}{2\sqrt{3}} \\
        -\frac{1}{2\sqrt{6}} & \frac{1}{2} - \frac{1}{2\sqrt{2}} & -\frac{1}{2\sqrt{3}} \\
        -\frac{1}{\sqrt{6}} & 0 & \frac{1}{2} + \frac{1}{2\sqrt{3}}
        \end{bmatrix},
    \end{displaymath}
    for $i = 2,3,4$.
\end{enumerate}
\end{example}

\begin{theorem}[Tubal Rayleigh Quotient]\label{thrm:Tubal-Rayleigh-Quotient}
Let $\mathcal{A} \in \mathbb{K}_M^{k \times k}$ be a symmetric tensor. The $i$th principal tubal eigenvalue $\lambda_i$, such that $\lambda_1 \succeq \lambda_2 \succeq \dots \succeq \lambda_k$, can be characterized by the tubal Rayleigh quotient,
\begin{equation}
    \lambda_i(\mathcal{A}) = \max_{\mathrm{Dim}_{\mathbb{K}_M}(\mathscr{W})=i\ast \mathbint{1}} \min_{\substack{\vec{\mathcal{X}} \in \mathscr{W} \\ \mathrm{Ann}(\vec{\mathcal{X}}) = \{\mathbint{0}\}\\ \|\vec{\mathcal{X}}\|_2 = \mathbint{1}}} (\vec{\mathcal{X}}^\top \star_M \vec{\mathcal{X}})^\dagger \circ_M (\vec{\mathcal{X}}^\top \star_M \mathcal{A} \star_M \vec{\mathcal{X}})
\end{equation}
where $\mathscr{W}$ is a finitely generated projective submodule of $\mathbb{K}_M^{k \times 1}$.
\end{theorem}
\begin{proof}
By \Cref{thrm:tensor-extended-spectral-thrm} and \Cref{lem:ortho-basis-existence}, there exists an orthogonal tensor $\mathcal{Q} = [\vec{\mathcal{Q}}_1, \dots, \vec{\mathcal{Q}}_k]$ spanning $\mathbb{K}_M^k$ such that $\mathcal{A} \star_M \vec{\mathcal{Q}}_j = \lambda_j \circ_M \vec{\mathcal{Q}}_j$. Any regular tubal vector $\vec{\mathcal{X}} \in \mathbb{K}_M^{k \times 1}$ can be expressed as, $\vec{\mathcal{X}} = \sum_{j=1}^k \mathbf{c}_j \circ_M \vec{\mathcal{Q}}_j$, where $\mathbf{c}_j \in \mathbb{K}_M$. Because $\vec{\mathcal{X}}$ is a regular vector, that is, $\mathrm{Ann}(\vec{\mathcal{X}}) = \{\mathbint{0}\}$, its squared tubal norm $\mathbf{n_x} = \vec{\mathcal{X}}^\top \star_M \vec{\mathcal{X}} = \sum_{j=1}^k \mathbf{c}_j^2$ is strictly invertible in $\mathbb{K}_M$. Also note that its tubal Moore-Penrose pseudoinverse coincides with its exact inverse, $\mathbf{n_x}^\dagger = \mathbf{n_x}^{-1}$. Hence, the tubal Rayleigh quotient $R(\vec{\mathcal{X}})$ evaluates to,
\begin{displaymath}
    R(\vec{\mathcal{X}}) = \mathbf{n_x}^{-1} \circ_M \left( \sum_{j=1}^k \lambda_j \circ_M \mathbf{c}_j^2 \right)
\end{displaymath}
Let $\mathscr{W}_i = \texttt{span}_{\mathbb{K}_M}\{\vec{\mathcal{Q}}_1, \dots, \vec{\mathcal{Q}}_i\}$ be an $i$-dimensional submodule. For any regular vector $\vec{\mathcal{X}} \in \mathscr{W}_i$, the coefficients $\mathbf{c}_j = \mathbint{0}$ for $j > i$. As $\lambda_j \succeq \lambda_i$ for all $j \leq i$, it follows that, $\sum_{j=1}^i \lambda_j \circ_M \mathbf{c}_j^2 \succeq \lambda_i \circ_M \sum_{j=1}^i \mathbf{c}_j^2$. Now, multiplying both sides by $\mathbf{n_x}^{-1}$, we get,
\begin{equation}\label{eq:rayleigh-proof-eq1}
    R(\vec{\mathcal{X}}) \succeq \lambda_i \circ_M \mathbf{n_x} \circ_M \mathbf{n_x}^{-1} = \lambda_i
\end{equation}
Now, let $\mathscr{V}$ be any arbitrary $i$-dimensional submodule, and let $\mathscr{U} = \texttt{span}_{\mathbb{K}_M}\{\vec{\mathcal{Q}}_i, \dots, \vec{\mathcal{Q}}_k\}$ be a submodule. Using \Cref{thrm:general-rank-nullity-thrm}, we get, $\mathrm{Dim}_{\mathbb{K}_M}(\mathscr{V} \cap \mathscr{U}) \succeq i + (k - i + 1) - k = \mathbint{1}$. Now, because the tubal dimension of the intersection is at least $\mathbint{1}$, there exists at least one globally regular tubal vector $\vec{\mathcal{Y}} \in \mathscr{V} \cap \mathscr{U}$ such that $\mathrm{Ann}(\vec{\mathcal{Y}}) = \{\mathbint{0}\}$. As $\vec{\mathcal{Y}} \in \mathscr{U}$, its coefficients $\mathbf{c}_j = \mathbint{0}$ for $j < i$. Hence, evaluating the Rayleigh quotient, we get,
\begin{equation}\label{eq:rayleigh-proof-eq2}
    R(\vec{\mathcal{Y}}) = \mathbf{n_x}^{-1} \circ_M \left( \sum_{j=i}^k \lambda_j \circ_M \mathbf{c}_j^2 \right) \preceq \lambda_i \circ_M \mathbf{n_x}^{-1} \circ_M \left( \sum_{j=i}^k \mathbf{c}_j^2 \right) = \lambda_i
\end{equation}
Now, using \Cref{eq:rayleigh-proof-eq1} and \Cref{eq:rayleigh-proof-eq2}, we get,
\begin{equation}
    \lambda_i(\mathcal{A}) = \max_{\mathrm{Dim}_{\mathbb{K}_M}(\mathscr{W})=i\ast \mathbint{1}} \min_{\substack{\vec{\mathcal{X}} \in \mathscr{W} \\ \mathrm{Ann}(\vec{\mathcal{X}}) = \{\mathbint{0}\}\\ \|\vec{\mathcal{X}}\|_2 = \mathbint{1}}} (\vec{\mathcal{X}}^\top \star_M \vec{\mathcal{X}})^\dagger \circ_M (\vec{\mathcal{X}}^\top \star_M \mathcal{A} \star_M \vec{\mathcal{X}})
\end{equation}
\end{proof}

\begin{lemma}[Tubal Min-Max Theorem]\label{lem:tubal-min-max}
Let $\mathcal{A} \in \mathbb{K}_M^{k \times k}$. The $i$th tubal singular value $\sigma_i(\mathcal{A})$ admits,
\begin{equation}\label{eq:primal-min-max}
    \sigma_i(\mathcal{A}) = \max_{\substack{\mathrm{Dim}_{\mathbb{K}_M}(\mathscr{W}) = i \ast \mathbint{1}}} \min_{\substack{\vec{\mathcal{X}} \in \mathscr{W} \\ \|\vec{\mathcal{X}}\|_2 = \mathbint{1} \\ \mathrm{Ann}(\vec{\mathcal{X}}) = \{\mathbint{0}\}}} \|\mathcal{A} \star_M \vec{\mathcal{X}}\|_2,
\end{equation}
and similarly,
\begin{equation}\label{eq:dual-min-max}
    \sigma_i(\mathcal{A}) = \min_{\substack{\mathrm{Dim}_{\mathbb{K}_M}(\mathscr{Y}) = (k - i + 1) \ast \mathbint{1}}} \max_{\substack{\vec{\mathcal{X}} \in \mathscr{Y} \\ \|\vec{\mathcal{X}}\|_2 = \mathbint{1} \\ \mathrm{Ann}(\vec{\mathcal{X}}) = \{\mathbint{0}\}}} \|\mathcal{A} \star_M \vec{\mathcal{X}}\|_2,
\end{equation}
where $\mathscr{W}$ and $\mathscr{Y}$ are finitely generated projective submodules of $\mathbb{K}_M^k$.
\end{lemma}
\begin{proof}
Constructing the self-adjoint positive semidefinite tensor $\mathcal{H} = \mathcal{A}^\top \star_M \mathcal{A} \in \mathbb{K}_M^{k \times k}$. By \Cref{thrm:tensor-extended-spectral-thrm}, $\mathcal{H}$ possesses eigenvalues $\lambda_i(\mathcal{H}) = \sigma_i^2(\mathcal{A})$ and an orthonormal basis of regular tubal eigenvectors $\{\vec{\mathcal{V}}_1, \dots, \vec{\mathcal{V}}_k\}$. For any regular unit vector $\vec{\mathcal{X}} \in \mathbb{K}_M^k$ satisfying $\|\vec{\mathcal{X}}\|_2 = \mathbint{1}$ and $\mathrm{Ann}(\vec{\mathcal{X}}) = \{\mathbint{0}\}$, the tubal Rayleigh quotient evaluates to $R_{\mathcal{H}}(\vec{\mathcal{X}}) = \vec{\mathcal{X}}^\top \star_M \mathcal{H} \star_M \vec{\mathcal{X}} = (\mathcal{A} \star_M \vec{\mathcal{X}})^\top \star_M (\mathcal{A} \star_M \vec{\mathcal{X}}) = \|\mathcal{A} \star_M \vec{\mathcal{X}}\|_2^2$. Applying \Cref{thrm:Tubal-Rayleigh-Quotient} to $\mathcal{H}$, we get,
\begin{equation}
    \lambda_i(\mathcal{H}) = \max_{\mathrm{Dim}_{\mathbb{K}_M}(\mathscr{W}) = i \ast \mathbint{1}} \min_{\substack{\vec{\mathcal{X}} \in \mathscr{W}\\ \|\vec{\mathcal{X}}\|_2 = \mathbint{1} \\ \mathrm{Ann}(\vec{\mathcal{X}}) = \{\mathbint{0}\}}} R_{\mathcal{H}}(\vec{\mathcal{X}}).
\end{equation}
As $\mathbb{K}_M$ is Von Neumann regular and $\sigma_i(\mathcal{A}) \succeq \mathbint{0}$, taking the unique positive tubal square root across both sides yields \Cref{eq:primal-min-max}.

Now, define the submodule $\mathscr{W}_i^\perp = \mathrm{span}_{\mathbb{K}_M}\{\vec{\mathcal{V}}_i, \dots, \vec{\mathcal{V}}_k\}$ with $\mathrm{Dim}_{\mathbb{K}_M}(\mathscr{W}_i^\perp) = (k - i + 1) \ast \mathbint{1}$. Any regular unit vector $\vec{\mathcal{Y}} \in \mathscr{W}_i^\perp$ satisfies $\vec{\mathcal{Y}} = \sum_{j=i}^k \mathbf{c}_j \circ_M \vec{\mathcal{V}}_j$ with $\sum_{j=i}^k \mathbf{c}_j^2 = \mathbint{1}$. Since $\sigma_i \succeq \sigma_j$ for all $j \ge i$, we get,
\begin{equation}
    \|\mathcal{A} \star_M \vec{\mathcal{Y}}\|_2^2 = \sum_{j=i}^k \sigma_j^2 \circ_M \mathbf{c}_j^2 \preceq \sigma_i^2 \circ_M \sum_{j=i}^k \mathbf{c}_j^2 = \sigma_i^2,
\end{equation}
implying $\max_{\vec{\mathcal{Y}} \in \mathscr{W}_i^\perp} \|\mathcal{A} \star_M \vec{\mathcal{Y}}\|_2 = \sigma_i(\mathcal{A})$. Now let $\mathscr{Y} \subseteq \mathbb{K}_M^k$ be an arbitrary projective submodule of dimension $(k - i + 1) \ast \mathbint{1}$, and define $\mathscr{U}_i = \mathrm{span}_{\mathbb{K}_M}\{\vec{\mathcal{V}}_1, \dots, \vec{\mathcal{V}}_i\}$. Now, using \Cref{thrm:general-rank-nullity-thrm}, we get $\mathrm{Dim}_{\mathbb{K}_M}(\mathscr{Y} \cap \mathscr{U}_i) \succeq \mathrm{Dim}_{\mathbb{K}_M}(\mathscr{Y}) + \mathrm{Dim}_{\mathbb{K}_M}(\mathscr{U}_i) - k \ast \mathbint{1} = \mathbint{1}$. That is, there exists at least one regular unit vector $\vec{\mathcal{Z}} \in \mathscr{Y} \cap \mathscr{U}_i$. Since $\vec{\mathcal{Z}} \in \mathscr{U}_i$, $\|\mathcal{A} \star_M \vec{\mathcal{Z}}\|_2 \succeq \sigma_i(\mathcal{A})$. Because $\vec{\mathcal{Z}} \in \mathscr{Y}$, $\max_{\vec{\mathcal{X}} \in \mathscr{Y}} \|\mathcal{A} \star_M \vec{\mathcal{X}}\|_2 \succeq \|\mathcal{A} \star_M \vec{\mathcal{Z}}\|_2 \succeq \sigma_i(\mathcal{A})$. Taking the minimum over all such submodules $\mathscr{Y}$ establishes \Cref{eq:dual-min-max}.
\end{proof}

\begin{theorem}[Tubal Weyl's Inequality]\label{thrm:tubal-weyl}
For any tensors $\mathcal{A}, \mathcal{B} \in \mathbb{K}_M^{k \times k}$ and any index $i \in \{1, \dots, k\}$, we have,
\begin{displaymath}
    |\sigma_i(\mathcal{A} + \mathcal{B}) - \sigma_i(\mathcal{A})| \preceq \|\mathcal{B}\|_2.
\end{displaymath}
\end{theorem}
\begin{proof}
Let $\mathscr{W}_i = \mathrm{span}_{\mathbb{K}_M}\{\vec{\mathcal{V}}_1, \dots, \vec{\mathcal{V}}_i\}$ denote an $i$-dimensional submodule of right singular vectors of $\mathcal{A}$. For any projective submodule $\mathscr{Y}$ of dimension $(k - i + 1) \ast \mathbint{1}$, \Cref{lem:tubal-min-max} ensures that $\mathrm{Dim}_{\mathbb{K}_M}(\mathscr{W}_i \cap \mathscr{Y}) \succeq \mathbint{1}$. Let $\vec{\mathcal{X}}_0 \in \mathscr{W}_i \cap \mathscr{Y}$ be a regular unit tubal vector. Applying the tubal triangle inequality, we get,
\begin{displaymath}
    \|(\mathcal{A} + \mathcal{B}) \star_M \vec{\mathcal{X}}_0\|_2 \preceq \|\mathcal{A} \star_M \vec{\mathcal{X}}_0\|_2 + \|\mathcal{B} \star_M \vec{\mathcal{X}}_0\|_2 \preceq \sigma_i(\mathcal{A}) + \|\mathcal{B}\|_2.
\end{displaymath}
Taking the maximum over $\vec{\mathcal{X}} \in \mathscr{Y}$ followed by the minimum over all projective submodules $\mathscr{Y}$ of dimension $(k - i + 1) \ast \mathbint{1}$, \Cref{lem:tubal-min-max} yields,
\begin{equation}\label{eq:weyl1}
    \sigma_i(\mathcal{A} + \mathcal{B}) \preceq \sigma_i(\mathcal{A}) + \|\mathcal{B}\|_2.
\end{equation}
Replacing $\mathcal{A} + \mathcal{B}$ by $\mathcal{A}$ and $\mathcal{A}$ by $(\mathcal{A} + \mathcal{B}) - \mathcal{B}$ gives,
\begin{equation}\label{eq:weyl2}
    \sigma_i(\mathcal{A}) \preceq \sigma_i(\mathcal{A} + \mathcal{B}) + \|-\mathcal{B}\|_2 = \sigma_i(\mathcal{A} + \mathcal{B}) + \|\mathcal{B}\|_2.
\end{equation}
Combining \Cref{eq:weyl1} and \Cref{eq:weyl2}, we get, $-\|\mathcal{B}\|_2 \preceq \sigma_i(\mathcal{A} + \mathcal{B}) - \sigma_i(\mathcal{A}) \preceq \|\mathcal{B}\|_2$, which implies $|\sigma_i(\mathcal{A} + \mathcal{B}) - \sigma_i(\mathcal{A})| \preceq \|\mathcal{B}\|_2$.
\end{proof}

\begin{definition}[Tubal Inner Product]
    Let $\mathcal{A}, \mathcal{B} \in \mathbb{K}_M^{m \times n}$. The tubal inner product, denoted $\langle \mathcal{A}, \mathcal{B} \rangle_{\mathbb{K}_M} \in \mathbb{K}_M$, is defined as,
    \begin{displaymath}
        \langle \mathcal{A}, \mathcal{B} \rangle_{\mathbb{K}_M} = \texttt{t-Tr}(\mathcal{A}^\top \star_M \mathcal{B}) = \sum_{i=1}^n (\mathcal{A}^\top \star_M \mathcal{B})^{[i,i]}.
    \end{displaymath}
\end{definition}

\begin{definition}[Tubal Norms]
    For a tensor $\mathcal{A} \in \mathbb{K}_M^{m \times n}$, the tubal norms are tubal scalars in $\mathbb{K}_M$ satisfying positive semi-definiteness ($\|\mathcal{A}\| \succeq \mathbint{0}$),
    \begin{enumerate}
        \item Tubal Frobenius Norm: $\|\mathcal{A}\|_F = \sqrt{\langle \mathcal{A}, \mathcal{A} \rangle_{\mathbb{K}_M}}$, where the square root is the unique positive tubal square root.
        \item Tubal Spectral Norm: $\|\mathcal{A}\|_2 = \sigma_1(\mathcal{A})$, where $\sigma_1$ is the largest principal tubal singular value defined in \Cref{thrm:tensor-extended-svd}.
    \end{enumerate}
\end{definition}

\begin{remark}[Real-Valued Global Norms]
    For practical uses, we map the tubal norms to standard real-valued norms using the supremum over the transformed ideals via the linear mapping $\Phi: \mathbb{K}_M \to \mathbb{R}^N$,
    \begin{displaymath}
        \|\mathcal{A}\|_{\mathbb{R}, 2} = \max(\Phi(\|\mathcal{A}\|_2)), \quad \text{and} \quad \|\mathcal{A}\|_{\mathbb{R}, F} = \sqrt{\Phi(\|\mathcal{A}\|_F)_j^2}.
    \end{displaymath}
\end{remark}

These definitions and inequalities are required for further perturbation analysis of the tubal tensor algebra, including perturbation bounds for the Moore-Penrose inverse computations, done in \Cref{sec:applications}.

\subsection{Gershgorin Hypersphere Theorem}
One of the important results in matrix theory with field elements is the Gershgorin theorem~\cite{gershgorin}, which bounds the eigenvalues and hence the energy of the system. However, there is no such result present in tubal tensor algebra. Hence, we propose the Gershgorin hypersphere theorem for tubal tensor algebra for bounding the tubal eigenvalues.

\begin{definition}[Strict tubal diagonal dominance]
A tensor $\mathcal{A} \in \mathbb{K}_M^{k \times k}$ is strictly tubally diagonally dominant if, for every row $i$,
\begin{displaymath}
    |\mathcal{A}^{[i,i]}| \succ \sum_{j \neq i} |\mathcal{A}^{[i,j]}|,
\end{displaymath}
where $|\cdot|$ is the tubal absolute value mapping satisfying positive definiteness, multiplicativity, and the triangle inequality.
\end{definition}

\begin{theorem}[Gershgorin hypersphere theorem]\label{thrm:Gershgorin-Hypersphere-thrm}
Every valid tubal eigenvalue $\lambda \in \mathbb{K}_M$ of the tensor $\mathcal{A} \in \mathbb{K}_M^{k \times k}$ resides within a composite Gershgorin hypersphere. Specifically, there exists a complete set of mutually orthogonal idempotents $\{\mathbf{u}_1, \dots, \mathbf{u}_k\} \subset \mathbb{K}_M$ such that,
\begin{displaymath}
    \sum_{i=1}^k \mathbf{u}_i \circ_M |\lambda - \mathcal{A}^{[i,i]}| \preceq \sum_{i=1}^k \mathbf{u}_i \circ_M \mathbf{r}_i,
\end{displaymath}
where $\mathbf{r}_i = \sum_{j \neq i} |\mathcal{A}^{[i,j]}|$ is the tubal radius of the $i$th row.
Consequently, every tubal eigenvalue $\lambda$ satisfies the global absolute bound, $|\lambda| \preceq \mathbf{r}_{max}$, where $\mathbf{r}_{max} = \sum_{i=1}^k \mathbf{v}_i \circ_M (\mathbf{r}_i + |\mathcal{A}^{[i,i]}|)$, constructed using a set of complete mutually orthogonal idempotents $\{\mathbf{v}_1, \dots, \mathbf{v}_k\} \subset \mathbb{K}_M$ such that $\mathbf{v}_m\circ_M(\mathbf{r}_m + |\mathcal{A}^{[m,m]}|) \succeq \mathbf{v}_m\circ_M(\mathbf{r}_j + |\mathcal{A}^{[j,j]}|)$ for all $j = 1,2,\dots, k$.
\end{theorem}
\begin{proof}
Let $\lambda \in \mathbb{K}_M$ be a valid tubal eigenvalue of $\mathcal{A}$, meaning there exists a regular tubal eigenvector $\vec{\mathcal{X}} \in \mathbb{K}_M^{k \times 1}$, such that $\sum_{j=1}^k \mathcal{A}^{[i,j]} \circ_M \vec{\mathcal{X}}^{[j]} = \lambda \circ_M \vec{\mathcal{X}}^{[i]}$, for all $i \in \{1, \dots, k\}$. Now, rearranging the terms, taking the tubal absolute value, and applying the tubal triangle inequality gives
\begin{equation}\label{eq:gershgorin-proof-1}
    |\lambda - \mathcal{A}^{[i,i]}| \circ_M |\vec{\mathcal{X}}^{[i]}| \preceq \sum_{j \neq i} |\mathcal{A}^{[i,j]}| \circ_M |\vec{\mathcal{X}}^{[j]}|.
\end{equation}
Now, define the local infinity norm $\mathbf{n}_\infty = \max_j |\vec{\mathcal{X}}^{[j]}| \in \mathbb{K}_M$. $\mathbf{n}_\infty$ may or may not occur at a $j$, instead, there exists a set of mutually orthogonal idempotents $\{\mathbf{u}_1, \dots, \mathbf{u}_k\}$ partitioning the ring ($\sum_{i=1}^k \mathbf{u}_i = \mathbint{1}$) such that the maximum is realized locally on each ideal, that is, $\mathbf{u}_i \circ_M \mathbf{n}_\infty = \mathbf{u}_i \circ_M |\vec{\mathcal{X}}^{[i]}|$. As $\vec{\mathcal{X}}$ is a regular tubal vector, $\mathrm{Ann}(\mathbf{n}_\infty) = \{\mathbint{0}\}$, that is, $\mathbf{n}_\infty$ is strictly invertible. Multiplying \Cref{eq:gershgorin-proof-1} by the idempotent $\mathbf{u}_i$, we have, $\mathbf{u}_i \circ_M |\lambda - \mathcal{A}^{[i,i]}| \circ_M \mathbf{n}_\infty \preceq \mathbf{u}_i \circ_M \sum_{j \neq i} |\mathcal{A}^{[i,j]}| \circ_M \mathbf{n}_\infty$. Substituting the $\mathbf{r}_i$ and multiplying both sides by $\mathbf{n}_\infty^{-1}$, we get,
\begin{equation}\label{eq:gershgorin-proof-2}
    \mathbf{u}_i \circ_M |\lambda - \mathcal{A}^{[i,i]}| \preceq \mathbf{u}_i \circ_M \mathbf{r}_i.
\end{equation}
Summing \Cref{eq:gershgorin-proof-2} over all indices $i$ from $1$ to $k$ proves the composite hypersphere bound,
\begin{displaymath}
    \sum_{i=1}^k \mathbf{u}_i \circ_M |\lambda - \mathcal{A}^{[i,i]}| \preceq \sum_{i=1}^k \mathbf{u}_i \circ_M \mathbf{r}_i.
\end{displaymath}
Now, by adding and subtracting $\mathcal{A}^{[i,i]}$ to $|\lambda|$ and applying the triangle inequality, we have $|\lambda| \preceq |\lambda - \mathcal{A}^{[i,i]}| + |\mathcal{A}^{[i,i]}|$. Multiplying by $\mathbf{u}_i$ and substituting \Cref{eq:gershgorin-proof-2}, we get, $\mathbf{u}_i \circ_M |\lambda| \preceq \mathbf{u}_i \circ_M (\mathbf{r}_i + |\mathcal{A}^{[i,i]}|)$. By the definition, $\mathbf{r}_i + |\mathcal{A}^{[i,i]}| \preceq \mathbf{r}_{max}$, holds for all $i$ across all ideals. Hence, we have, $\mathbf{u}_i \circ_M |\lambda| \preceq \mathbf{u}_i \circ_M \mathbf{r}_{max}$. Now, summing across all $k$ idempotents gives $(\sum_{i=1}^k \mathbf{u}_i) \circ_M |\lambda| \preceq (\sum_{i=1}^k \mathbf{u}_i) \circ_M \mathbf{r}_{max}$, which simplifies to $|\lambda| \preceq \mathbf{r}_{max}$.
\end{proof}

\begin{figure}[htbp]
  \subfloat[Tensor $\hat{\mathcal{A}}_2\in\mathbb{R}^{3\times 3\times 2}$.]{
      \includegraphics[width=0.275\textwidth]{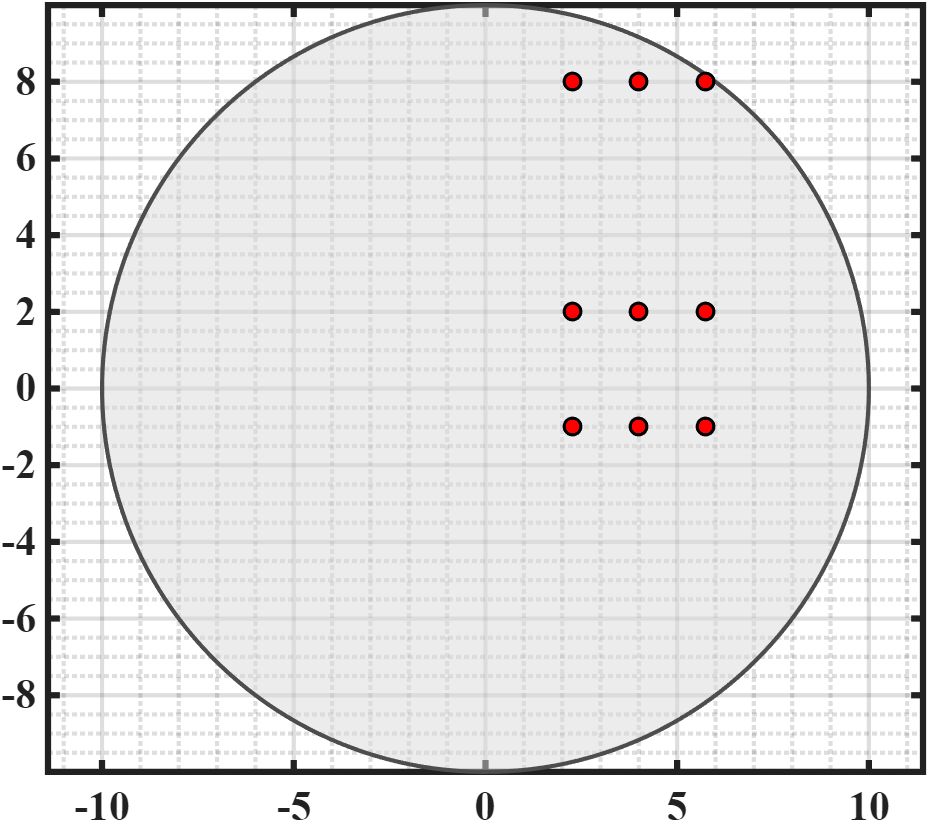}
  \label{fig:greshgorin2D}
  }
  \hfill
  \subfloat[Tensor $\hat{\mathcal{A}}_3\in\mathbb{R}^{3\times 3\times 3}$.]{
      \includegraphics[width=0.275\textwidth]{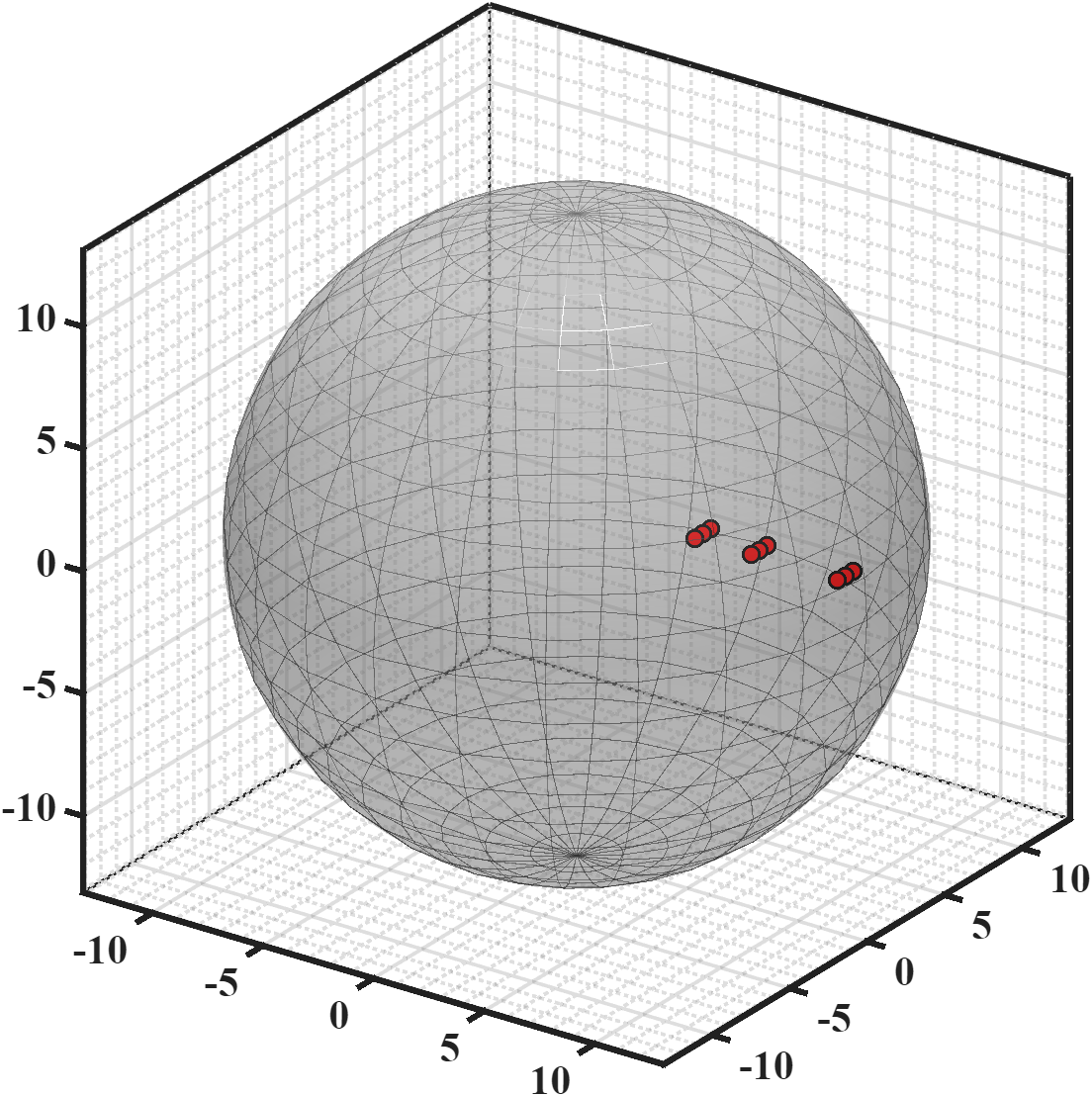}
  \label{fig:greshgorin3D}
  }
  \hfill
  \subfloat[Tensor $\hat{\mathcal{A}}_4\in\mathbb{R}^{3\times 3\times 4}$.]{
      \includegraphics[width=0.325\textwidth]{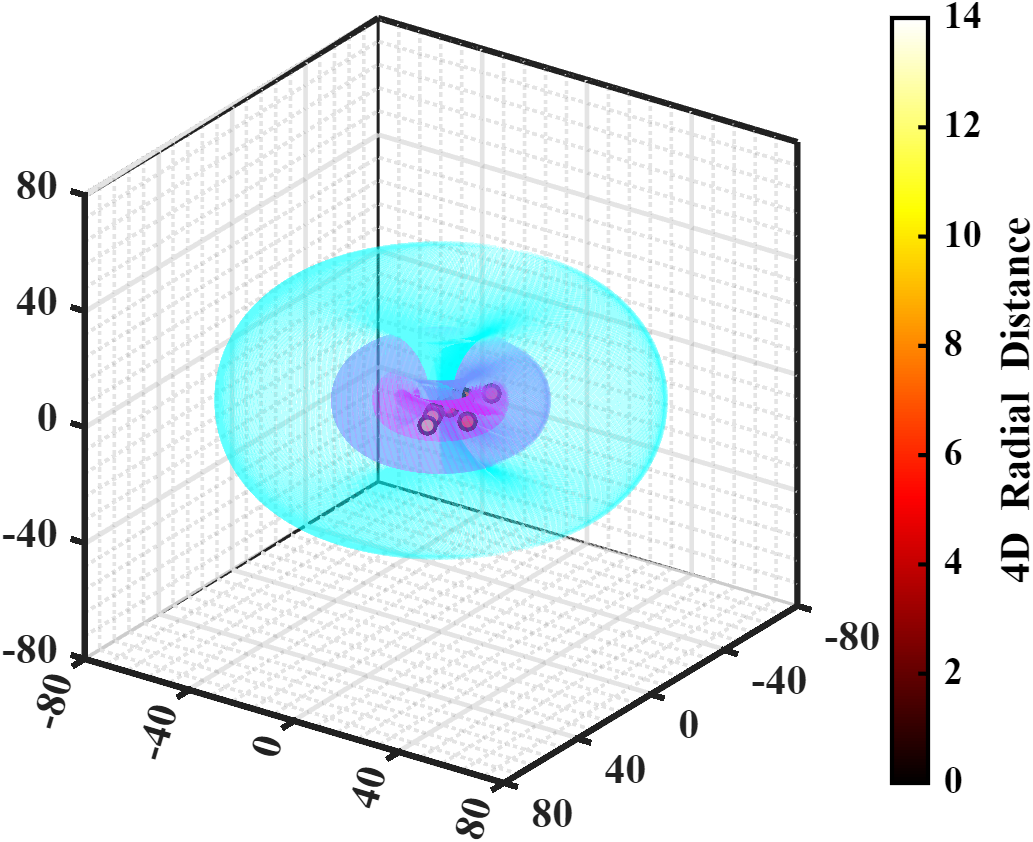}
  \label{fig:greshgorin4D}
  }
  \caption{Visual Analysis of \Cref{thrm:Gershgorin-Hypersphere-thrm} on $\Phi$-space using \Cref{example:gershgorin}.}
  \label{fig:Gershgorin-hypersphere-thrm}
\end{figure}

\begin{example}\label{example:gershgorin}
    To demonstrate \Cref{thrm:Gershgorin-Hypersphere-thrm}, we calculate the hypersphere bounds for three different tensor cases.
    \begin{enumerate}
        \item Assuming tensor $\mathcal{A}_2 \in \mathbb{R}^{3\times 3\times 2}$, that is, $N=2$ with the real transformation matrix $M = \frac{1}{\sqrt{2}}\left[\begin{smallmatrix} 1 & 1 \\ 1 & -1 \end{smallmatrix}\right]$,
        \begin{displaymath}
            \mathcal{A}_2[:,:,1] = \frac{1}{\sqrt{2}} \begin{bmatrix} 7 & 1 & 0 \\ 1 & 12 & 1 \\ 0 & 1 & 2 \end{bmatrix}, \quad 
            \mathcal{A}_2[:,:,2] = \frac{1}{\sqrt{2}} \begin{bmatrix} 3 & 1 & 0 \\ 1 & -4 & 1 \\ 0 & 1 & 4 \end{bmatrix}.
        \end{displaymath}
        Applying the transformation $\hat{\mathcal{A}} = \mathfrak{M}(\mathcal{A}_2)$,
        \begin{displaymath}
            \hat{\mathcal{A}}_2[:,:,1] = \begin{bmatrix} 5 & 1 & 0 \\ 1 & 4 & 1 \\ 0 & 1 & 3 \end{bmatrix}, \quad 
            \hat{\mathcal{A}}_2[:,:,2] = \begin{bmatrix} 2 & 0 & 0 \\ 0 & 8 & 0 \\ 0 & 0 & -1 \end{bmatrix}.
        \end{displaymath}
        In the transform domain, we have $|\hat{\mathcal{A}}_2^{[1,1]}| + \hat{\mathbf{r}}_1 = [6, 2]^\top$, $|\hat{\mathcal{A}}_2^{[2,2]}| + \hat{\mathbf{r}}_2 = [6, 8]^\top$ and $|\hat{\mathcal{A}}_2^{[3,3]}| + \hat{\mathbf{r}}_3 = [4, 1]^\top$. Now, computing $\hat{\mathbf{r}}_{max}$, across the natural fundamental independent ideals, $\hat{\mathbf{r}}_{max} = [\max(6, 6, 4), \max(2, 8, 1)]^\top = [6, 8]^\top$. Hence, $\mathbf{r}_{max} = M [6, 8]^\top$, which effectively encloses the tubal eigenvalues. Now, here we have $3^2 = 9$ possible eigenvalues. However, all of them follow, $|\lambda|\preceq \mathbf{r}_{max}$. The eigenvalues evaluated in the transform domain are $\hat{\Lambda}_1 = \{4, 4+\sqrt{3}, 4-\sqrt{3}\}$ and $\hat{\Lambda}_2 = \{2, 8, -1\}$. And the $\lambda$ is formed as, $\lambda = M[\hat{\Lambda}_1, \hat{\Lambda}_2]^T$. The same is represented in \Cref{fig:greshgorin2D}.

        \item Assuming tensor $\mathcal{A}_3 \in \mathbb{R}^{3\times 3\times 3}$, that is, $N=3$ with the real transformation matrix $M = I_3$,
        \begin{displaymath}
            \hat{\mathcal{A}}_3[:,:,1] = \begin{bmatrix} 10 & 2 & 1 \\ 1 & 8 & 1 \\ 1 & 1 & 6 \end{bmatrix}, \quad 
            \hat{\mathcal{A}}_3[:,:,2] = \begin{bmatrix} 1 & 0.5 & 0 \\ 0.5 & 1 & 0 \\ 0 & 0 & 1 \end{bmatrix}, \quad 
            \hat{\mathcal{A}}_3[:,:,3] = I_3.
        \end{displaymath}
        In the transform domain, we have $|\hat{\mathcal{A}}_3^{[1,1]}| + \hat{\mathbf{r}}_1 = [13, 1.5, 1]^\top$, $|\hat{\mathcal{A}}_3^{[2,2]}| + \hat{\mathbf{r}}_2 = [10, 1.5, 1]^\top$ and $|\hat{\mathcal{A}}_3^{[3,3]}| + \hat{\mathbf{r}}_3 = [8, 1, 1]^\top$.
        Now, computing $\hat{\mathbf{r}}_{max}$, across the natural fundamental independent ideals, $\hat{\mathbf{r}}_{max} = [13,1.5,1]^\top$. Hence, $\mathbf{r}_{max} = [13, 1.5, 1]^\top$, which effectively encloses the tubal eigenvalues. Now, here we have $3^3 = 27$ possible eigenvalues. However, all of them follow, $|\lambda|\preceq \mathbf{r}_{max}$. The eigenvalues evaluated in the transform domain are $\hat{\Lambda}_1 = \{11.454, 7.551, 4.995\}$, $\hat{\Lambda}_2 = \{1.5, 1, 0.5\}$, and $\hat{\Lambda}_3 = \{1, 1, 1\}$. And the $\lambda$ is formed as, $\lambda = [\hat{\Lambda}_1, \hat{\Lambda}_2, \hat{\Lambda}_3]^T$. The same is represented in \Cref{fig:greshgorin3D}.

        \item Assuming tensor $\mathcal{A}_4 \in \mathbb{R}^{3\times 3\times 4}$, that is, $N=4$ with the real transformation matrix $M$ same as \Cref{example:determinant},
        \begin{displaymath}
            \mathcal{A}_4[:,:,1] = \begin{bmatrix} 6 & 1.5 & 0 \\ 0.5 & 0 & 0 \\ 0 & 1 & 5.5 \end{bmatrix}, ~ 
            \mathcal{A}_4[:,:,2] = \begin{bmatrix} 4 & 0.5 & 0 \\ 0.5 & -5 & 0 \\ 0 & 0 & -3.5 \end{bmatrix}, ~
        \end{displaymath}
        \begin{displaymath}
            \mathcal{A}_4[:,:,3] = \begin{bmatrix} 8 & 0.5 & 0 \\ 0.5 & -1 & 0 \\ 0 & 0 & 0.5 \end{bmatrix}, ~ 
            \mathcal{A}_4[:,:,4] = \begin{bmatrix} 2 & 1.5 & 0 \\ 0.5 & -4 & 0 \\ 0 & 1 & 1.5 \end{bmatrix}.
        \end{displaymath}
        Applying the transformation $\hat{\mathcal{A}} = \mathfrak{M}(\mathcal{A}_4)$,
        \begin{displaymath}
            \hat{\mathcal{A}}_4[:,:,1] = \begin{bmatrix} 10 & 2 & 0 \\ 1 & -5 & 0 \\ 0 & 1 & 2 \end{bmatrix},~ 
            \hat{\mathcal{A}}_4[:,:,2] = \begin{bmatrix} 4 & 0 & 0 \\ 0 & 4 & 0 \\ 0 & 0 & 4 \end{bmatrix},~
        \end{displaymath}
        \begin{displaymath}
            \hat{\mathcal{A}}_4[:,:,3] = \begin{bmatrix} 0 & 0 & 0 \\ 0 & 0 & 0 \\ 0 & 0 & 0 \end{bmatrix},~ 
            \hat{\mathcal{A}}_4[:,:,4] = \begin{bmatrix} -2 & 1 & 0 \\ 0 & 1 & 0 \\ 0 & 1 & 5 \end{bmatrix}.
        \end{displaymath}
        In the transform domain, we have $|\hat{\mathcal{A}}_4^{[1,1]}| + \hat{\mathbf{r}}_1 = [12, 4, 0, 3]^\top$, $|\hat{\mathcal{A}}_4^{[2,2]}| + \hat{\mathbf{r}}_2 = [6, 4, 0, 1]^\top$ and $|\hat{\mathcal{A}}_4^{[3,3]}| + \hat{\mathbf{r}}_3 = [3, 4, 0, 6]^\top$.
        Now, computing $\hat{\mathbf{r}}_{max}$, across the natural fundamental independent ideals, $\hat{\mathbf{r}}_{max} = [12, 4, 0, 6]^\top$. Hence, $\mathbf{r}_{max} = M[12, 4, 0, 6]^\top$, which effectively encloses the tubal eigenvalues. Now, here we have $3^4 = 81$ possible eigenvalues. However, all of them follow, $|\lambda|\preceq \mathbf{r}_{max}$. The eigenvalues evaluated in the transform domain are $\hat{\Lambda}_1 = \{10.13, -5.13, 2\}$, $\hat{\Lambda}_2 = \{4, 4, 4\}$, $\hat{\Lambda}_3 = \{0, 0, 0\}$, and $\hat{\Lambda}_4 = \{-2, 1, 5\}$. And the $\lambda$ is formed as, $\lambda = [\hat{\Lambda}_1, \hat{\Lambda}_2, \hat{\Lambda}_3, \hat{\Lambda}_4]^T$. The same is represented in \Cref{fig:greshgorin4D}.
    \end{enumerate}
\end{example}

\section{Numerical Illustrations in Tensor Computations}\label{sec:applications}
To demonstrate the empirical utility of the proposed tubal algebra, we apply our theoretical results to two fundamental computational scenarios. To ensure statistical robustness, all numerical experiments are averaged over $10$ seeds from $42$ to $51$. We utilize the discrete cosine transform (DCT) matrix for all the experiments, which is generally known as the $c$-product~\cite{kernfeld2015tensor}, a variant of the $M$-product. Experiments are performed using an Intel Core $i9$-$12900$K CPU, $32$GB RAM, and $16$ cores with an NVIDIA GeForce GTX $1660$ SUPER GPU with $6$GB of GDDR$6$ memory. The code is written and executed in Python $3.11.0$. The plots are generated using MATLAB R2025a. 

\subsection{Numerical Illustrations in Tensor Optimization}
Gradient descent is one of the classical and most important algorithms in computational mathematics. However, even with so many advancements, choosing an optimal learning rate $\eta$ is an open area in tensor computations due to the existence of tubal tensor algebra. Considering the unconstrained tensor least-squares optimization problem as follows,
\begin{equation}\label{eq:gradient_descent}
    \min_{\mathcal{X}} f(\mathcal{X}) = \frac{1}{2} \|\mathcal{A} \star_M \mathcal{X} - \mathcal{B}\|_F^2,
\end{equation}
where $\mathcal{A} \in \mathbb{K}_M^{m \times k}$, $\mathcal{X} \in \mathbb{K}_M^{k \times \ell}$, and $\mathcal{B} \in \mathbb{K}_M^{m \times \ell}$. A standard tensor gradient descent update takes the form $\mathcal{X}_{t+1} = \mathcal{X}_t - \eta \ast \nabla f(\mathcal{X}_t)$, where the step size is $\eta \in \mathbb{R}$ and the tubal gradient is derived algebraically as $\nabla f(\mathcal{X}) = \mathcal{A}^\top \star_M (\mathcal{A} \star_M \mathcal{X} - \mathcal{B})$. 

The strict stability condition for gradient descent requires the step size to satisfy $\eta < 2 / \max(\Phi(\lambda_{max}(\mathcal{H})))$, where $\mathcal{H} = \mathcal{A}^\top \star_M \mathcal{A} \in \mathbb{K}_M^{k \times k}$ is the symmetric tubal Hessian, and $\Phi: \mathbb{K}_M \to \mathbb{R}^N$ is the linear isomorphism mapping defined in \Cref{def:transfomed-domain-mapping}. Computing the exact maximal eigenvalue $\lambda_{max}(\mathcal{H})$ via \Cref{thrm:tensor-extended-svd} provides the exact stability limit, but incurs a severe initialization cost of $\mathbf{O}((mk + k^2)N^2 + mk^2N + k^3 N)$. We formulate the exact step size dynamically near the absolute edge of stability as $\eta_{exact} = 1.99 / \max(\Phi(\lambda_{max}))$, alongside a theoretically divergent step $\eta_{div} = 2.05 / \max(\Phi(\lambda_{max}))$. Alternatively, applying the \Cref{thrm:Gershgorin-Hypersphere-thrm} allows us to compute a rigorously safe upper bound on $\lambda_{max}$ without entering the spectral domain. We compute the tubal radius algebraically as $\mathbf{r}_{max} = \sum_{i=1}^k \mathbf{v}_i \circ_M (\sum_{j \neq i} |\mathcal{H}^{[i,j]}| + |\mathcal{H}^{[i,i]}|)$, which scales strictly at $\mathbf{O}((mk + k^2)N^2 + mk^2N + k^2 N)$. The Gershgorin-derived step size is defined as $\eta_{g} = 1.99 / \max(\Phi(\mathbf{r}_{max}))$.

\subsubsection{Gershgorin and Tubal-SVD Computation Time}\label{subsubsec:gershgorin-vs-svd}
\begin{figure*}[htbp]
    \begin{minipage}{0.54\textwidth}
        \resizebox{\linewidth}{!}{%
        \begin{tabular}{ccc}
        \hline
        Tensor Size $k$ & Gershgorin Time (s) & SVD Time (s) \\
        \hline
        2 & $1.15 \times 10^{-4} \pm 5.56 \times 10^{-6}$ & $6.05 \times 10^{-5} \pm 5.58 \times 10^{-6}$ \\
        4 & $6.44 \times 10^{-5} \pm 3.28 \times 10^{-6}$ & $4.79 \times 10^{-5} \pm 8.80 \times 10^{-7}$ \\
        8 & $6.19 \times 10^{-5} \pm 1.86 \times 10^{-6}$ & $7.35 \times 10^{-5} \pm 5.86 \times 10^{-7}$ \\
        16 & $6.55 \times 10^{-5} \pm 8.96 \times 10^{-7}$ & $1.58 \times 10^{-4} \pm 1.98 \times 10^{-6}$ \\
        32 & $8.12 \times 10^{-5} \pm 1.10 \times 10^{-6}$ & $4.58 \times 10^{-4} \pm 1.04 \times 10^{-5}$ \\
        64 & $1.34 \times 10^{-4} \pm 2.02 \times 10^{-5}$ & $1.56 \times 10^{-3} \pm 4.00 \times 10^{-5}$ \\
        128 & $6.57 \times 10^{-4} \pm 9.24 \times 10^{-5}$ & $7.81 \times 10^{-3} \pm 2.80 \times 10^{-4}$ \\
        256 & $2.13 \times 10^{-3} \pm 8.25 \times 10^{-5}$ & $1.02 \times 10^{-1} \pm 4.72 \times 10^{-3}$ \\
        512 & $1.94 \times 10^{-2} \pm 8.19 \times 10^{-4}$ & $8.48 \times 10^{-1} \pm 2.22 \times 10^{-2}$ \\
        1024 & $9.97 \times 10^{-2} \pm 5.03 \times 10^{-3}$ & $3.75 \times 10^{0} \pm 4.38 \times 10^{-2}$ \\
        \hline
        \end{tabular}}
        \captionof{table}{Execution Time Comparison for experiment given in \Cref{subsubsec:gershgorin-vs-svd}.}
        \label{tab:time}
    \end{minipage}\hfill
    \begin{minipage}{0.39\textwidth}
        \centering
            \includegraphics[width=\linewidth]{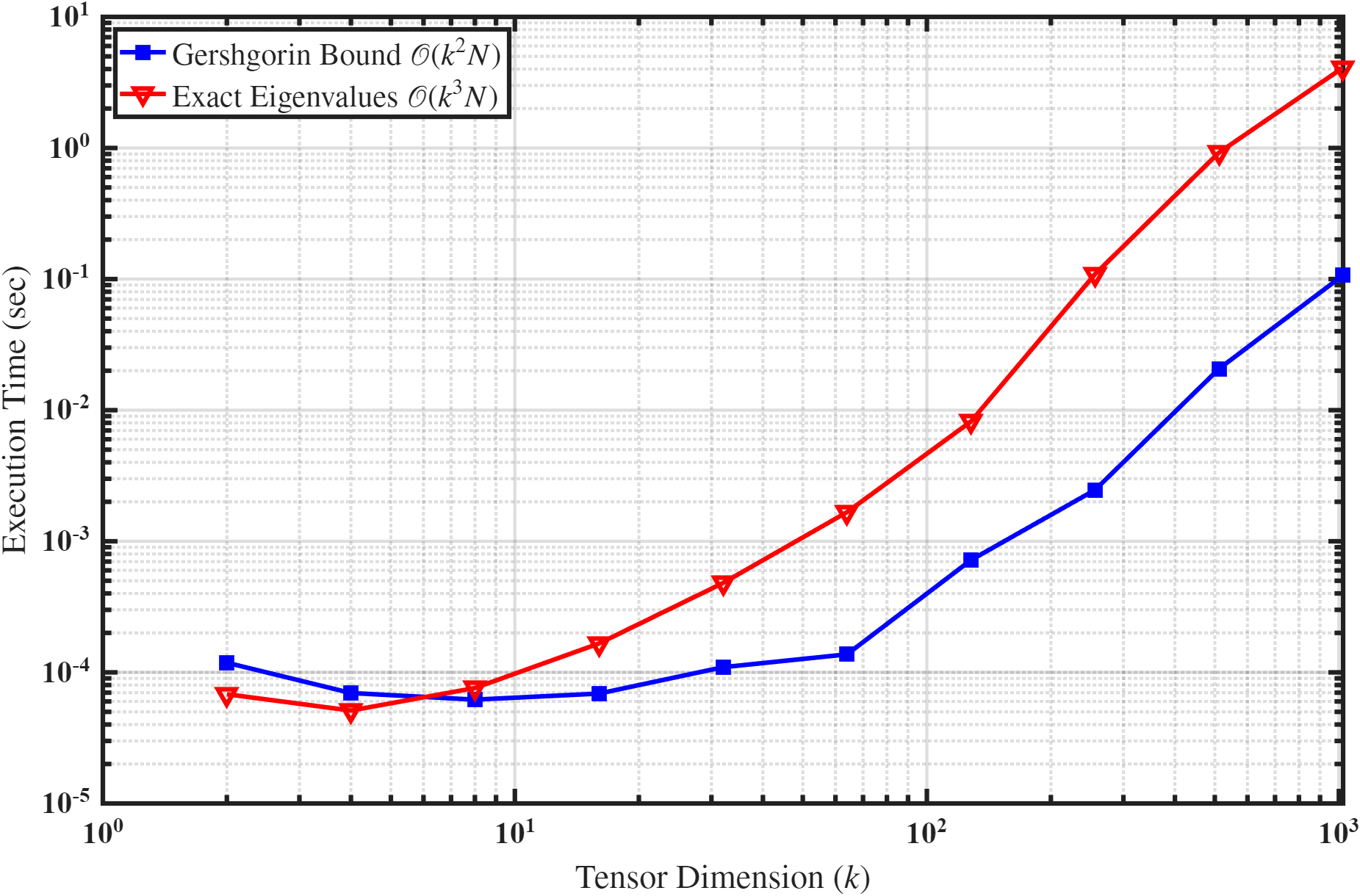}
            \captionof{figure}{Visual Analysis for experiment given in \Cref{subsubsec:gershgorin-vs-svd}.}
        \label{fig:time}
    \end{minipage}
\end{figure*}

Assuming tensors of size $\mathcal{A}\in\mathbb{R}^{k \times k \times 16}$ with various values of $k$ as $\{2^1, 2^2, \dots, 2^{10}\}$. We have compared the timings using the \texttt{numpy} library against \texttt{np.linalg.eigvalsh}. \Cref{tab:time} and \Cref{fig:time} details the execution time divergence as the spatial dimension $k$ scales, confirming that the proposed hypersphere bound, $\mathbf{O}(272k^2)$ works faster than classic eigenvalue decomposition, $\mathbf{O}(256k^2 + 16k^3)$ for large-scale tensors.

\subsubsection{Gradient Descent Algorithm and Stability}\label{subsubsec:gradient-descent}
To guarantee monotonic convergence, we require $\eta < \frac{2}{\lambda_{\max}(\mathcal{H})}$ where, $\lambda_{\max}(\mathcal{H}) = \max_{j, i} \lambda_i(\widehat{\mathcal{H}}^{(j)})$ where $\widehat{\mathcal{H}}^{(j)}$ is the $j$th slice in transformed domain. For testing, we assumed $3$ different learning rates as follows,
\begin{enumerate}
    \item $\eta_{svd} = \frac{1.99}{\lambda_{\max}(\mathcal{H})}$. The decay rate of the principal eigenvector becomes $|1 - 1.99| = 0.99$, making the steepest direction structurally the slowest to converge.
    \item $\eta_{div} = \frac{2.05}{\lambda_{\max}(\mathcal{H})}$. The spectral radius exceeds $1$. The error along the principal eigenvector amplifies by a factor of $|1 - 2.05| = 1.05$ per iteration, causing signed oscillatory divergence.
    \item $\eta_g = \frac{1.99}{\mathbf{r}_{max}}$. This intrinsic overestimation strictly bounds the step size away from the stability edge, instantaneously annihilating high-frequency errors.
\end{enumerate}

\begin{figure}[htbp]
  \subfloat[Convergence Analysis.]{
      \includegraphics[width=0.3\textwidth]{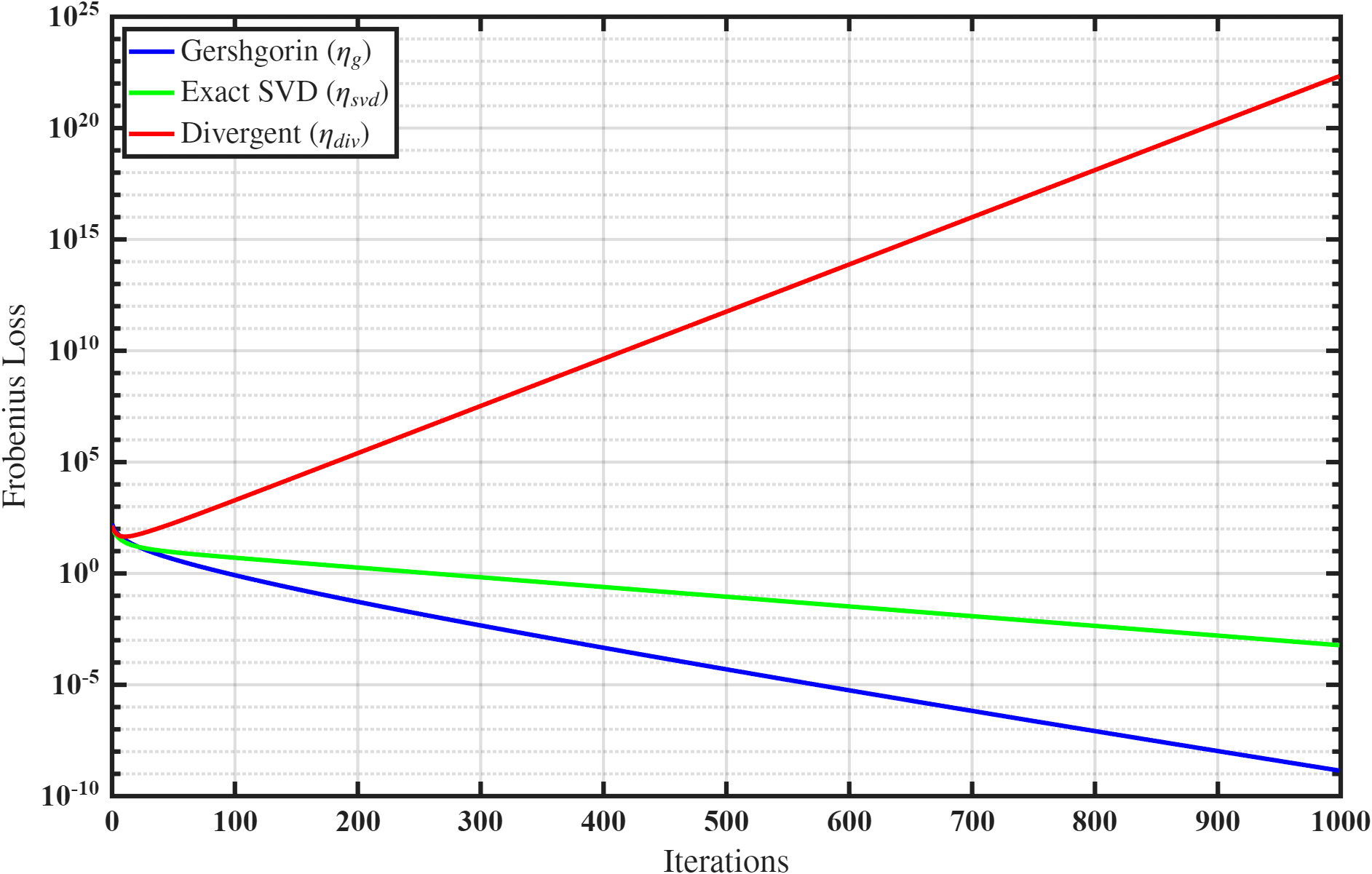}
  \label{fig:cvg}
  }
  \hfill
  \subfloat[Oscillations Analysis.]{
      \includegraphics[width=0.3\textwidth]{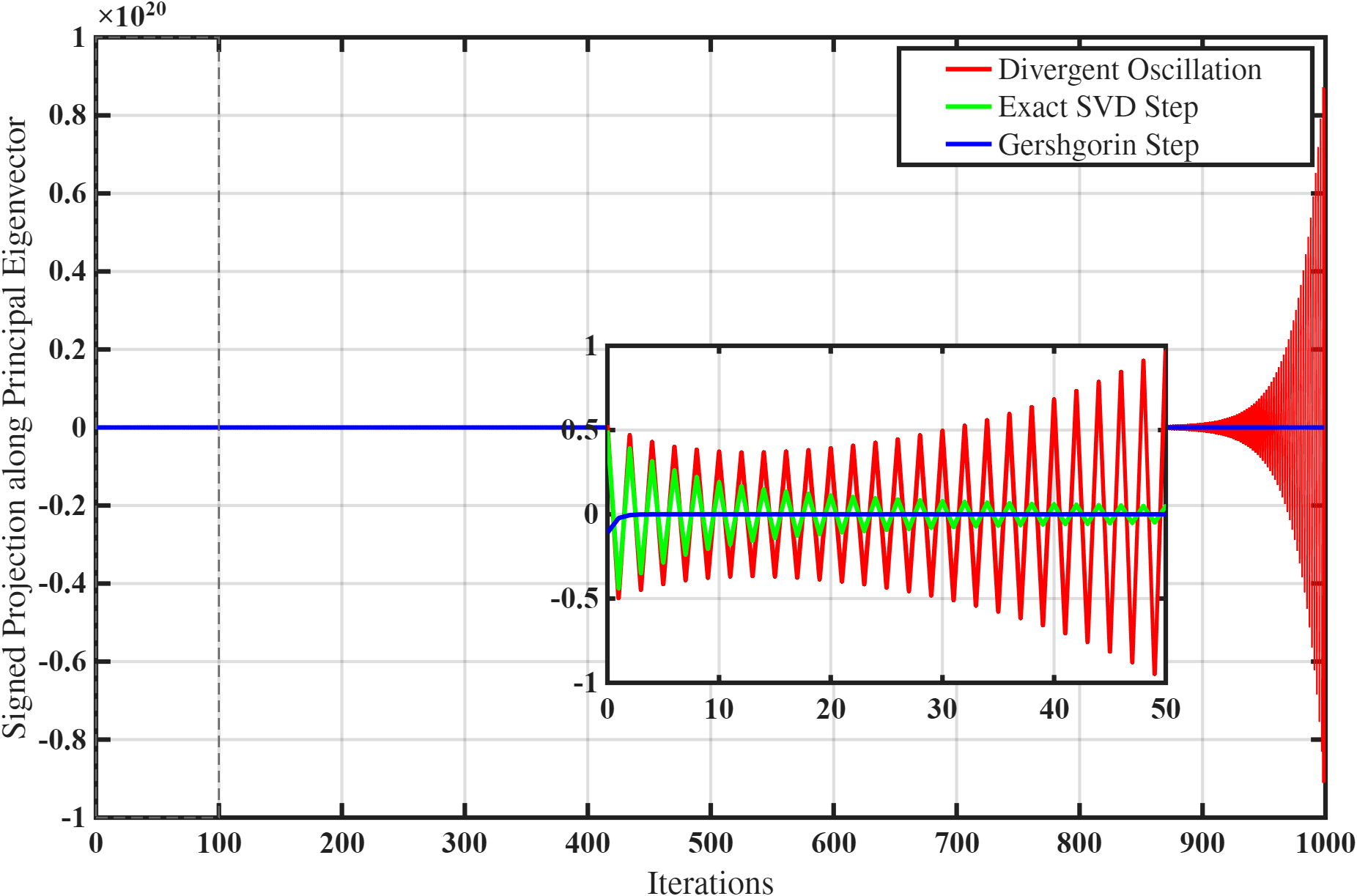}
  \label{fig:oscillate}
  }
  \hfill
  \subfloat[Rayleigh Analysis.]{
      \includegraphics[width=0.3\textwidth]{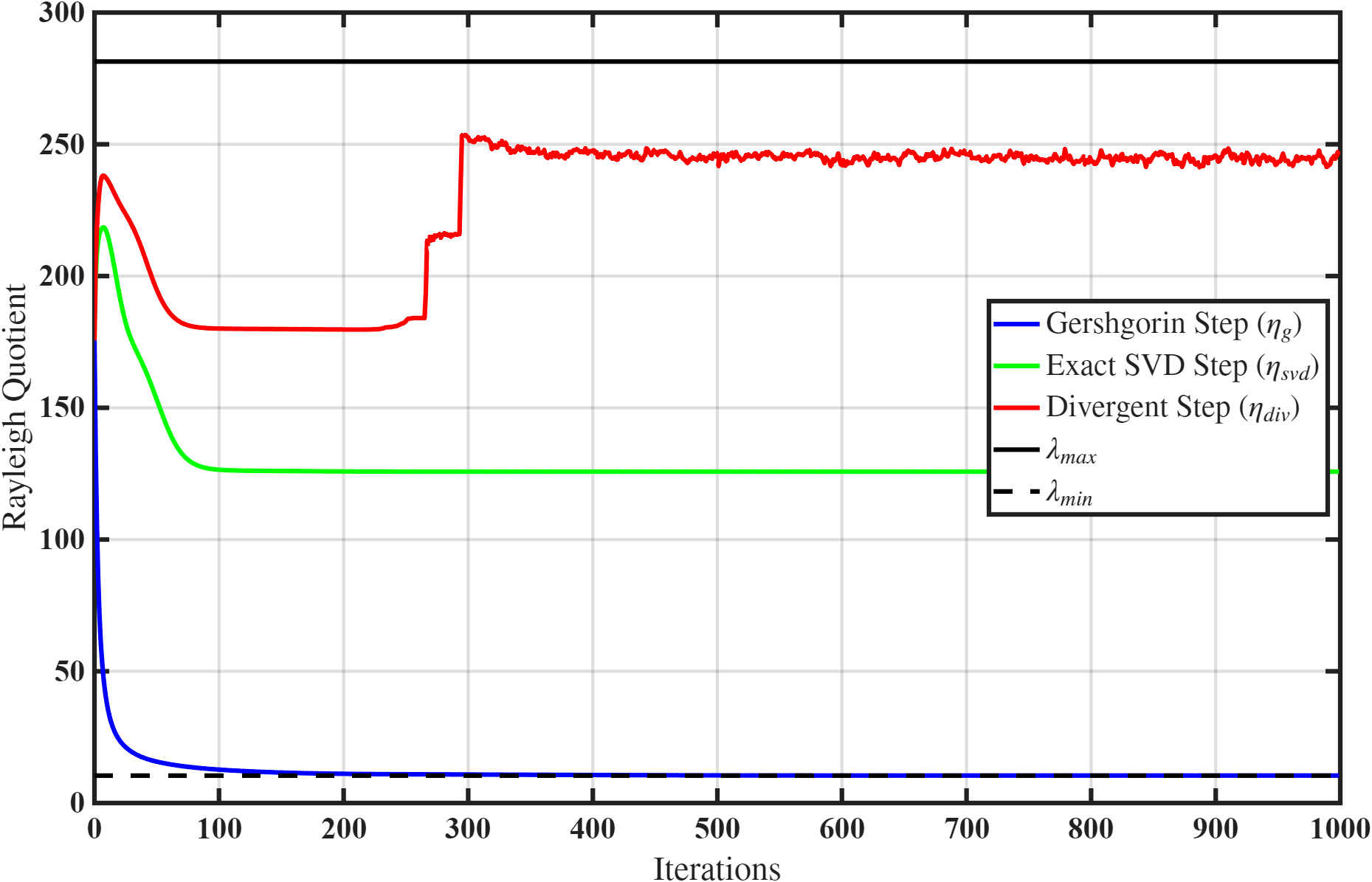}
  \label{fig:rayleigh}
  }
  \caption{Visual Analysis of the experiment given in \Cref{subsubsec:gradient-descent}.}
  \label{fig:gradient-descent}
\end{figure}

Assuming tensor $\mathcal{A}\in\mathbb{R}^{100 \times 50 \times 16}$ with target tensor $\vec{\mathcal{X}}_{\mathrm{true}}\in\mathbb{R}^{50 \times 1 \times 16}$ in \Cref{eq:gradient_descent} with maximum $1000$ iterations. For observation and analysis, we take only the $1$st slice of the tensor. To empirically observe the active subspace dominating the optimizer, we track the Rayleigh quotient of the gradient. As the optimizer solves the flat subspaces given by $\lambda_{\min}$, $R(\nabla f(\vec{\mathcal{X}})) \to \lambda_{\max}$ if the step size oscillates, but smoothly transitions to $R(\nabla f(\vec{\mathcal{X}})) \to \lambda_{\min}$ under the Gershgorin regime as shown in \Cref{fig:rayleigh}. The exact oscillatory state is tracked via the signed projection, which is given by, $\text{Proj}(t) = \mathbf{v}_{\max}^T \hat{\mathcal{E}}_t^{(0)}$, where $\mathbf{v}_{\max}$ is the principal eigenvector of $\hat{\mathcal{H}}^{(0)}$.

In \Cref{fig:oscillate}, it can be seen for learning rates $\eta_g$, $\eta_{svd}$ and $\eta_{div}$ have total oscillatory variations as $1.22 \times 10^{-1}$, $4.35 \times 10^{1}$ and $1.10 \times 10^{22}$ respectively. In \Cref{fig:cvg}, the final losses for $\eta_g$, $\eta_{svd}$ and $\eta_{div}$ are $1.40 \times 10^{-9} \pm 1.70 \times 10^{-9}$, $5.98 \times 10^{-4} \pm 3.60 \times 10^{-4}$ and $2.14 \times 10^{22} \pm 1.29 \times 10^{22}$ respectively. From \Cref{fig:cvg}, it is clear that $\eta_{div}$ strictly diverges. However, the Gershgorin step size converges faster and more smoothly than the mathematically maximal SVD step size. Also, \Cref{fig:oscillate} shows high oscillations for $\eta_{div}$ but very minimal oscillations for $\eta_g$ for $\text{Proj}(t) = \mathbf{v}_{\max}^T \hat{\mathcal{E}}_t^{(0)}$.

\subsection{Numerical Illustrations in Tensor Pseudo-Inversion}
Computing the Moore-Penrose pseudoinverse $\mathcal{A}^\dagger$ of a noise-corrupted tensor $\mathcal{A}_{noisy} = \mathcal{A}_{clean} + \Delta\mathcal{A} \in \mathbb{K}_M^{k \times k}$ is required in many application fields of mathematics. When the underlying tensor suffers from a McCoy defect ($\mathbint{r}(\mathcal{A}_{clean}) < \mathfrak{r}_t(\mathcal{A}_{clean})$), specific tubal singular values exist purely as zero-divisors on a defective ideal $\mathbf{u}_{defect}$. Under naive multidimensional unfolding, injected measurement noise ($\|\Delta\mathcal{A}\|_F = \epsilon$) causes these exact zero-divisors to be infinitely small, inducing a huge numerical explosion.

By leveraging \Cref{thrm:general-rank-nullity-thrm}, we define an exact perturbation bound by algebraically isolating the regular submodules from the defective ones. We extract the true effective singular values by formulating idempotents $\mathbf{e}_i \in \mathbb{K}_M$ that strictly map to ideals where $\Phi(\sigma_i) > \tau$, where $\tau$ is an algebraic noise threshold. The algebraic tubal pseudoinverse is rigorously constructed as follows,
\begin{equation}
    \mathcal{A}_{alg}^\dagger = \sum_{i=1}^k \left(\mathbf{e}_i \circ_M \sigma_i^\dagger\right) \star_M \vec{\mathcal{V}}_i \star_M \vec{\mathcal{U}}_i^\top.
\end{equation}

Now, we prove the importance of taking algebraic properties into consideration for the calculation of the Moore-Penrose inverse.
\begin{theorem}\label{thrm:tubal-resolvent}
Assume $\mathcal{A}, \tilde{\mathcal{A}} \in \mathbb{K}_M^{k \times k}$ satisfying $\mathscr{R}(\tilde{\mathcal{A}}) = \mathscr{R}(\mathcal{A})$ and $\mathscr{N}(\tilde{\mathcal{A}}) = \mathscr{N}(\mathcal{A})$. Then we have,
\begin{displaymath}
    \tilde{\mathcal{A}}^\dagger - \mathcal{A}^\dagger = -\tilde{\mathcal{A}}^\dagger \star_M (\tilde{\mathcal{A}} - \mathcal{A}) \star_M \mathcal{A}^\dagger.
\end{displaymath}
\end{theorem}
\begin{proof}
As $\mathscr{R}(\tilde{\mathcal{A}}) = \mathscr{R}(\mathcal{A})$ and $\mathscr{N}(\tilde{\mathcal{A}}) = \mathscr{N}(\mathcal{A})$, the orthogonal projectors onto the range and row submodules are identical, which means
\begin{displaymath}
    \tilde{\mathcal{A}} \star_M \tilde{\mathcal{A}}^\dagger = \mathcal{A} \star_M \mathcal{A}^\dagger = \mathcal{P}_{\mathscr{R}}, \quad \tilde{\mathcal{A}}^\dagger \star_M \tilde{\mathcal{A}} = \mathcal{A}^\dagger \star_M \mathcal{A} = \mathcal{I} - \mathcal{P}_{\mathscr{N}}.
\end{displaymath}
Using these and Moore-Penrose condition $\mathcal{X} \star_M \mathcal{Y} \star_M \mathcal{X} = \mathcal{X}$, we get the required results.
\end{proof}

\begin{lemma}\label{lem:tubal-inversion-bound}
Assuming $\Delta\mathcal{A} = \tilde{\mathcal{A}} - \mathcal{A}$ with $\mathscr{R}(\tilde{\mathcal{A}}) = \mathscr{R}(\mathcal{A})$, $\mathscr{N}(\tilde{\mathcal{A}}) = \mathscr{N}(\mathcal{A})$, and $\|\mathcal{A}^\dagger\|_2 \|\Delta\mathcal{A}\|_2 \prec \mathbint{1}$. Then we have,
\begin{displaymath}
    \|\tilde{\mathcal{A}}^\dagger - \mathcal{A}^\dagger\|_F \preceq \frac{\|\mathcal{A}^\dagger\|_2^2\circ_M \|\Delta\mathcal{A}\|_F}{\mathbint{1} - \|\mathcal{A}^\dagger\|_2\circ_M \|\Delta\mathcal{A}\|_2}.
\end{displaymath}
\end{lemma}
\begin{proof}
Using \Cref{thrm:tubal-resolvent}, $\tilde{\mathcal{A}}^\dagger \star_M \left( \mathcal{I} + \Delta\mathcal{A} \star_M \mathcal{A}^\dagger \right) = \mathcal{A}^\dagger$. Now, as $\|\Delta\mathcal{A} \star_M \mathcal{A}^\dagger\|_2 \preceq \|\Delta\mathcal{A}\|_2 \circ_M \|\mathcal{A}^\dagger\|_2 \prec \mathbint{1}$, the tensor $(\mathcal{I} + \Delta\mathcal{A} \star_M \mathcal{A}^\dagger)$ is invertible via the tubal Neumann series $\sum_{p=0}^\infty (-\Delta\mathcal{A} \star_M \mathcal{A}^\dagger)^p$. Hence, we can write,
\begin{displaymath}
    \tilde{\mathcal{A}}^\dagger - \mathcal{A}^\dagger = -\mathcal{A}^\dagger \star_M \Delta\mathcal{A} \star_M \mathcal{A}^\dagger \star_M \left( \mathcal{I} + \Delta\mathcal{A} \star_M \mathcal{A}^\dagger \right)^{-1}.
\end{displaymath}
Taking the Frobenius norm on both sides, we get,
\begin{displaymath}
    \|\tilde{\mathcal{A}}^\dagger - \mathcal{A}^\dagger\|_F \preceq \|\mathcal{A}^\dagger\|_2^2\circ_M \|\Delta\mathcal{A}\|_F\circ_M \left\| \sum_{p=0}^\infty (-\Delta\mathcal{A} \star_M \mathcal{A}^\dagger)^p \right\|_2 \nonumber \preceq \frac{\|\mathcal{A}^\dagger\|_2^2\circ_M \|\Delta\mathcal{A}\|_F}{\mathbint{1} - \|\mathcal{A}^\dagger\|_2\circ_M \|\Delta\mathcal{A}\|_2}.
\end{displaymath}
\end{proof}

\begin{theorem}[Perturbation Bounds]\label{thrm:perturbation-bounds}
Let $\mathcal{A} \in \mathbb{K}_M^{k \times k}$ possess a McCoy defect such that $\sigma_j = \mathbint{0}$ on ideal $\mathbf{u}_{\mathrm{defect}}$. Let $\mathcal{A}_{\mathrm{noisy}} = \mathcal{A} + \Delta\mathcal{A}$ with noise bounded by $\|\Delta\mathcal{A}\|_{\mathbb{R},F} = \epsilon$.
\begin{enumerate}
    \item Under naive unfolding, the local error for a perturbed zero-divisor $\tilde{\sigma}_j \sim \mathbf{O}(\epsilon \ast \mathbint{1})$ diverges unboundedly as noise approaches zero,
    \begin{displaymath}
        \lim_{\epsilon \to 0} \|\mathcal{A}^\dagger - \mathcal{A}_{\mathrm{naive}}^\dagger\|_F \succeq \left\| \mathbint{0} - \frac{1}{\mathbf{O}(\epsilon \ast \mathbint{1})} \right\|_F = \mathbf{O}(\epsilon^{-1} \ast \mathbint{1})\approx \mathbf{O}(\epsilon^{-1}).
    \end{displaymath}
    \item Using \Cref{thrm:general-rank-nullity-thrm}, the algebraic projector explicitly annihilates the null-space submodule $\mathscr{N}(\mathcal{A})$, perfectly nullifying the degenerate noise ($\|\mathbint{0} - \mathbint{0}\|_{\mathbb{R},F} = 0$). On the regular submodule $\mathscr{R}(\mathcal{A})$, the perturbation limits applies as,
    \begin{displaymath}
        \lim_{\epsilon \to 0} \|\mathcal{A}^\dagger - \mathcal{A}_{\mathrm{alg}}^\dagger\|_F \preceq \frac{\|\Delta \mathcal{A}\|_F}{\sigma_{\min}^2 - \mathbf{O}(\epsilon^2\ast\mathbint{1})} = \mathbf{O}(\epsilon\ast\mathbint{1})\approx \mathbf{O}(\epsilon).
    \end{displaymath}
\end{enumerate}
\end{theorem}
\begin{proof}
By hypothesis, $\mathbf{u}_{\mathrm{defect}} \circ_M \sigma_j(\mathcal{A}) = \mathbint{0}$. Under pseudoinverse $\mathcal{A}^\dagger$, $\mathbf{u}_{\mathrm{defect}} \circ_M \sigma_j^\dagger(\mathcal{A}) = \mathbint{0}$. Applying \Cref{thrm:tubal-weyl} to the perturbed tensor $\mathcal{A}_{\mathrm{noisy}} = \mathcal{A} + \Delta\mathcal{A}$, we get,
\begin{displaymath}
    |\tilde{\sigma}_j - \sigma_j| \preceq \|\Delta\mathcal{A}\|_2 \ast \mathbint{1} \preceq \epsilon \ast \mathbint{1}.
\end{displaymath}
Multiplying by the idempotent $\mathbf{u}_{\mathrm{defect}}$, we get,
\begin{displaymath}
    \mathbf{u}_{\mathrm{defect}} \circ_M \tilde{\sigma}_j = \mathbf{u}_{\mathrm{defect}} \circ_M |\tilde{\sigma}_j - \mathbint{0}| \preceq \epsilon \ast \mathbf{u}_{\mathrm{defect}}.
\end{displaymath}
For a simple perturbation, $\mathbf{u}_{\mathrm{defect}} \circ_M \tilde{\sigma}_j = \mathbf{O}(\epsilon) \ast \mathbf{u}_{\mathrm{defect}} \succ \mathbint{0}$. Naive inversion attempts a direct tubal inversion of the perturbed component as follows,
\begin{displaymath}
    \mathbf{u}_{\mathrm{defect}} \circ_M \tilde{\sigma}_j^\dagger = \frac{1}{\mathbf{O}(\epsilon)} \ast \mathbf{u}_{\mathrm{defect}} = \mathbf{O}(\epsilon^{-1}) \ast \mathbf{u}_{\mathrm{defect}}.
\end{displaymath}
Projecting the difference tensor onto $\vec{\mathcal{V}}_j$ and $\vec{\mathcal{U}}_j^\top$, we get,
\begin{align*}
    \|\mathcal{A}^\dagger - \mathcal{A}_{\mathrm{naive}}^\dagger\|_F 
    &\succeq \|\vec{\mathcal{V}}_j^\top \star_M (\mathcal{A}^\dagger - \mathcal{A}_{\mathrm{naive}}^\dagger) \star_M \vec{\mathcal{U}}_j\|_F \succeq \|\mathbf{u}_{\mathrm{defect}} \circ_M (\sigma_j^\dagger - \tilde{\sigma}_j^\dagger)\|_F \nonumber \\
    &= \left\| \mathbint{0} - \frac{1}{\mathbf{O}(\epsilon)} \ast \mathbf{u}_{\mathrm{defect}} \right\|_F = \mathbf{O}(\epsilon^{-1}).
\end{align*}
which completes the required proof for the pseudoinverse bound under naive unfolding.

Now, by \Cref{thrm:general-rank-nullity-thrm}, $\mathbb{K}_M^k = \mathscr{R}(\mathcal{A}) \oplus \mathscr{N}(\mathcal{A})$, with orthogonal projectors $\mathcal{P}_{\mathscr{R}} = \mathcal{A} \star_M \mathcal{A}^\dagger$ and $\mathcal{P}_{\mathscr{N}} = \mathcal{I} - \mathcal{A}^\dagger \star_M \mathcal{A}$. Due to the orthogonality of the submodules,
\begin{equation}\label{eq:error-decomposition}
    \|\mathcal{A}^\dagger - \mathcal{A}_{\mathrm{alg}}^\dagger\|_F^2 = \|\mathcal{P}_{\mathscr{R}} \star_M (\mathcal{A}^\dagger - \mathcal{A}_{\mathrm{alg}}^\dagger)\|_F^2 + \|\mathcal{P}_{\mathscr{N}} \star_M (\mathcal{A}^\dagger - \mathcal{A}_{\mathrm{alg}}^\dagger)\|_F^2.
\end{equation}
The proposed algebraic pseudoinverse uses idempotents $\mathbf{e}_i$ defined by the noise threshold $\epsilon < \tau < \sigma_{\min}$ as follows,
\begin{displaymath}
    \mathcal{A}_{\mathrm{alg}}^\dagger = \sum_{i=1}^k (\mathbf{e}_i \circ_M \tilde{\sigma}_i^\dagger) \star_M \vec{\mathcal{V}}_i \star_M \vec{\mathcal{U}}_i^\top,
\end{displaymath}
where $\mathbf{e}_i \circ_M \mathbf{u} = \mathbf{u}$ when $\tilde{\sigma}_i \succeq \tau \ast \mathbint{1}$, and $\mathbint{0}$ otherwise. For any component belonging to $\mathscr{N}(\mathcal{A})$, the true singular value satisfies $\sigma = \mathbint{0}$. Applying \Cref{thrm:tubal-weyl},
\begin{displaymath}
    \tilde{\sigma} \preceq \sigma + \|\Delta\mathcal{A}\|_2 \ast \mathbint{1} \preceq \epsilon \ast \mathbint{1} \prec \tau \ast \mathbint{1}.
\end{displaymath}
Hence, $\mathbf{e}_i \circ_M \mathcal{P}_{\mathscr{N}} = \mathcal{O}$, which enforces $\mathcal{P}_{\mathscr{N}} \star_M \mathcal{A}_{\mathrm{alg}}^\dagger = \mathcal{O}$. Because $\mathcal{P}_{\mathscr{N}} \star_M \mathcal{A}^\dagger = \mathcal{O}$ by \Cref{def:tubal-null-space}, the error on $\mathscr{N}(\mathcal{A})$ vanishes as follows,
\begin{equation}\label{eq:null-annihilation}
    \|\mathcal{P}_{\mathscr{N}} \star_M (\mathcal{A}^\dagger - \mathcal{A}_{\mathrm{alg}}^\dagger)\|_F = \|\mathcal{O} - \mathcal{O}\|_F = \mathbint{0}.
\end{equation}
For any component in the regular submodule $\mathscr{R}(\mathcal{A})$, $\sigma_i \succeq \sigma_{\min} \ast \mathbint{1}$. Applying \Cref{thrm:tubal-weyl}, we get,
\begin{displaymath}
    \tilde{\sigma}_i \succeq \sigma_i - \|\Delta\mathcal{A}\|_2 \ast \mathbint{1} \succeq (\sigma_{\min} - \epsilon) \ast \mathbint{1} \succ \tau \ast \mathbint{1}.
\end{displaymath}
Therefore, $\mathbf{e}_i \circ_M \mathcal{P}_{\mathscr{R}} = \mathcal{P}_{\mathscr{R}}$, leaving all components on $\mathscr{R}(\mathcal{A})$ unperturbed by the truncation.  Restricting $\mathcal{A}$ to $\mathscr{R}(\mathcal{A})$, let $\mathcal{A}_{\mathscr{R}} = \mathcal{A} \star_M \mathcal{P}_{\mathscr{R}}$. Then $\|\mathcal{A}_{\mathscr{R}}^\dagger\|_2 = \sigma_{\min}^{-1}$. Applying \Cref{lem:tubal-inversion-bound} on $\mathscr{R}(\mathcal{A})$, we get,
\begin{equation}\label{eq:range-bound}
    \|\mathcal{P}_{\mathscr{R}} \star_M (\mathcal{A}^\dagger - \mathcal{A}_{\mathrm{alg}}^\dagger)\|_F \preceq \frac{\|\mathcal{A}_{\mathscr{R}}^\dagger\|_2^2 \circ_M \|\Delta\mathcal{A}\|_F}{\mathbint{1} - \|\mathcal{A}_{\mathscr{R}}^\dagger\|_2\circ_M \|\Delta\mathcal{A}\|_2} = \frac{\sigma_{\min}^{-2}\circ_M \|\Delta\mathcal{A}\|_F}{\mathbint{1} - \epsilon\ast\sigma_{\min}^{-1}} = \frac{\|\Delta\mathcal{A}\|_F}{\sigma_{\min}^2 - \epsilon\ast\sigma_{\min}}.
\end{equation}
Substituting \Cref{eq:null-annihilation} and \Cref{eq:range-bound} in \Cref{eq:error-decomposition}, we get,
\begin{displaymath}
    \|\mathcal{A}^\dagger - \mathcal{A}_{\mathrm{alg}}^\dagger\|_F \preceq \frac{\|\Delta\mathcal{A}\|_F}{\sigma_{\min}^2 - \mathbf{O}(\epsilon^2)} = \mathbf{O}(\epsilon),
\end{displaymath}
which completes the required proof for the tubal algebraic pseudoinverse bound.
\end{proof}

\subsubsection{Algebraic Tubal Moore-Penrose Inverse}\label{subsubsec:mp-inv-experiment}
\begin{figure}[htbp]
  \subfloat[Noise Variance Analysis.]{
      \includegraphics[width=0.45\textwidth]{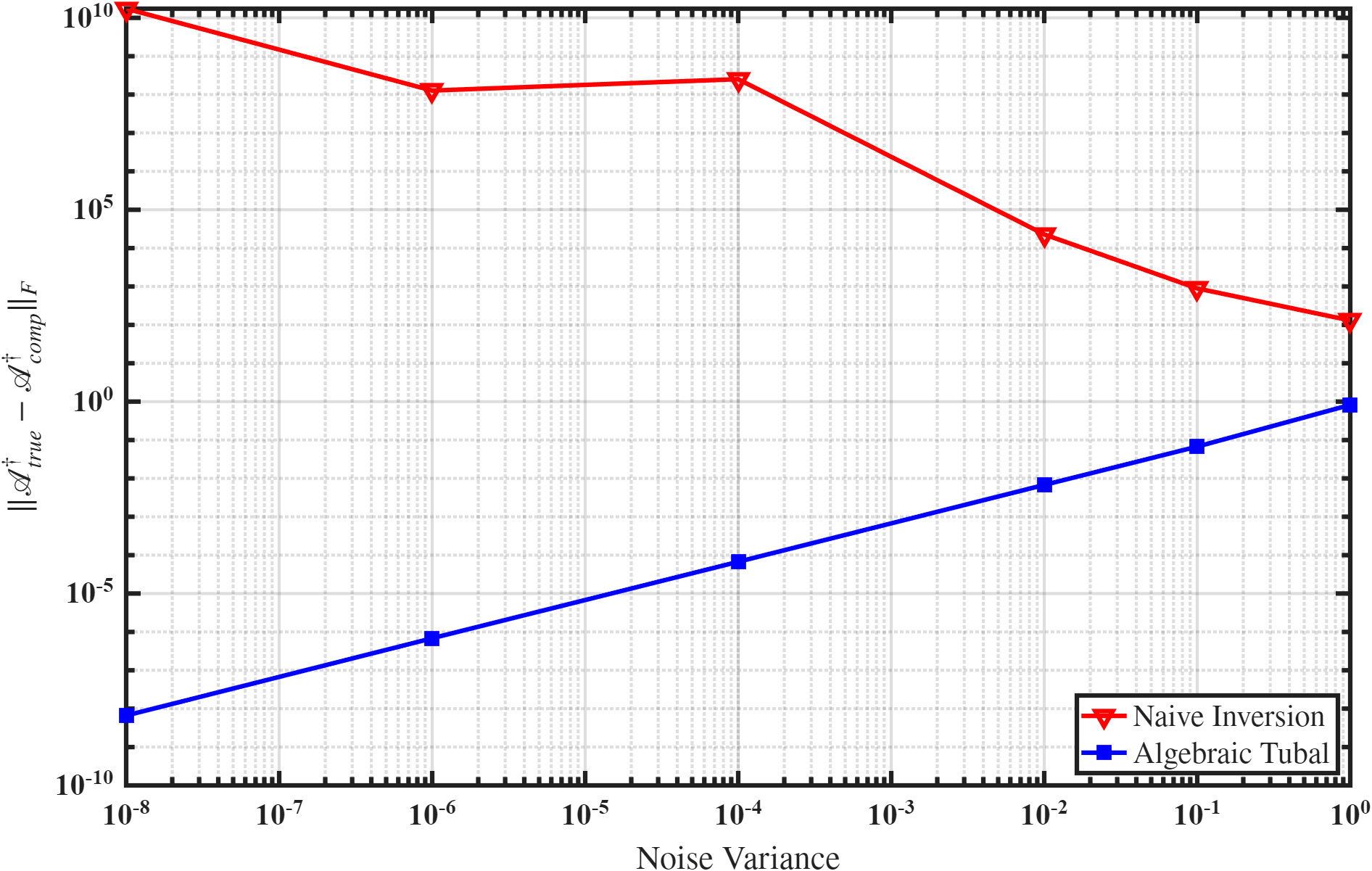}
  \label{fig:pinvNoise}
  }
  \hfill
  \subfloat[Condition Number Analysis.]{
      \includegraphics[width=0.45\textwidth]{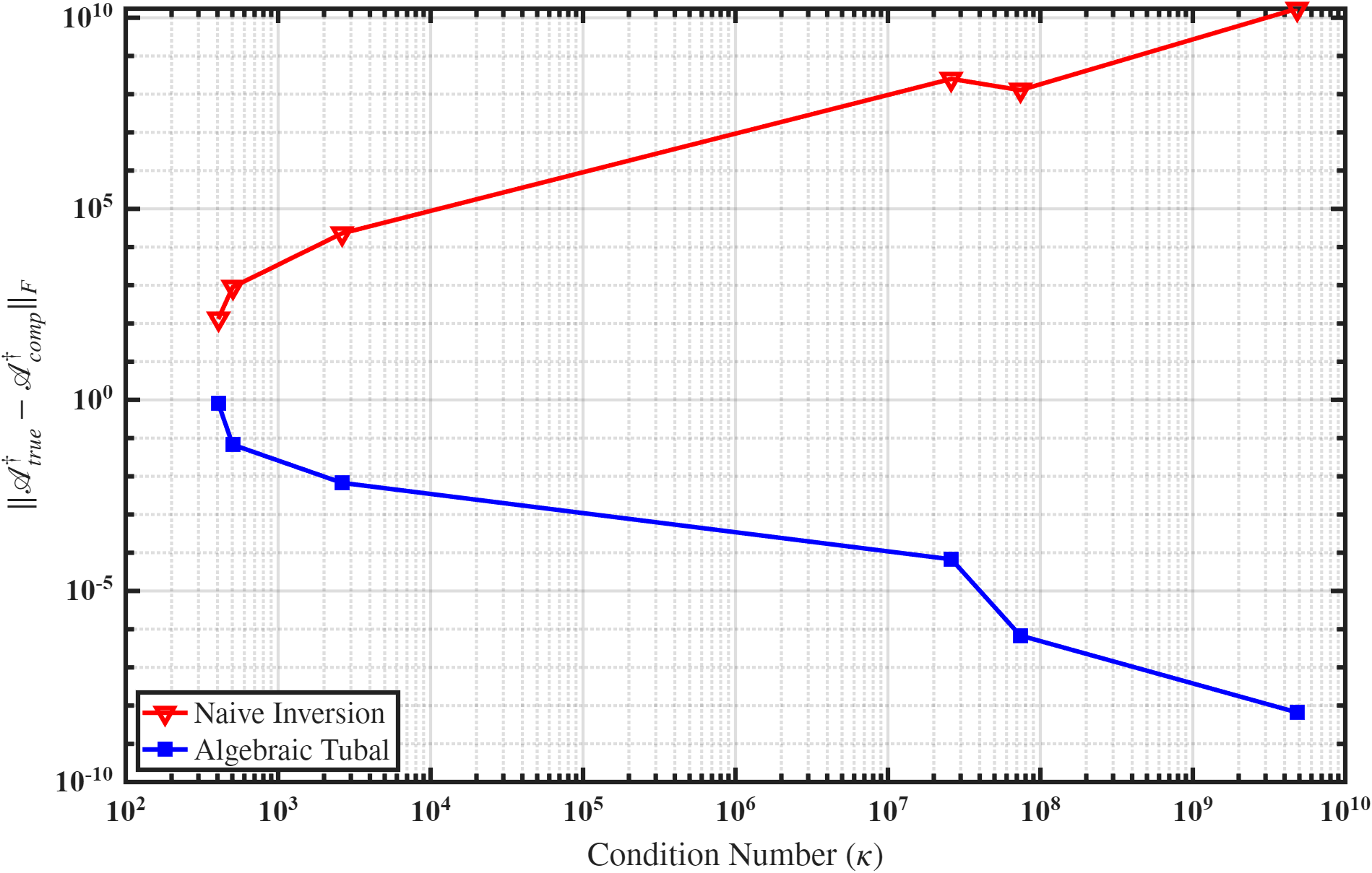}
  \label{fig:pinvCond}
  }
  \caption{Visual Analysis of the Pseudoinverse Stability along varying variances and condition numbers given in \Cref{subsubsec:mp-inv-experiment}.}
  \label{fig:pinv}
\end{figure}

Assuming a tensor of dimension $\mathcal{A}\in\mathbb{R}^{5 \times 5 \times 16}$ along with different noise variance levels $\epsilon\in\{10^{-8}, 10^{-6}, 10^{-4},$ $10^{-2}, 10^{-1}, 1\}$. We have intentionally injected a McCoy defect by replacing some of the lowest singular values by $0$ differently in slices. Also, we set slices $\geq 3$ as zero slices. We have assumed the algebraic truncation threshold to be $\tau \leq 10 \epsilon$. The analytical bound given by \Cref{thrm:perturbation-bounds} is confirmed experimentally and given in \Cref{tab:pinv_stability}. The algebraic tubal projection tracks actual perturbation in $\mathbf{O}(\epsilon)$, staying safely below the noise floor, whereas the naive approach suffers a catastrophic error of $\mathbf{O}(\epsilon^{-1})$ at small noise scales, failing to preserve the foundational topology of the system. The same can be visualized through \Cref{fig:pinv}, and the computational framework is given by \Cref{alg:tubal_pinv}.

\begin{algorithm}[htbp]
\caption{Computational Algebraic Tubal Pseudoinverse}
\label{alg:tubal_pinv}
\begin{algorithmic}[1]
\State \textbf{Require:} $\mathcal{A} \in \mathbb{R}^{n_1 \times n_2 \times \dots \times n_d}$, Threshold $\tau$
\State Compute $N = n_3n_4\dots n_d$
\State Compute $\hat{\mathcal{A}} = \Phi(\mathcal{A})$
\For{$j = 1, 2, \dots, N$}
    \State $[\mathbf{U}, \mathbf{S}, \mathbf{V}] = \texttt{svd}(\hat{\mathcal{A}}^{[j]})$
    \State Initialize $\mathbf{S}^\dagger = O\in\mathbb{R}^{n_2 \times n_1}$
    \For{$i = 1, 2, \dots, \min(n_1, n_2)$}
        \If{$\mathbf{S}_{ii} > \tau$}
            \State $\mathbf{S}^\dagger_{ii} = \frac{1}{\mathbf{S}_{ii}}$
        \EndIf
    \EndFor
    \State $\hat{\mathcal{A}}^{\dagger [j]} = \mathbf{V} \mathbf{S}^\dagger \mathbf{U}^T$
\EndFor
\State Compute $\mathcal{A}^\dagger = \Phi^{-1}(\hat{\mathcal{A}}^\dagger)$
\end{algorithmic}
\end{algorithm}

\begin{table}[htbp]
\centering
\begin{tabular}{ccc}
\hline
Noise Var & Naive Inversion Error & Algebraic Tubal Error \\
\hline
$10^{-8}$ & $1.73 \times 10^{10} \pm 1.72 \times 10^{10}$ & $6.62 \times 10^{-9} \pm 8.77 \times 10^{-10}$ \\
$10^{-6}$ & $1.27 \times 10^{8} \pm 2.41 \times 10^{8}$ & $6.83 \times 10^{-7} \pm 9.11 \times 10^{-8}$ \\
$10^{-4}$ & $2.52 \times 10^{8} \pm 7.55 \times 10^{8}$ & $6.71 \times 10^{-5} \pm 7.24 \times 10^{-6}$ \\
$10^{-2}$ & $2.28 \times 10^{4} \pm 4.59 \times 10^{4}$ & $6.75 \times 10^{-3} \pm 8.39 \times 10^{-4}$ \\
$10^{-1}$ & $8.93 \times 10^{2} \pm 8.05 \times 10^{2}$ & $6.73 \times 10^{-2} \pm 1.65 \times 10^{-2}$ \\
$10^{0}$ & $1.30 \times 10^{2} \pm 7.98 \times 10^{1}$ & $8.23 \times 10^{-1} \pm 1.11 \times 10^{-16}$ \\
\hline
\end{tabular}
\caption{Pseudoinverse Stability Analysis in the presence of noise through experiment given in \Cref{subsubsec:mp-inv-experiment}.}
\label{tab:pinv_stability}
\end{table}

\subsubsection{Tensor Randomized Block Kaczmarz Algorithm}\label{subsubsec:kaczmarz}
Tensor systems play an important role in various applications across science and engineering. Kaczmarz algorithms~\cite{kaczmarz1, kaczmarz2} are one of the prominent iterative solvers for linear systems of equations that project the solution iteratively onto the range space. Assuming an overdetermined and noisy linear tensor system is partitioned into contiguous row blocks $\tau \in \{1, \dots, N\}$ as follows,
\begin{equation}
    \mathcal{A}_\tau \star_M \vec{\mathcal{X}} = \vec{\mathcal{B}}_{\tau, \text{clean}} + \vec{\mathcal{E}}_{\tau, \text{noise}}
\end{equation}
The randomized tensor block kaczmarz algorithm~\cite{kaczmarz2} iteratively projects the solution onto the affine subspace of block $\tau$ as follows,
\begin{displaymath}
    \vec{\mathcal{X}}_{k+1} = \vec{\mathcal{X}}_k - \mathcal{A}_\tau^\dagger \star_M (\mathcal{A}_\tau \star_M \vec{\mathcal{X}}_k - \vec{\mathcal{B}}_\tau)
\end{displaymath}
The naive tubal pseudoinverse computes the Moore-Penrose inverse for all slices independently. The solution update incorporates the noise term regardless of whether it is a McCoy defect, $\Delta \vec{\mathcal{X}}_{\mathrm{noise}} = \mathcal{A}_\tau^\dagger \star_M \vec{\mathcal{E}}_{\tau, \mathrm{noise}}$. If slice $j$ contains the McCoy defect $\epsilon$, the naive pseudoinverse operator yields $\|\hat{\mathcal{A}}_\tau^{\dagger [j]}\|_2 = \mathbf{O}(\epsilon^{-1})$.

\begin{algorithm}[htbp]
\caption{Tensor Randomized Block Kaczmarz Algorithm}
\label{alg:tensor_kaczmarz}
\begin{algorithmic}[1]
\State \textbf{Require:} $\mathcal{A} \in \mathbb{R}^{m \times k \times n_3\times \dots\times n_d}$, $\vec{\mathcal{B}} \in \mathbb{R}^{m \times 1 \times n_3\times \dots\times n_d}$, Block size $b$, $\texttt{iters}$, $\tau$
\State Initialize $\vec{\mathcal{X}} = \vec{\mathcal{O}} \in \mathbb{R}^{k \times 1 \times n_3\times \dots\times n_d}$
\State Partition rows of $\mathcal{A}$ and $\vec{\mathcal{B}}$ into $N = \lceil m/b \rceil$ blocks
\For{\texttt{it} $= 1, 2, \dots, \texttt{iters}$}
    \State Permute block ordering $P = \text{randperm}(N)$
    \For{\textbf{each} index $i \in P$}
        \State $\vec{\mathcal{R}} = \mathcal{A}_i \star_M \vec{\mathcal{X}} - \vec{\mathcal{B}}_i$
        \State $\mathcal{A}_i^\dagger = \texttt{pinv}(\mathcal{A}_i, \mathbf{M}, \tau)$
        \State $\vec{\mathcal{X}} \leftarrow \vec{\mathcal{X}} - \mathcal{A}_i^\dagger \star_M \vec{\mathcal{R}}$
    \EndFor
\EndFor
\end{algorithmic}
\end{algorithm}

\begin{figure}[htbp]
  \subfloat[Kacmarz Convergence Analysis.]{
      \includegraphics[width=0.4\textwidth]{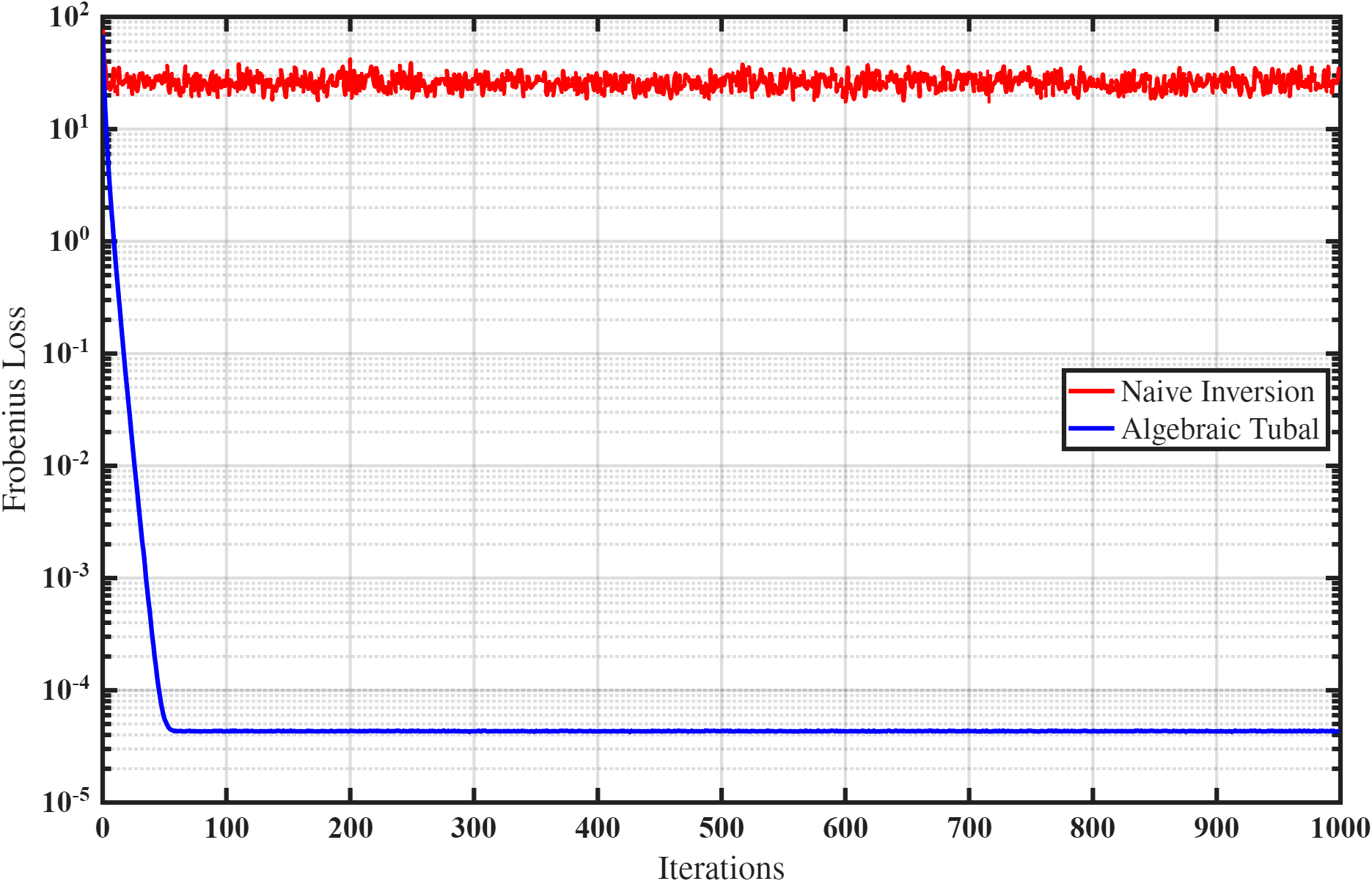}
  }
  \hfill
  \subfloat[Kaczmarz Stability Analysis.]{
      \includegraphics[width=0.4\textwidth]{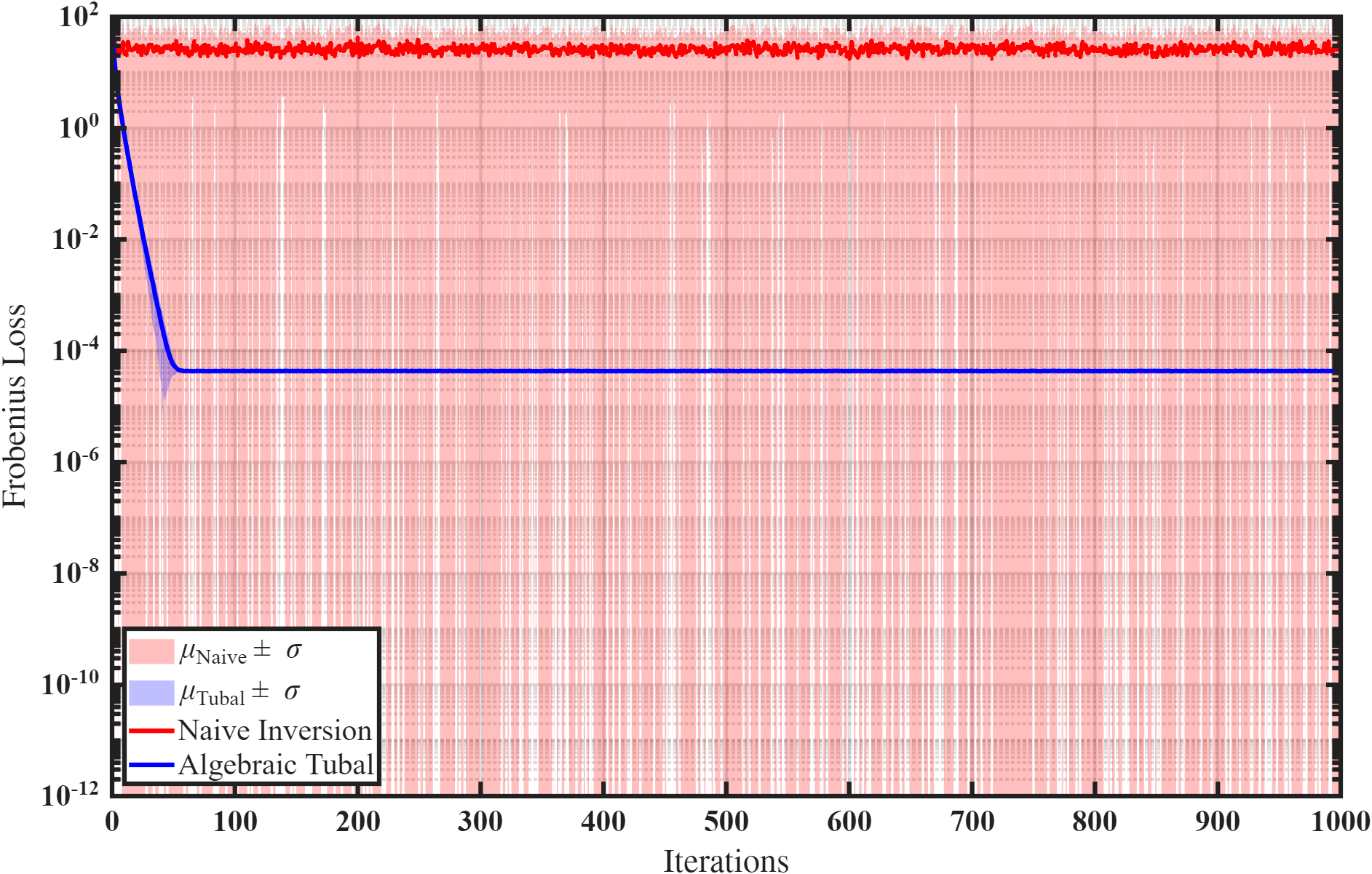}
  }
  \caption{Kaczmarz convergence and stability analysis of the iterations for the experiment given in \Cref{subsubsec:kaczmarz}.}
  \label{fig:kaczmarz}
\end{figure}

Assuming tensor $\mathcal{A}\in\mathbb{R}^{100\times 50\times 16}$ with block size, $b = 10$ and maximum number of iterations as $1000$. We have introduced $\vec{\mathcal{E}}_{\tau, \text{noise}} = \mathbf{O}(10^{-6})$ and McCoy defect is injected similar to \Cref{subsubsec:mp-inv-experiment} with noise $\epsilon = \mathbf{O}(10^{-7})$. If we are using naive inversion in \Cref{alg:tensor_kaczmarz}, we get the final loss as $3.57 \times 10^{1} \pm 4.69 \times 10^{1}$. However, if we utilize \Cref{alg:tubal_pinv} for finding \texttt{pinv} in \Cref{alg:tensor_kaczmarz}, we get the final loss as $4.30 \times 10^{-5} \pm 9.85 \times 10^{-7}$. It is also clear from \Cref{fig:kaczmarz} that the naive inversion leads to very high loss and large variance across runs. However, utilizing \Cref{alg:tubal_pinv} leads to convergence and much less variance across multiple runs.

\section{Conclusion}\label{sec:conclusions}
In this paper, we have formalized a rigorous algebra for $d$th-order matrix-mimicking tensors. We introduced a unique tubal determinant function and utilized the tubal McCoy rank to provide exact algebraic conditions for tensor invertibility, completely characterizing the defects caused by zero-divisor tubes as well. Furthermore, we resolved the ambiguities of linear dependence over rings by distinguishing between absolute and tubal dependence, which enabled the proof of a tubal rank-nullity theorem using orthogonal projectors onto projective submodules. In the spectral domain, we formalized regular tubal eigenvectors and proved the Gershgorin hypersphere theorem, offering a crucial analytical tool for bounding tensor eigenvalues. We have proved the algebraic applicability of the Gershgorin theorem for finding the learning rate of the gradient descent algorithm and the tubal rank-nullity theorem for finding the pseudoinverse in the Kazmarz algorithm.

The algebraic properties established in this work provide the necessary theoretical foundations to analyze the stability, convergence, and exactness of tensor-based algorithms in various applications. Future work will extend this framework to complex settings ($\mathbb{C}$) as well as non-commutative algebraic structures, such as quaternions ($\mathbb{H}$), to handle highly multidimensional and correlated data. Additionally, we plan to explore tubal tensor calculus, generalizing derivatives and manifolds over $\mathbb{K}_M$ for tensor optimization techniques.


\section*{Declarations}
\subsection*{Ethical approval}
Not applicable.
\subsection*{Availability of supporting data} 
Not applicable.
\subsection*{Competing interests}
The authors declare that they have no conflicts of interest.
\subsection*{Funding}
No research funding was received.


\bibliographystyle{elsarticle-num}
\bibliography{references}

\end{document}